\documentclass[journal]{IEEEtran}

\usepackage[utf8]{inputenc}
\usepackage{color}
\usepackage{array}
\usepackage{verbatim}
\usepackage{float}
\usepackage{amsmath}
\usepackage{amsthm}
\usepackage{amssymb}
\usepackage{graphicx}
\usepackage{longtable}
\usepackage{multirow}
\usepackage[unicode=true,
bookmarks=false,
breaklinks=false,pdfborder={0 0 1},colorlinks=false]
{hyperref}
\hypersetup{
	colorlinks,bookmarksopen,bookmarksnumbered,citecolor=blue,urlcolor=blue}
\usepackage{cite}

\makeatletter
\let\oldforeign@language\foreign@language
\DeclareRobustCommand{\foreign@language}[1]{%
	\lowercase{\oldforeign@language{#1}}}

\let\oldforeign@language\foreign@language
\DeclareRobustCommand{\foreign@language}[1]{%
	\lowercase{\oldforeign@language{#1}}}

\ifCLASSINFOpdf
\else
\fi

\newcommand{\MYfooter}{\smash{
		\hfil\parbox[t][\height][t]{\textwidth}{\centering
			\thepage}\hfil\hbox{}}}

\makeatletter

\def\ps@IEEEtitlepagestyle{%
	\def\@oddhead{\parbox[t][\height][t]{\textwidth}{\centering \scriptsize
			Personal use of this material is permitted. Permission from the author(s) and/or copyright holder(s), must be obtained for all other uses. Please contact us and provide details if you believe this document breaches copyrights.\\
			\noindent\makebox[\linewidth]{}
		}\hfil\hbox{}}%
	\def\@evenhead{\scriptsize\thepage \hfil \leftmark\mbox{}}%
	\def\@oddfoot{\parbox[t][\height][l]{\textwidth}{
			\vspace{-20pt}{\rule{\textwidth}{0.4pt}}\\ \footnotesize\underline{To cite this article:}
			{\bf{\textcolor{red}{H. A. Hashim, L. J. Brown, and K. McIsaac, "Nonlinear Pose Filters on the Special Euclidean Group SE(3) with Guaranteed Transient and Steady-state Performance," IEEE Transactions on Systems, Man, and Cybernetics: Systems, vol. 51, no. 5, pp. 2949-2962, 2021.}}} doi: \href{https://doi.org/10.1109/TSMC.2019.2920114}{10.1109/TSMC.2019.2920114}\\
			\noindent\makebox[\linewidth]
		}\hfil\hbox{}}%
	\def\@evenfoot{\MYfooter}}

\makeatother
\newtheorem{lem}{Lemma}

\newtheorem{prop}{Proposition}
\newtheorem{thm}{Theorem}
\newtheorem{rem}{Remark}

\newtheorem{assum}{Assumption}

\begin{document}
	\bstctlcite{IEEEexample:BSTcontrol}

\title{Nonlinear Pose Filters on the Special Euclidean Group SE(3) with Guaranteed Transient and Steady-state Performance}

\author{Hashim~A.~Hashim, Lyndon J. Brown, and~Kenneth McIsaac
\thanks{This work was supported in part by the NSERC Discovery Grant program.}
\thanks{H. A. Hashim, L. J. Brown and K. McIsaac are with the Department of Electrical and Computer Engineering,
University of Western Ontario, London, ON, Canada, N6A-5B9, e-mail: hmoham33@uwo.ca, lbrown@uwo.ca and kmcisaac@uwo.ca.}
}



\maketitle

\begin{abstract}
Two novel nonlinear pose (\textit{i.e}, attitude and position) filters developed directly on the Special
Euclidean Group $\mathbb{SE}\left(3\right)$ able to guarantee prescribed
characteristics of transient and steady-state performance are proposed.
The position error and normalized Euclidean distance of attitude error
are trapped to arbitrarily start within a given large set and converge
systematically and asymptotically to the origin from almost any initial
condition. The transient error is guaranteed not to exceed a prescribed
value while the steady-state error is bounded by a predefined small
value. %
{} The first pose filter operates based on a set of vectorial measurements
coupled with a group of velocity vectors and requires preliminary
pose reconstruction. The second filter, on the contrary, is able to
perform its function using a set of vectorial measurements and a group
of velocity vectors directly. Both proposed filters provide reasonable
pose estimates with superior convergence properties while being able
to use measurements obtained from low-cost inertial measurement, landmark
measurement, and velocity measurement units. The equivalent quaternion representation and complete implementation
steps of the proposed filters are presented.  Simulation results demonstrate
effectiveness and robustness of the proposed filters considering large
error in initialization and high level of uncertainties in velocity
vectors as well as in the set of vector measurements.
\end{abstract}

\begin{IEEEkeywords}
Attitude, position, pose estimation, nonlinear observer, special orthogonal
group, special Euclidean group, SO(3), SE(3), prescribed performance
function, transient, steady-state error, transformed error, feature
measurement, PPF, IMU.
\end{IEEEkeywords}

\IEEEpeerreviewmaketitle{}

\section{Introduction}

\IEEEPARstart{R}{obotics} and engineering applications such as aerial and underwater
vehicles, satellites and space crafts are concerned with accurately
estimating the pose of a rigid-body in 3D space. In essence, the pose
of a rigid-body consists of two elements: orientation and position.
The orientation of a rigid-body in 3D space is often referred to as
attitude, therefore, orientation and attitude will be used interchangeably.
One of the basic methods of attitude reconstruction is the algebraic
approach. It allows to reconstruct the attitude given the availability
of two or more non-collinear inertial-frame vectors and their body-frame
vectors utilizing algorithms such as QUEST \cite{shuster1981three},
or singular value decomposition (SVD) \cite{markley1988attitude}.
However, the process of attitude reconstruction is vulnerable to the
effects of noise and bias contaminating the body-frame measurements
which causes \cite{shuster1981three,markley1988attitude} to produce
unsatisfactory results. This is particularly true in the context of
a rigid-body fitted with low-cost inertial measurement unit (IMU)
\cite{hashim2018SO3Stochastic,mahony2008nonlinear,mohamed2019filter}.

Gaussian filters or nonlinear deterministic filters have been used
historically to address the challenge of attitude estimation \cite{hashim2018SO3Stochastic}.
The family of Gaussian filters, which includes Kalman filter (KF)
\cite{choukroun2006novel}, extended KF (EKF) \cite{lefferts1982kalman},
and multiplicative EKF (MEKF) \cite{markley2003attitude}, often consider
the unit quaternion in attitude representation \cite{hashim2018SO3Stochastic,mohamed2019filter}.
For good survey of Gaussian filters visit \cite{hashim2018SO3Stochastic}.
However, it is crucial to note the nonlinear nature of the attitude
problem. Nonlinear attitude filters such as \cite{hashim2018SO3Stochastic,mahony2008nonlinear,hashim2018Conf1,grip2012attitude,liu2018complementary}
are evolved directly on the Special Orthogonal Group $\mathbb{SO}\left(3\right)$.
In particular, nonlinear deterministic attitude filters outperform
the Gaussian filters in many respects, namely they are simpler in
derivation and representation, they demand less processing power,
and they show better tracking convergence \cite{hashim2018SO3Stochastic,mahony2008nonlinear}.
Attitude estimation is an essential part of the pose estimation problem.
Taking into consideration the remarkable advantages of nonlinear attitude
filters, attitude-position (pose) filtering problem is best approached
in a nonlinear sense. 

The pose estimation problem relies on filters evolved on the Special
Euclidean Group $\mathbb{SE}\left(3\right)$ which require a measurement
derived from a group velocity vector, vectorial measurements that
could be provided by IMU, landmark measurements collected, for example,
by a vision system and an estimate of the bias associated with velocity
measurements.%
{} Pose estimation commonly involves a computer vision system with a
monocular camera and IMU \cite{rehbinder2003pose,baldwin2007complementary,hashim2018SE3Stochastic,hashim2019Conf1}.
The pose filter described in \cite{baldwin2007complementary} was
developed directly on $\mathbb{SE}\left(3\right)$ and its performance
has been proven to be exponentially stable. Although, the filter in
\cite{baldwin2007complementary} requires pose reconstruction for
the implementation, the nonlinear filter can be modified to function
based solely on a set of vectorial measurements avoiding the need
for pose reconstruction \cite{baldwin2009nonlinear,hua2011observer}.
In spite of the simplicity of the filter design in \cite{baldwin2007complementary,baldwin2009nonlinear,hua2011observer},
numerical results show high sensitivity to noise and bias attached
to the measurements. In addition, no systematic convergence is observed
in \cite{rehbinder2003pose,baldwin2007complementary,baldwin2009nonlinear,hua2011observer,vasconcelos2010nonlinear,dominguez2017simultaneous,hua2017riccati,hashim2018SE3Stochastic},
such that the tracking error does not follow a predefined transient
and steady-state measures. Accordingly, successful pose estimation
for spacecraft control applications, such as \cite{tanaka2017practical,santoso2018robust,sun2018disturbance,lee2018geometric},
cannot be achieved without pose filters which are robust against uncertain
measurements, demonstrate fast tracking performance, and satisfy a
certain level of transient and steady-state characteristics.

Prescribed performance implies confining the error to initially start
within a predefined large set and decay systematically and smoothly
to a predefined small residual set \cite{bechlioulis2008robust}.
The error trajectory is constrained by a prescribed performance function
(PPF) to satisfy transient as well as steady-state performance. The
main objective of prescribed performance is to relax the constrained
error to its unconstrained form, termed transformed error, which allows
to keep the error within the decaying dynamic boundaries, and thereby
achieve successful estimation or control applications. These applications
include but are not limited to two degrees of freedom planar robots
\cite{bechlioulis2008robust,wang2017dynamic}, uncertain nonlinear
systems \cite{yang2018prescribed}, servo mechanism with friction
compensation \cite{na2014adaptive}, and uncertain multi-agent system
\cite{hashim2017neuro,hashim2017adaptive}.

In this paper two robust nonlinear pose filters on $\mathbb{SE}\left(3\right)$
with predefined transient as well as steady-state measures are proposed.
The main contributions are as follows:
\begin{enumerate}
	\item[\textbf{1)}] The proposed filters guarantee boundedness of the closed loop error
	signals with constrained error and unconstrained transformed error
	being proven to be almost globally asymptotically stable such that
	the error in the homogeneous transformation matrix is regulated asymptotically
	to the identity from almost any initial condition. Most significantly,
	the exceptional performance is guaranteed even when the measurements
	are supplied by a low-cost measurement unit, for instance, an IMU
	module equipped with a gyroscope, a vision unit, and a GPS.
	\item[\textbf{2)}] The proposed filters guarantee systematic convergence of the error
	controlled by the dynamic reducing boundaries forcing the error to
	start within a predefined large set and decrease systematically and
	smoothly to a residual small set%
	, unlike \cite{rehbinder2003pose,baldwin2007complementary,baldwin2009nonlinear,hua2011observer,vasconcelos2010nonlinear,hashim2018SE3Stochastic}.
	\item[\textbf{3)}] The proposed pose filters are more efficient at ensuring fast convergence
	compared to similar estimators described in the literature, for instance
	\cite{rehbinder2003pose,baldwin2007complementary,baldwin2009nonlinear,hua2011observer,vasconcelos2010nonlinear,hashim2018SE3Stochastic}. 
\end{enumerate}
The fast convergence is mainly attributed to the dynamic behavior
of the estimator gains. The first filter requires a group of velocity
vectors and a set of measurements to obtain an online algebraic reconstruction
of the pose. The second filter uses the group of velocity vector and
the set of vectorial measurements directly. 

The remainder of the paper is organized as follows: Section \ref{sec:SE3PPF_Math-Notations}
gives an overview of $\mathbb{SO}\left(3\right)$ and $\mathbb{SE}\left(3\right)$,
mathematical notation and identities. The pose problem is formulated,
vector measurements are demonstrated and prescribed performance is
introduced in Section \ref{sec:SE3PPF_Problem-Formulation-in}. The
two proposed filters and the related stability analysis are presented
in Section \ref{sec:SE3PPF-Filters}. Section \ref{sec:SO3PPF_Simulations}
elaborates on the effectiveness and robustness of the proposed filters.
Finally, Section \ref{sec:SO3PPF_Conclusion} draws a conclusion of
this work.

\section{Preliminaries and Mathematical Identities \label{sec:SE3PPF_Math-Notations}}

In this paper $\mathbb{R}_{+}$ refers to the set of nonnegative real
numbers. $\mathbb{R}^{n}$ and $\mathbb{R}^{n\times m}$ denote a
real $n$-dimensional space column vector and real $n\times m$ dimensional
space, respectively. The Euclidean norm of $x\in\mathbb{R}^{n}$ is
$\left\Vert x\right\Vert =\sqrt{x^{\top}x}$ with $^{\top}$ being
the transpose of the component. $\lambda\left(\cdot\right)$ denotes
a set of singular values of a matrix with $\underline{\lambda}\left(\cdot\right)$
being its minimum value. $\mathbf{I}_{n}$ stands for an $n$-by-$n$
identity matrix, while $\underline{\mathbf{0}}_{n}\in\mathbb{R}^{n}$
is a zero column vector. The frame notation is as follows: $\left\{ \mathcal{B}\right\} $
refers to the body-frame and $\left\{ \mathcal{I}\right\} $ represents
the inertial-frame.

Define $\mathbb{GL}\left(3\right)$ as a 3-dimensional general linear
group which is a Lie group with smooth multiplication and inversion.
The orthogonal group, denoted by $\mathbb{O}\left(3\right)$, is a
subgroup of $\mathbb{GL}\left(3\right)$ defined by
\[
\mathbb{O}\left(3\right)=\left\{ \left.M\in\mathbb{R}^{3\times3}\right|M^{\top}M=MM^{\top}=\mathbf{I}_{3}\right\} 
\]
with $\mathbf{I}_{3}$ being a 3-by-3 identity matrix. Let $\mathbb{SO}\left(3\right)$
denote the Special Orthogonal Group which is a subgroup of $\mathbb{O}\left(3\right)$.
The orientation of a rigid-body in 3D space is termed attitude, denoted
by $R$, and defined as follows: 
\[
\mathbb{SO}\left(3\right)=\left\{ \left.R\in\mathbb{R}^{3\times3}\right|RR^{\top}=R^{\top}R=\mathbf{I}_{3}\text{, }{\rm det}\left(R\right)=+1\right\} 
\]
with ${\rm det\left(\cdot\right)}$ being the determinant of the associated
matrix. $\mathbb{SE}\left(3\right)$ stands for the Special Euclidean
Group, a subset of the affine group $\mathbb{GA}\left(3\right)=\mathbb{SO}\left(3\right)\times\mathbb{R}^{3}$
defined by
\[
\mathbb{SE}\left(3\right)=\left\{ \left.\boldsymbol{T}\in\mathbb{R}^{4\times4}\right|R\in\mathbb{SO}\left(3\right),P\in\mathbb{R}^{3}\right\} 
\]
where $\boldsymbol{T}\in\mathbb{SE}\left(3\right)$, termed a homogeneous
transformation matrix, represents the pose of a rigid-body in 3D space
with
\begin{equation}
\boldsymbol{T}=\left[\begin{array}{cc}
R & P\\
\underline{\mathbf{0}}_{3}^{\top} & 1
\end{array}\right]\in\mathbb{SE}\left(3\right)\label{eq:SE3STCH_T_matrix}
\end{equation}
where $P\in\mathbb{R}^{3}$ and $R\in\mathbb{SO}\left(3\right)$ denote
position and attitude of a rigid-body in 3D space, respectively, and
$\underline{\mathbf{0}}_{3}^{\top}$ is a zero row. $\mathfrak{so}\left(3\right)$
is a Lie-algebra related to $\mathbb{SO}\left(3\right)$ defined by
\[
\mathfrak{so}\left(3\right)=\left\{ \left.A\in\mathbb{R}^{3\times3}\right|A^{\top}=-A\right\} 
\]
where $A$ is a skew symmetric matrix. Define the map $\left[\cdot\right]_{\times}:\mathbb{R}^{3}\rightarrow\mathfrak{so}\left(3\right)$
as 
\[
\left[\alpha\right]_{\times}=\left[\begin{array}{ccc}
0 & -\alpha_{3} & \alpha_{2}\\
\alpha_{3} & 0 & -\alpha_{1}\\
-\alpha_{2} & \alpha_{1} & 0
\end{array}\right]\in\mathfrak{so}\left(3\right),\hspace{1em}\alpha=\left[\begin{array}{c}
\alpha_{1}\\
\alpha_{2}\\
\alpha_{3}
\end{array}\right]
\]
For any $\alpha,\beta\in\mathbb{R}^{3}$, we define $\left[\alpha\right]_{\times}\beta=\alpha\times\beta$
with $\times$ being the cross product. The wedge operator is denoted
by $\wedge$, and for any $\mathcal{Y}=\left[y_{1}^{\top},y_{2}^{\top}\right]^{\top}$
with $y_{1},y_{2}\in\mathbb{R}^{3}$ the wedge map $\left[\cdot\right]_{\wedge}:\mathbb{R}^{6}\rightarrow\mathfrak{se}\left(3\right)$
is defined by
\[
\left[\mathcal{Y}\right]_{\wedge}=\left[\begin{array}{cc}
\left[y_{1}\right]_{\times} & y_{2}\\
\underline{\mathbf{0}}_{3}^{\top} & 0
\end{array}\right]\in\mathfrak{se}\left(3\right)
\]
$\mathfrak{se}\left(3\right)$ is a Lie algebra of $\mathbb{SE}\left(3\right)$
and can be expressed as{\small{}
	\begin{align*}
	\mathfrak{se}\left(3\right) & =\left\{ \left.\left[\mathcal{Y}\right]_{\wedge}\in\mathbb{R}^{4\times4}\right|\exists y_{1},y_{2}\in\mathbb{R}^{3}:\left[\mathcal{Y}\right]_{\wedge}=\left[\begin{array}{cc}
	\left[y_{1}\right]_{\times} & y_{2}\\
	\underline{0}_{3}^{\top} & 0
	\end{array}\right]\right\} 
	\end{align*}
}The inverse of $\left[\cdot\right]_{\times}$ is defined by $\mathbf{vex}:\mathfrak{so}\left(3\right)\rightarrow\mathbb{R}^{3}$,
and for $\alpha\in\mathbb{R}^{3}$ and $\left[\alpha\right]_{\times}\in\mathfrak{so}\left(3\right)$
we have
\begin{equation}
\mathbf{vex}(\left[\alpha\right]_{\times})=\alpha\in\mathbb{R}^{3}\label{eq:SE3STCH_VEX}
\end{equation}
$\boldsymbol{\mathcal{P}}_{a}$ stands for an anti-symmetric projection
operator on the Lie-algebra $\mathfrak{so}\left(3\right)$ while its
mapping is given by $\boldsymbol{\mathcal{P}}_{a}:\mathbb{R}^{3\times3}\rightarrow\mathfrak{so}\left(3\right)$
such that
\begin{equation}
\boldsymbol{\mathcal{P}}_{a}\left(M\right)=\frac{1}{2}(M-M^{\top})\in\mathfrak{so}\left(3\right),\,M\in\mathbb{R}^{3\times3}\label{eq:SE3STCH_Pa}
\end{equation}
The normalized Euclidean distance of the attitude matrix $R\in\mathbb{SO}\left(3\right)$
is given by
\begin{equation}
\left\Vert R\right\Vert _{I}=\frac{1}{4}{\rm Tr}\{\mathbf{I}_{3}-R\}\label{eq:SE3STCH_Ecul_Dist}
\end{equation}
with ${\rm Tr}\left\{ \cdot\right\} $ being a trace of a matrix,
while $\left\Vert R\right\Vert _{I}\in\left[0,1\right]$. To reconstruct
the orientation of any rigid-body in 3D space it is sufficient to
know unit-axis $u\in\mathbb{R}^{3}$ in the sphere $\mathbb{S}^{2}$
and angle of rotation $\alpha\in\mathbb{R}$ about $u$. This type
of parameterization is termed angle-axis parameterization and its
mapping to $\mathbb{SO}\left(3\right)$ is given by $\mathcal{R}_{\alpha}:\mathbb{R}\times\mathbb{R}^{3}\rightarrow\mathbb{SO}\left(3\right)$
such that
\begin{equation}
\mathcal{R}_{\alpha}(\alpha,u)=\mathbf{I}_{3}+\sin\left(\alpha\right)\left[u\right]_{\times}+\left(1-\cos\left(\alpha\right)\right)\left[u\right]_{\times}^{2}\label{eq:SE3STCH_att_ang}
\end{equation}
For $\alpha,\beta\in{\rm \mathbb{R}}^{3}$, $R\in\mathbb{SO}\left(3\right)$,
$A\in\mathbb{R}^{3\times3}$ and $B=B^{\top}\in\mathbb{R}^{3\times3}$
the following mathematical identities
\begin{align}
\left[\alpha\times\beta\right]_{\times}= & \beta\alpha^{\top}-\alpha\beta^{\top}\label{eq:SO3PPF_Identity1}\\
\left[R\alpha\right]_{\times}= & R\left[\alpha\right]_{\times}R^{\top}\label{eq:SO3PPF_Identity2}\\
\left[\alpha\right]_{\times}^{2}= & -\alpha^{\top}\alpha\mathbf{I}_{3}+\alpha\alpha^{\top}\label{eq:SO3PPF_Identity3}\\
B\left[\alpha\right]_{\times}+\left[\alpha\right]_{\times}B= & {\rm Tr}\left\{ B\right\} \left[\alpha\right]_{\times}-\left[B\alpha\right]_{\times}\label{eq:SO3PPF_Identity4}\\
{\rm Tr}\{B\left[\alpha\right]_{\times}\}= & 0\label{eq:SO3PPF_Identity6}
\end{align}
{\small{}
	\begin{align}
	{\rm Tr}\{A\left[\alpha\right]_{\times}\}= & {\rm Tr}\{\boldsymbol{\mathcal{P}}_{a}\left(A\right)\left[\alpha\right]_{\times}\}=-2\mathbf{vex}(\boldsymbol{\mathcal{P}}_{a}(A))^{\top}\alpha\label{eq:SO3PPF_Identity7}
	\end{align}
}will be used in the subsequent derivations.

\section{Problem Formulation with Prescribed Performance \label{sec:SE3PPF_Problem-Formulation-in}}

Pose estimator relies on a set of vectorial measurements made on inertial-frame
and body-frame. This section aims to define the pose problem and present
the associated measurements. Next, the pose error and its reformulation
are geared towards attaining desired characteristics of transient
and steady-state performance. 

\subsection{Pose Kinematics and Measurements\label{subsec:SE3PPF_Pose-Kinematics}}

\begin{figure}[h]
	\centering{}\includegraphics[scale=0.43]{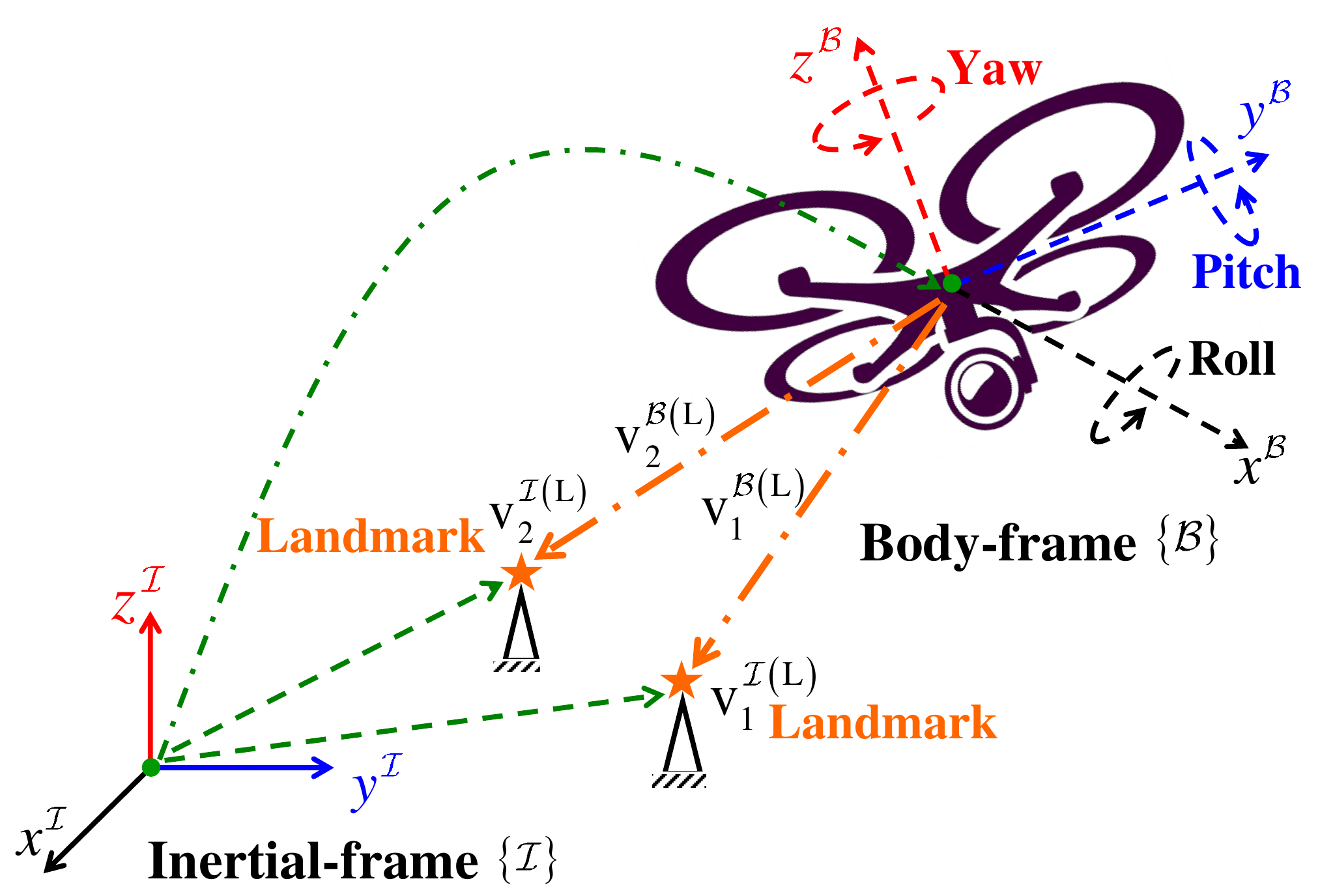}\caption{Pose estimation problem of a rigid-body in 3D space.}
	\label{fig:SE3PPF_1} 
\end{figure}

The pose of any rigid-body in 3D space consists of two elements: attitude
and position, and this work aims to estimate both elements. The attitude
of a rigid-body is commonly represented by a rotational matrix $R\in\mathbb{SO}\left(3\right)$
defined relative to the body-frame such that $R\in\left\{ \mathcal{B}\right\} $.
Position of a rigid-body is, on the contrary, defined by $P\in\mathbb{R}^{3}$
with respect to the inertial-frame $P\in\left\{ \mathcal{I}\right\} $.
The pose problem can be characterized by the homogeneous transformation
matrix $\boldsymbol{T}\in\mathbb{SE}\left(3\right)$ as
\begin{equation}
\boldsymbol{T}=\left[\begin{array}{cc}
R & P\\
\underline{\mathbf{0}}_{3}^{\top} & 1
\end{array}\right]\label{eq:SE3PPF_T_matrix2}
\end{equation}
The pose estimation problem of a rigid-body in 3D space is depicted
in Fig. \ref{fig:SE3PPF_1}.

Let the components associated with body-frame and inertial-frame be
assigned superscripts $\mathcal{B}$ and $\mathcal{I}$, respectively.
The attitude can be obtained given $N_{{\rm R}}$ known non-collinear
inertial vectors, available for measurements at a coordinate fixed
to the moving body. IMU exemplify sensors, which could provide those
measurements. The $i$th body-frame vector measurement is given by
\[
\left[\begin{array}{c}
{\rm v}_{i}^{\mathcal{B}\left({\rm R}\right)}\\
0
\end{array}\right]=\boldsymbol{T}^{-1}\left[\begin{array}{c}
{\rm v}_{i}^{\mathcal{I}\left({\rm R}\right)}\\
0
\end{array}\right]+\left[\begin{array}{c}
{\rm b}_{i}^{\mathcal{B}\left({\rm R}\right)}\\
0
\end{array}\right]+\left[\begin{array}{c}
\omega_{i}^{\mathcal{B}\left({\rm R}\right)}\\
0
\end{array}\right]
\]
such that
\begin{equation}
{\rm v}_{i}^{\mathcal{B}\left({\rm R}\right)}=R^{\top}{\rm v}_{i}^{\mathcal{I}\left({\rm R}\right)}+{\rm b}_{i}^{\mathcal{B}\left({\rm R}\right)}+\omega_{i}^{\mathcal{B}\left({\rm R}\right)}\label{eq:SE3STCH_Vect_R}
\end{equation}
with ${\rm v}_{i}^{\mathcal{I}\left({\rm R}\right)}$ being the $i$th
known vector in the inertial-frame, and ${\rm b}_{i}^{\mathcal{B}\left({\rm R}\right)}$
and $\omega_{i}^{\mathcal{B}\left({\rm R}\right)}$ being unknown
bias and noise components added to the $i$th measurement, respectively,
for all ${\rm v}_{i}^{\mathcal{B}\left({\rm R}\right)},{\rm v}_{i}^{\mathcal{I}\left({\rm R}\right)},{\rm b}_{i}^{\mathcal{B}\left({\rm R}\right)},\omega_{i}^{\mathcal{B}\left({\rm R}\right)}\in\mathbb{R}^{3}$
and $i=1,2,\ldots,N_{{\rm R}}$. %
The known inertial vector ${\rm v}_{i}^{\mathcal{I}\left({\rm R}\right)}$
and the available body-frame measurement ${\rm v}_{i}^{\mathcal{B}\left({\rm R}\right)}$
in \eqref{eq:SE3STCH_Vect_R} can be normalized such that
\begin{equation}
\upsilon_{i}^{\mathcal{I}\left({\rm R}\right)}=\frac{{\rm v}_{i}^{\mathcal{I}\left({\rm R}\right)}}{\left\Vert {\rm v}_{i}^{\mathcal{I}\left({\rm R}\right)}\right\Vert },\hspace{1em}\upsilon_{i}^{\mathcal{B}\left({\rm R}\right)}=\frac{{\rm v}_{i}^{\mathcal{B}\left({\rm R}\right)}}{\left\Vert {\rm v}_{i}^{\mathcal{B}\left({\rm R}\right)}\right\Vert }\label{eq:SE3STCH_Vector_norm}
\end{equation}
Thus, the attitude of a rigid-body can be extracted using $\upsilon_{i}^{\mathcal{I}\left({\rm R}\right)}$
and $\upsilon_{i}^{\mathcal{B}\left({\rm R}\right)}$ in \eqref{eq:SE3STCH_Vector_norm}
rather than ${\rm v}_{i}^{\mathcal{I}\left({\rm R}\right)}$ and ${\rm v}_{i}^{\mathcal{B}\left({\rm R}\right)}$.
Let us introduce the following two sets
\begin{align}
\upsilon^{\mathcal{I}\left({\rm R}\right)} & =\left[\upsilon_{1}^{\mathcal{I}\left({\rm R}\right)},\ldots,\upsilon_{N_{{\rm R}}}^{\mathcal{I}\left({\rm R}\right)}\right]\in\left\{ \mathcal{I}\right\} \nonumber \\
\upsilon^{\mathcal{B}\left({\rm R}\right)} & =\left[\upsilon_{1}^{\mathcal{B}\left({\rm R}\right)},\ldots,\upsilon_{N_{{\rm R}}}^{\mathcal{B}\left({\rm R}\right)}\right]\in\left\{ \mathcal{B}\right\} \label{eq:SE3STCH_Set_R_Norm}
\end{align}
where the two sets in \eqref{eq:SE3STCH_Set_R_Norm} include the normalized
vectors in \eqref{eq:SE3STCH_Vector_norm} for all $\upsilon^{\mathcal{I}\left({\rm R}\right)},\upsilon^{\mathcal{B}\left({\rm R}\right)}\in\mathbb{R}^{3\times N_{{\rm R}}}$.
The position of the moving body can be extracted if its attitude $R$
has already been determined and there exist $N_{{\rm L}}$ known landmarks
(feature points) obtained, for example, by a vision system. The $i$th
body-frame landmark measurement is given by
\[
\left[\begin{array}{c}
{\rm v}_{i}^{\mathcal{B}\left({\rm L}\right)}\\
1
\end{array}\right]=\boldsymbol{T}^{-1}\left[\begin{array}{c}
{\rm v}_{i}^{\mathcal{I}\left({\rm L}\right)}\\
1
\end{array}\right]+\left[\begin{array}{c}
{\rm b}_{i}^{\mathcal{B}\left({\rm L}\right)}\\
0
\end{array}\right]+\left[\begin{array}{c}
\omega_{i}^{\mathcal{B}\left({\rm L}\right)}\\
0
\end{array}\right]
\]
such that
\begin{equation}
{\rm v}_{i}^{\mathcal{B}\left({\rm L}\right)}=R^{\top}\left({\rm v}_{i}^{\mathcal{I}\left({\rm L}\right)}-P\right)+{\rm b}_{i}^{\mathcal{B}\left({\rm L}\right)}+\omega_{i}^{\mathcal{B}\left({\rm L}\right)}\label{eq:SE3STCH_Vec_Landmark}
\end{equation}
where ${\rm v}_{i}^{\mathcal{I}\left({\rm L}\right)}$ is the $i$th
known fixed landmark located in the inertial-frame, ${\rm b}_{i}^{\mathcal{B}\left({\rm L}\right)}$
and $\omega_{i}^{\mathcal{B}\left({\rm L}\right)}$ are the additive
unknown bias and noise vectors of the $i$th measurement, respectively,
for all ${\rm v}_{i}^{\mathcal{B}\left({\rm L}\right)},{\rm v}_{i}^{\mathcal{I}\left({\rm L}\right)},{\rm b}_{i}^{\mathcal{B}\left({\rm L}\right)},\omega_{i}^{\mathcal{B}\left({\rm L}\right)}\in\mathbb{R}^{3}$
and $i=1,2,\ldots,N_{{\rm L}}$. Define the set of inertial-frame
and body-frame vectors associated with landmarks by
\begin{align}
{\rm v}^{\mathcal{B}\left({\rm L}\right)} & =\left[{\rm v}_{1}^{\mathcal{B}\left({\rm L}\right)},\ldots,{\rm v}_{N_{{\rm L}}}^{\mathcal{B}\left({\rm L}\right)}\right]\in\left\{ \mathcal{B}\right\} \nonumber \\
{\rm v}^{\mathcal{I}\left({\rm L}\right)} & =\left[{\rm v}_{1}^{\mathcal{I}\left({\rm L}\right)},\ldots,{\rm v}_{N_{{\rm L}}}^{\mathcal{I}\left({\rm L}\right)}\right]\in\left\{ \mathcal{I}\right\} \label{eq:SE3STCH_Set_L}
\end{align}
In case when more than one landmark is available for measurement,
it is common to obtain a weighted geometric center of all the landmarks,
which can be calculated as follows: 
\begin{align}
\mathcal{G}_{c}^{\mathcal{I}} & =\frac{1}{\sum_{i=1}^{N_{{\rm L}}}k_{i}^{{\rm L}}}\sum_{i=1}^{N_{{\rm L}}}k_{i}^{{\rm L}}{\rm v}_{i}^{\mathcal{I}\left({\rm L}\right)}\label{eq:SE3STCH_Center_Landmark_I}\\
\mathcal{G}_{c}^{\mathcal{B}} & =\frac{1}{\sum_{i=1}^{N_{{\rm L}}}k_{i}^{{\rm L}}}\sum_{i=1}^{N_{{\rm L}}}k_{i}^{{\rm L}}{\rm v}_{i}^{\mathcal{B}\left({\rm L}\right)}\label{eq:SE3STCH_Center_Landmark_B}
\end{align}
such that $k_{i}^{{\rm L}}$ is the confidence level of the $i$th
measurement. 
\begin{assum}
	\label{Assum:SE3STCH_1} (Rigid-body pose observability) The pose
	of a rigid-body in 3D space can be extracted given the availability
	of at least two non-collinear vectors from the sets in \eqref{eq:SE3STCH_Set_R_Norm}
	($N_{{\rm R}}\geq2$) and at least one feature point from the sets
	in \eqref{eq:SE3STCH_Set_L} with $N_{{\rm L}}\geq1$. In the case
	when $N_{{\rm R}}=2$, the third vector can be obtained by the means
	of cross multiplication: $\upsilon_{3}^{\mathcal{I}\left({\rm R}\right)}=\upsilon_{1}^{\mathcal{I}\left({\rm R}\right)}\times\upsilon_{2}^{\mathcal{I}\left({\rm R}\right)}$
	and $\upsilon_{3}^{\mathcal{B}\left({\rm R}\right)}=\upsilon_{1}^{\mathcal{B}\left({\rm R}\right)}\times\upsilon_{2}^{\mathcal{B}\left({\rm R}\right)}$. 
\end{assum}
According to Assumption \ref{Assum:SE3STCH_1} a set of vectorial
measurement described in \eqref{eq:SE3STCH_Set_R_Norm} is sufficient
to have rank 3. Accordingly, the
homogeneous transformation matrix $\boldsymbol{T}$ can be extracted
if Assumption \ref{Assum:SE3STCH_1} is met. For simplicity, the body-frame
vectors ${\rm v}_{i}^{\mathcal{B}\left({\rm R}\right)}$ and ${\rm v}_{i}^{\mathcal{B}\left({\rm L}\right)}$
are considered to be noise and bias free in the stability analysis.
In the Simulation Section, on the contrary, the noise and bias corrupting
the measurements of ${\rm v}_{i}^{\mathcal{B}\left({\rm R}\right)}$
and ${\rm v}_{i}^{\mathcal{B}\left({\rm L}\right)}$ are taken into
consideration. The pose kinematics of the homogeneous transformation
matrix $\boldsymbol{T}$ in \eqref{eq:SE3PPF_T_matrix2} are given
by
\[
\left[\begin{array}{cc}
\dot{R} & \dot{P}\\
\underline{\mathbf{0}}_{3}^{\top} & 0
\end{array}\right]=\left[\begin{array}{cc}
R & P\\
\underline{\mathbf{0}}_{3}^{\top} & 1
\end{array}\right]\left[\begin{array}{cc}
\left[\Omega\right]_{\times} & V\\
\underline{\mathbf{0}}_{3}^{\top} & 0
\end{array}\right]
\]
such that
\begin{align}
\dot{P} & =RV\nonumber \\
\dot{R} & =R\left[\Omega\right]_{\times}\label{eq:SE3PPF_R_Dynamics}\\
\dot{\boldsymbol{T}} & =\boldsymbol{T}\left[\mathcal{Y}\right]_{\wedge}\label{eq:SE3PPF_T_Dynamics}
\end{align}
with $\Omega\in\mathbb{R}^{3}$ and $V\in\mathbb{R}^{3}$ being the
true angular and translational velocity of the moving body, respectively,
and $\mathcal{Y}=\left[\Omega^{\top},V^{\top}\right]^{\top}\in\mathbb{R}^{6}$
being the group velocity vector. The angular velocity can be measured
by a gyroscope, for example, and expressed as follows:
\begin{equation}
\Omega_{m}=\Omega+b_{\Omega}+\omega_{\Omega}\in\left\{ \mathcal{B}\right\} \label{eq:SE3PPF_Angular}
\end{equation}
where $b_{\Omega}$ is an unknown constant or slowly time-varying
bias, and $\omega_{\Omega}$ is an unknown random noise attached to
the measurement, for all $b_{\Omega},\omega_{\Omega}\in\mathbb{R}^{3}$.
Likewise, the translational velocity measurement of a moving body
can be obtained using a GPS, for instance, and defined by 
\begin{equation}
V_{m}=V+b_{V}+\omega_{V}\in\left\{ \mathcal{B}\right\} \label{eq:SE3PPF_V_Trans}
\end{equation}
with $b_{V}\in\mathbb{R}^{3}$ denoting an unknown constant or slowly
time-varying  bias, and $\omega_{V}\in\mathbb{R}^{3}$ being random
noise attached to the translational velocity measurements. The group
of velocity measurements and bias associated with it can be defined
by $\mathcal{Y}_{m}=\left[\Omega_{m}^{\top},V_{m}^{\top}\right]^{\top}\in\mathbb{R}^{6}$
and $b=\left[b_{\Omega}^{\top},b_{V}^{\top}\right]^{\top}\in\mathbb{R}^{6}$,
respectively. For the sake of simplicity, we consider $\omega_{\Omega}=\omega_{V}=\underline{\mathbf{0}}_{3}$
in the analysis. However, in the implementation it is used $\omega_{\Omega}\neq\underline{\mathbf{0}}_{3}$
and $\omega_{V}\neq\underline{\mathbf{0}}_{3}$. Considering the normalized
Euclidean distance of the rotational matrix $R$ in \eqref{eq:SE3STCH_Ecul_Dist}
and the identity in \eqref{eq:SO3PPF_Identity7}, the true attitude
kinematics in \eqref{eq:SE3PPF_R_Dynamics} can be expressed in view
of \eqref{eq:SE3STCH_Ecul_Dist} as
\begin{align}
||\dot{R}||_{I} & =-\frac{1}{4}{\rm Tr}\{\dot{R}\}\nonumber \\
& =-\frac{1}{4}{\rm Tr}\{\boldsymbol{\mathcal{P}}_{a}\left(R\right)\left[\Omega\right]_{\times}\}\nonumber \\
& =\frac{1}{2}\mathbf{vex}(\boldsymbol{\mathcal{P}}_{a}(R))^{\top}\Omega\label{eq:SE3PPF_NormR_dynam}
\end{align}
Accordingly, the problem of pose kinematics in \eqref{eq:SE3PPF_T_Dynamics}
can be reformulated and expressed in vector form as
\begin{equation}
\left[\begin{array}{c}
||\dot{R}||_{I}\\
\dot{P}
\end{array}\right]=\left[\begin{array}{cc}
\frac{1}{2}\mathbf{vex}(\boldsymbol{\mathcal{P}}_{a}(R))^{\top} & \underline{\mathbf{0}}_{3}^{\top}\\
\mathbf{0}_{3\times3} & R
\end{array}\right]\left[\begin{array}{c}
\Omega_{m}-b_{\Omega}\\
V_{m}-b_{V}
\end{array}\right]\label{eq:SE3PPF_T_VEC_Dyn}
\end{equation}

\noindent with $\mathbf{0}_{3\times3}$ being a zero matrix and $\omega_{\Omega}=\omega_{V}=\underline{\mathbf{0}}_{3}$.
Let the estimate of the homogeneous transformation matrix in \eqref{eq:SE3PPF_T_matrix2},
denoted by $\hat{\boldsymbol{T}}$, be given by

\noindent 
\begin{equation}
\hat{\boldsymbol{T}}=\left[\begin{array}{cc}
\hat{R} & \hat{P}\\
\underline{\mathbf{0}}_{3}^{\top} & 1
\end{array}\right]\label{eq:SE3PPF_Test_matrix}
\end{equation}
with $\hat{R}$ and $\hat{P}$ being the estimates of $R$ and $P$,
respectively. Let us define the error in the homogeneous transformation
matrix from body-frame to estimator-frame by
\begin{align}
\tilde{\boldsymbol{T}} & =\hat{\boldsymbol{T}}\boldsymbol{T}^{-1}=\left[\begin{array}{cc}
\tilde{R} & \hat{P}-\tilde{R}P\\
\underline{\mathbf{0}}_{3}^{\top} & 1
\end{array}\right]=\left[\begin{array}{cc}
\tilde{R} & \tilde{P}\\
\underline{\mathbf{0}}_{3}^{\top} & 1
\end{array}\right]\label{eq:SE3PPF_Terr_matrix}
\end{align}

\noindent where $\tilde{R}=\hat{R}R^{\top}$ and $\tilde{P}$ are
the errors associated with attitude and position, respectively. The
aim of this work is to drive $\hat{\boldsymbol{T}}\rightarrow\boldsymbol{T}$
which in turn guarantees driving $\tilde{P}\rightarrow\underline{\mathbf{0}}_{3}$,
$\tilde{R}\rightarrow\mathbf{I}_{3}$, and $\tilde{\boldsymbol{T}}\rightarrow\mathbf{I}_{4}$.
Lemma \ref{Lemm:SE3PPF_1} presented below will prove useful in the
subsequent filter derivation. 
\begin{lem}
	\label{Lemm:SE3PPF_1}Let $R\in\mathbb{SO}\left(3\right)$, $M=M^{\top}\in\mathbb{R}^{3\times3}$,
	$M$ have rank 3, ${\rm Tr}\left\{ M\right\} =3$, and $\bar{\mathbf{M}}={\rm Tr}\left\{ M\right\} \mathbf{I}_{3}-M$,
	while the minimum singular value of $\bar{\mathbf{M}}$ is $\underline{\lambda}:=\underline{\lambda}(\bar{\mathbf{M}})$.
	Then, the following holds:
	\begin{align}
	||\mathbf{vex}(\boldsymbol{\mathcal{P}}_{a}(R))||^{2} & =4(1-||R||_{I})||R||_{I}\label{eq:SE3PPF_lemm1_1}\\
	\frac{2}{\underline{\lambda}}\frac{||\mathbf{vex}(\boldsymbol{\mathcal{P}}_{a}(RM))||^{2}}{1+{\rm Tr}\{RMM^{-1}\}} & \geq\left\Vert RM\right\Vert _{I}\label{eq:SE3PPF_lemm1_2}
	\end{align}
	\textbf{Proof. See \nameref{sec:SO3STCH_EXPL_AppendixA}.} 
\end{lem}

\subsection{Prescribed Performance \label{subsec:SE3PPF_Prescribed-Performance}}

\begin{figure}[h!]
	\centering{}\includegraphics[scale=0.27]{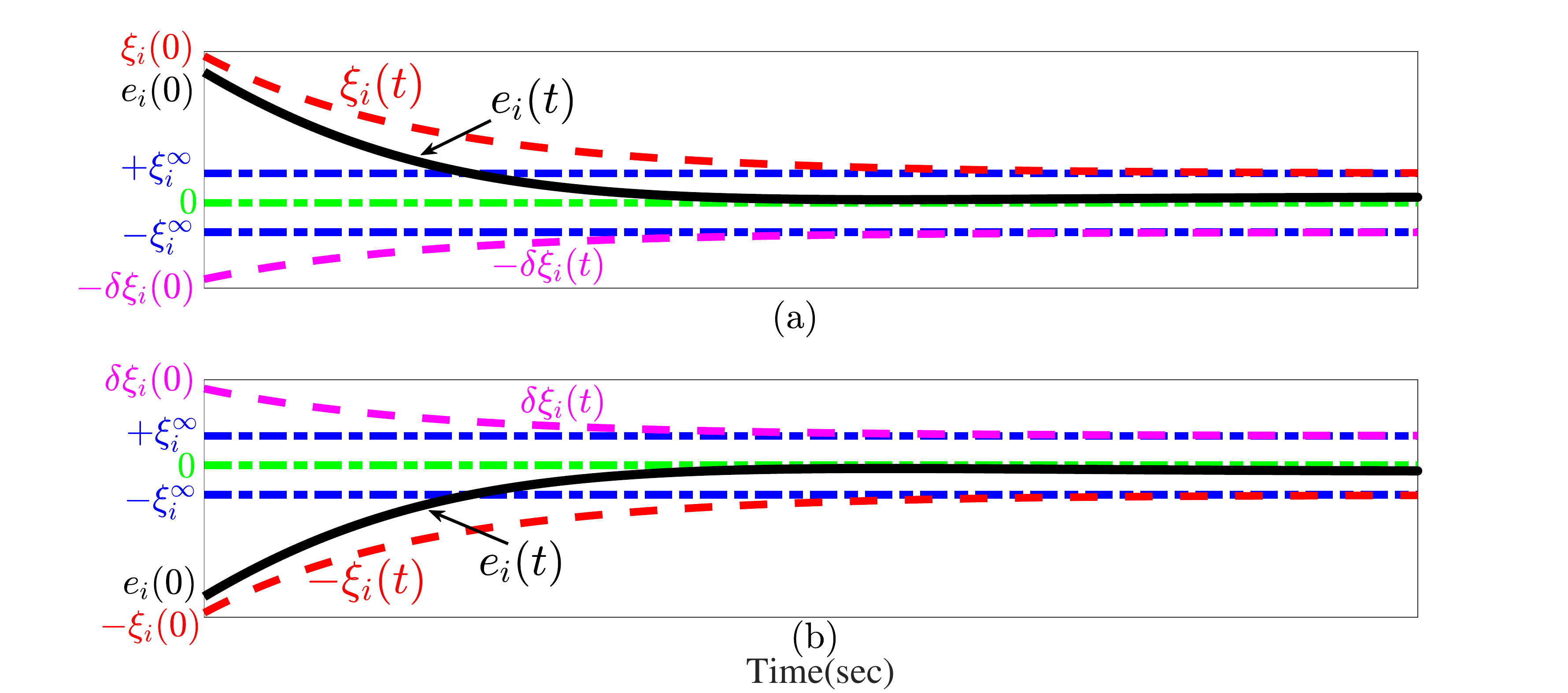} \caption{Graphical representation of the systematic convergence of tracking
		error $\boldsymbol{e}_{i}$ with PPF satisfying (a) Eq. \eqref{eq:SE3PPF_ePos};
		(b) Eq. \eqref{eq:SE3PPF_eNeg}.}
	\label{fig:SO3PPF_2} 
\end{figure}

Considering the error in the homogeneous transformation matrix as
in \eqref{eq:SE3PPF_Terr_matrix} and in view of the pose kinematics
in \eqref{eq:SE3PPF_T_VEC_Dyn}, let us define the error in vector
form by
\begin{equation}
\boldsymbol{e}=\left[\boldsymbol{e}_{1},\boldsymbol{e}_{2},\boldsymbol{e}_{3},\boldsymbol{e}_{4}\right]^{\top}=\left[||\tilde{R}||_{I},\tilde{P}^{\top}\right]^{\top}\in\mathbb{R}^{4}\label{eq:SE3PPF_Vec_error}
\end{equation}
The objective of this subsection is to reformulate the problem such
that the error in \eqref{eq:SE3PPF_Vec_error} satisfies transient
as well as steady-state measures predefined by the user. This can
be achieved by guiding the error vector $\boldsymbol{e}$ to initiate
within a large known set and after decaying smoothly and systematically
settle within a predefined small set using prescribed performance
function (PPF) \cite{bechlioulis2008robust,hashim2017neuro}. The
PPF is defined by $\xi_{i}\left(t\right)$ which is a positive smooth
time-decreasing function which satisfies $\xi_{i}:\mathbb{R}_{+}\to\mathbb{R}_{+}$
and $\lim\limits _{t\to\infty}\xi_{i}\left(t\right)=\xi_{i}^{\infty}>0$
and can be expressed by \cite{bechlioulis2008robust}
\begin{equation}
\xi_{i}\left(t\right)=\left(\xi_{i}^{0}-\xi_{i}^{\infty}\right)\exp\left(-\ell_{i}t\right)+\xi_{i}^{\infty}\label{eq:SE3PPF_Presc}
\end{equation}
with $\xi_{i}\left(0\right)=\xi_{i}^{0}$ being the initial value
of the PPF and the upper bound of the known large set, $\xi_{i}^{\infty}$
being the upper bound of the narrow set, and $\ell_{i}$ being a positive
constant controlling the convergence rate of $\xi\left(t\right)$
from $\xi_{i}^{0}$ to $\xi_{i}^{\infty}$ for all $i=1,\ldots,4$.
The error $\boldsymbol{e}_{i}\left(t\right)$ is guaranteed to follow
the predefined transient and steady-state boundaries, if the conditions
below are met:
\begin{align}
-\delta\xi_{i}\left(t\right)<\boldsymbol{e}_{i}\left(t\right)<\xi_{i}\left(t\right), & \text{ if }\boldsymbol{e}_{i}\left(0\right)>0\label{eq:SE3PPF_ePos}\\
-\xi_{i}\left(t\right)<\boldsymbol{e}_{i}\left(t\right)<\delta\xi_{i}\left(t\right), & \text{ if }\boldsymbol{e}_{i}\left(0\right)<0\label{eq:SE3PPF_eNeg}
\end{align}
with $\delta\in\left[0,1\right]$. For clarity, define $\boldsymbol{e}_{i}:=\boldsymbol{e}_{i}\left(t\right)$
and $\xi_{i}:=\xi_{i}\left(t\right)$. Also, let us define $\xi=[\xi_{1},\xi_{2},\xi_{3},\xi_{4}]^{\top}$,
$\ell=[\ell_{1},\ell_{2},\ell_{3},\ell_{4}]^{\top}$, $\xi^{0}=\left[\xi_{1}^{0},\xi_{2}^{0},\xi_{3}^{0},\xi_{4}^{0}\right]^{\top}$,
and $\xi^{\infty}=\left[\xi_{1}^{\infty},\xi_{2}^{\infty},\xi_{3}^{\infty},\xi_{4}^{\infty}\right]^{\top}$
for all $\xi,\ell,\xi^{0},\xi^{\infty}\in\mathbb{R}^{4}$. The systematic
convergence of the tracking error $\boldsymbol{e}_{i}$, from a given
large set to a given narrow set in accordance with \eqref{eq:SE3PPF_ePos}
and \eqref{eq:SE3PPF_eNeg} is depicted in Fig. \ref{fig:SO3PPF_2}.

\begin{rem}
	\label{SE3PPF_rem3} In accordance with the discussion in \cite{bechlioulis2008robust,hashim2017neuro},
	knowing the upper bound and the sign of $\boldsymbol{e}_{i}\left(0\right)$
	is sufficient to force the error to satisfy the performance constraints
	and maintain the error regulation within predefined dynamically reducing
	boundaries for all $t>0$. If the condition in \eqref{eq:SE3PPF_ePos}
	or \eqref{eq:SE3PPF_eNeg} is met, the maximum overshoot is sufficient
	to be bounded by $\pm\delta\xi_{i}$, the steady-state error is bounded
	by $\pm\xi_{i}^{\infty}$, and $|\boldsymbol{e}_{i}|$ is trapped
	between $\xi_{i}$ and $\delta\xi_{i}$ as presented in Fig. \ref{fig:SO3PPF_2}. 
\end{rem}
Define the error $\boldsymbol{e}_{i}$ by 
\begin{equation}
\boldsymbol{e}_{i}=\xi_{i}\mathcal{Z}(\mathcal{E}_{i})\label{eq:SE3PPF_e_Trans}
\end{equation}
where $\xi_{i}\in\mathbb{R}$ is defined in \eqref{eq:SE3PPF_Presc},
$\mathcal{E}_{i}\in\mathbb{R}$ is a relaxed form of the constrained
error referred to as transformed error, and $\mathcal{Z}(\mathcal{E}_{i})$
is a smooth function that behaves according to Assumption \ref{Assum:SE3PPF_1}: 
\begin{assum}
	\label{Assum:SE3PPF_1}The smooth function $\mathcal{Z}(\mathcal{E}_{i})$
	has the following properties \cite{bechlioulis2008robust}: 
\begin{enumerate}
	\item[\textbf{1. }] $\mathcal{Z}(\mathcal{E}_{i})$ is smooth and strictly increasing. 
	\item[\textbf{2. }] $\mathcal{Z}(\mathcal{E}_{i})$ is constrained by the following two
	bounds \\
	$-\underline{\delta}_{i}<\mathcal{Z}(\mathcal{E}_{i})<\bar{\delta}_{i},{\rm \text{ if }}\boldsymbol{e}_{i}\left(0\right)\geq0$\\
	$-\bar{\delta}_{i}<\mathcal{Z}(\mathcal{E}_{i})<\underline{\delta}_{i},{\rm \text{ if }}\boldsymbol{e}_{i}\left(0\right)<0$
	\\
	with $\bar{\delta}_{i}$ and $\underline{\delta}_{i}$ being positive
	constants satisfy $\underline{\delta}_{i}\leq\bar{\delta}_{i}$. 
	\item[\textbf{3. }] 
	\item[] $\left.\begin{array}{c}
	\underset{\mathcal{E}_{i}\rightarrow-\infty}{\lim}\mathcal{Z}(\mathcal{E}_{i})=-\underline{\delta}_{i}\\
	\underset{\mathcal{E}_{i}\rightarrow+\infty}{\lim}\mathcal{Z}(\mathcal{E}_{i})=\bar{\delta}_{i}
	\end{array}\right\} {\rm \text{ if }}\boldsymbol{e}_{i}\geq0$\\
	$\left.\begin{array}{c}
	\underset{\mathcal{E}_{i}\rightarrow-\infty}{\lim}\mathcal{Z}(\mathcal{E}_{i})=-\bar{\delta}_{i}\\
	\underset{\mathcal{E}_{i}\rightarrow+\infty}{\lim}\mathcal{Z}(\mathcal{E}_{i})=\underline{\delta}_{i}
	\end{array}\right\} {\rm \text{ if }}\boldsymbol{e}_{i}<0$ \\
	such that 
\end{enumerate}
\end{assum}
\begin{equation}
\mathcal{Z}\left(\mathcal{E}_{i}\right)=\begin{cases}
\frac{\bar{\delta}_{i}\exp(\mathcal{E}_{i})-\underline{\delta}_{i}\exp(-\mathcal{E}_{i})}{\exp(\mathcal{E}_{i})+\exp(-\mathcal{E}_{i})}, & \bar{\delta}_{i}\geq\underline{\delta}_{i}\text{ if }\boldsymbol{e}_{i}\geq0\\
\frac{\bar{\delta}_{i}\exp(\mathcal{E}_{i})-\underline{\delta}_{i}\exp(-\mathcal{E}_{i})}{\exp(\mathcal{E}_{i})+\exp(-\mathcal{E}_{i})}, & \underline{\delta}_{i}\geq\bar{\delta}_{i}\text{ if }\boldsymbol{e}_{i}<0
\end{cases}\label{eq:SE3PPF_Smooth}
\end{equation}

The transformed error could be extracted through the inverse transformation
of \eqref{eq:SE3PPF_Smooth} 
\begin{equation}
\mathcal{E}_{i}(\boldsymbol{e}_{i},\xi_{i})=\mathcal{Z}^{-1}(\boldsymbol{e}_{i}/\xi_{i})\label{eq:SO3PPF_Trans1}
\end{equation}
with $\mathcal{E}_{i}\in\mathbb{R}$, $\mathcal{Z}\in\mathbb{R}$
and $\mathcal{Z}^{-1}\in\mathbb{R}$ being smooth functions. For simplicity,
let $\mathcal{E}_{i}:=\mathcal{E}_{i}(\cdot,\cdot)$, $\bar{\delta}=[\bar{\delta}_{1},\bar{\delta}_{2},\bar{\delta}_{3},\bar{\delta}_{4}]^{\top}$,
$\underline{\delta}=[\underline{\delta}_{1},\underline{\delta}_{2},\underline{\delta}_{3},\underline{\delta}_{4}]^{\top}$,
$\mathcal{E}=\left[\mathcal{E}_{R},\mathcal{E}_{P}^{\top}\right]^{\top}$
for all $\bar{\delta},\underline{\delta},\mathcal{E}\in\mathbb{R}^{4}$
with $\mathcal{E}_{R}=\mathcal{E}_{1}\in\mathbb{R}$ and $\mathcal{E}_{P}=[\mathcal{E}_{2},\mathcal{E}_{3},\mathcal{E}_{4}]^{\top}\in\mathbb{R}^{3}$.
In fact, the transformed error $\mathcal{E}_{i}$ translates $\boldsymbol{e}_{i}$
from the given constrained form in \eqref{eq:SE3PPF_ePos} or \eqref{eq:SE3PPF_eNeg}
to its unconstrained form as in \eqref{eq:SO3PPF_Trans1}. From \eqref{eq:SE3PPF_Smooth},
the inverse transformation can be expressed as
\begin{equation}
\begin{aligned}\mathcal{E}_{i}= & \frac{1}{2}\begin{cases}
\text{ln}\frac{\underline{\delta}_{i}+\boldsymbol{e}_{i}/\xi_{i}}{\bar{\delta}_{i}-\boldsymbol{e}_{i}/\xi_{i}}, & \bar{\delta}_{i}\geq\underline{\delta}_{i}\text{ if }\boldsymbol{e}_{i}\geq0\\
\text{ln}\frac{\underline{\delta}_{i}+\boldsymbol{e}_{i}/\xi_{i}}{\bar{\delta}_{i}-\boldsymbol{e}_{i}/\xi_{i}}, & \underline{\delta}_{i}\geq\bar{\delta}_{i}\text{ if }\boldsymbol{e}_{i}<0
\end{cases}\end{aligned}
\label{eq:SE3PPF_trans2}
\end{equation}

\begin{rem}
	\label{rem:SO3PPF_1}
	Consider the transformed error in \eqref{eq:SE3PPF_trans2}. The
	transient and steady-state performance of the tracking error ($\boldsymbol{e}_{i}$)
	is bounded by the performance function $\xi_{i}$, and therefore,
	the prescribed performance is achieved if and only if $\mathcal{E}_{i}$
	is guaranteed to be bounded for all $t\geq0$. 
\end{rem}
\begin{prop}
	\label{Prop:SE3PPF_1}Consider the error vector in \eqref{eq:SE3PPF_Vec_error}
	with the normalized Euclidean distance error $||\tilde{R}||_{I}$
	being given by \eqref{eq:SE3STCH_Ecul_Dist}. From \eqref{eq:SE3PPF_e_Trans},
	\eqref{eq:SE3PPF_Smooth}, and \eqref{eq:SO3PPF_Trans1} let the transformed
	error be expressed as in \eqref{eq:SE3PPF_trans2} provided that $\underline{\delta}=\bar{\delta}$.
	Then the following statements are true. 
\begin{enumerate}
	\item[(i)] The only possible representation of $\mathcal{E}_{1}$ is as follows:
	\begin{equation}
	\mathcal{E}_{1}=\frac{1}{2}\text{ln}\frac{\underline{\delta}_{1}+\boldsymbol{e}_{1}/\xi_{1}}{\bar{\delta}_{1}-\boldsymbol{e}_{1}/\xi_{1}}=\frac{1}{2}\text{ln}\frac{\underline{\delta}_{1}+||\tilde{R}||_{I}/\xi_{1}}{\bar{\delta}_{1}-||\tilde{R}||_{I}/\xi_{1}}\label{eq:SE3PPF_trans3}
	\end{equation}
	\item[(ii)] The transformed error $\mathcal{E}_{1}>0\forall||\tilde{R}||_{I}\neq0$. 
	\item[(iii)] $\mathcal{E}=\underline{\mathbf{0}}_{4}$ only at $\boldsymbol{e}=\underline{\mathbf{0}}_{4}$
	and the critical point of $\mathcal{E}$ satisfies $\boldsymbol{e}=\underline{\mathbf{0}}_{4}$. 
	\item[(iv)] The only critical point of $\mathcal{E}$ is $\tilde{\boldsymbol{T}}=\mathbf{I}_{4}$. 
\end{enumerate}
\end{prop}
\textbf{Proof.} Given that $0\leq||\tilde{R}\left(t\right)||_{I}\leq1,\forall t\geq0$
as defined in \eqref{eq:SE3STCH_Ecul_Dist}, one can find that the
upper part of \eqref{eq:SE3PPF_trans2} holds $\forall t\geq0$ which
proves (i). Since $\underline{\delta}=\bar{\delta}$ with the constraint
$||\tilde{R}||_{I}\leq\xi_{1}$, the expression in \eqref{eq:SE3PPF_trans3}
is $(\underline{\delta}_{1}+||\tilde{R}||_{I}/\xi_{1})/(\bar{\delta}_{1}-||\tilde{R}||_{I}/\xi_{1})\geq1\forall||\tilde{R}||_{I}\neq0$.
Thus, $\mathcal{E}_{1}>0\forall||\tilde{R}||_{I}\neq0$ which confirms
(ii). Considering $\underline{\delta}=\bar{\delta}$ with the constraint
$\boldsymbol{e}_{i}\leq\xi_{i}$, it is clear that $(\underline{\delta}_{i}+\boldsymbol{e}_{i}/\xi_{i})/(\bar{\delta}_{i}-\boldsymbol{e}_{i}/\xi_{i})=1$
if and only if $\boldsymbol{e}_{i}=0$. Accordingly, $\mathcal{E}_{i}\neq0\forall\boldsymbol{e}_{i}\neq0$
and $\mathcal{E}_{i}=0$ only at $\boldsymbol{e}_{i}=0$ which proves
(iii). For (iv), from \eqref{eq:SE3STCH_Ecul_Dist} and \eqref{eq:SE3PPF_Terr_matrix},
$||\tilde{R}||_{I}=0$ and $\tilde{P}=0$ if and only if $\tilde{\boldsymbol{T}}=\mathbf{I}_{4}$.
Thus, the critical point of $\mathcal{E}$ satisfies $||\tilde{R}||_{I}=0$
and $\tilde{P}=0$ which in turn satisfies $\tilde{\boldsymbol{T}}=\mathbf{I}_{4}$
and proves (iv). Define $\mu_{i}:=\mu_{i}\left(\boldsymbol{e}_{i},\xi_{i}\right)$
such that
\begin{equation}
\begin{split}\mu_{i} & =\frac{1}{2\xi_{i}}\frac{\partial\mathcal{Z}^{-1}\left(\boldsymbol{e}_{i}/\xi_{i}\right)}{\partial\left(\boldsymbol{e}_{i}/\xi_{i}\right)}=\frac{1}{2\xi_{i}}\left(\frac{1}{\underline{\delta}_{i}+\boldsymbol{e}_{i}/\xi_{i}}+\frac{1}{\bar{\delta}_{i}-\boldsymbol{e}_{i}/\xi_{i}}\right)\end{split}
\label{eq:SE3PPF_mu}
\end{equation}
Hence, one can find that the derivative of $\dot{\mathcal{E}_{i}}$
is as follows:
\begin{align}
\dot{\mathcal{E}_{i}} & =\mu_{i}\left(\dot{\boldsymbol{e}}_{i}-\frac{\dot{\xi}_{i}}{\xi_{i}}\boldsymbol{e}_{i}\right)\label{eq:SE3PPF_Trans_dot1}
\end{align}
More simply, the expression in \eqref{eq:SE3PPF_Trans_dot1} is
\begin{equation}
\dot{\mathcal{E}}=\left[\begin{array}{cc}
\boldsymbol{\Psi}_{R} & \underline{\mathbf{0}}_{3}^{\top}\\
\underline{\mathbf{0}}_{3} & \boldsymbol{\Psi}_{P}
\end{array}\right]\left(\dot{\boldsymbol{e}}-\left[\begin{array}{cc}
\boldsymbol{\Lambda}_{R} & \underline{\mathbf{0}}_{3}^{\top}\\
\underline{\mathbf{0}}_{3} & \boldsymbol{\Lambda}_{P}
\end{array}\right]\boldsymbol{e}\right)\label{eq:SE3PPF_Trans_dot}
\end{equation}
with $\boldsymbol{\Lambda}_{R}=\frac{\dot{\xi}_{1}}{\xi_{1}}$, $\boldsymbol{\Lambda}_{P}={\rm diag}\left(\frac{\dot{\xi}_{2}}{\xi_{2}},\frac{\dot{\xi}_{3}}{\xi_{3}},\frac{\dot{\xi}_{4}}{\xi_{4}}\right)$,
$\boldsymbol{\Psi}_{R}=\mu_{1}$, and $\boldsymbol{\Psi}_{P}={\rm diag}(\mu_{2},\mu_{3},\mu_{4})$
for all $\boldsymbol{\Lambda}_{R},\boldsymbol{\Psi}_{R}\in\mathbb{R}$
and $\boldsymbol{\Lambda}_{P},\boldsymbol{\Psi}_{P}\in\mathbb{R}^{3\times3}$.
The following section introduces two nonlinear pose filters on $\mathbb{SE}\left(3\right)$
with prescribed performance characteristics which for $0\leq\left|\boldsymbol{e}_{i}\left(0\right)\right|<\xi_{i}\left(0\right)$
guarantee $\mathcal{E}_{i}\in\mathcal{L}_{\infty},\forall t\geq0$
and, therefore, satisfy \eqref{eq:SE3PPF_ePos} or \eqref{eq:SE3PPF_eNeg}.

\section{Nonlinear Complementary Pose Filters On $\mathbb{SE}\left(3\right)$
	with Prescribed Performance \label{sec:SE3PPF-Filters}}

This section aims to provide a comprehensive description of the two
nonlinear complementary pose filters evolved on $\mathbb{SE}\left(3\right)$
with the error vector, introduced in \eqref{eq:SE3PPF_Vec_error},
behaving in accordance with the predefined transient as well as steady-state
measures specified by the user. %
The first proposed filter is named a semi-direct pose filter with
prescribed performance and the second one is termed a direct pose
filter with prescribed performance. The difference between the two
lies in the fact that while the semi-direct filter requires both attitude
and position to be reconstructed through a set of vectorial measurements
given in \eqref{eq:SE3STCH_Set_R_Norm} and \eqref{eq:SE3STCH_Set_L}
combined with the measurement of the group velocity vector as described
in \eqref{eq:SE3PPF_Angular} and \eqref{eq:SE3PPF_V_Trans}, the
direct filter only utilizes the above-mentioned measurements in the
filter design. The structure of the proposed pose filters described
in the two subsequent subsections is nonlinear on the Lie group of
$\mathbb{SE}\left(3\right)$ and is given by 
\[
\dot{\hat{\boldsymbol{T}}}=\hat{\boldsymbol{T}}[\hat{\mathcal{Y}}]_{\wedge}
\]
with $\hat{\mathcal{Y}}=[\hat{\Omega}^{\top},\hat{V}^{\top}]\in\mathbb{R}^{6}$
such that $\dot{\hat{R}}=\hat{R}[\hat{\Omega}]_{\times}$ and $\dot{\hat{P}}=\hat{R}\hat{V}$.

\subsection{Semi-direct Pose Filter with Prescribed Performance \label{subsec:SE3PPF_Passive-Filter}}

Recall the error in \eqref{eq:SE3PPF_Vec_error} $\boldsymbol{e}=\left[||\tilde{R}||_{I},\tilde{P}^{\top}\right]^{\top}$.
Define $\boldsymbol{T}_{y}=\left[\begin{array}{cc}
R_{y} & P_{y}\\
\underline{\mathbf{0}}_{3}^{\top} & 1
\end{array}\right]$ as a reconstructed homogeneous transformation matrix of the true
$\boldsymbol{T}$. $R_{y}$ corrupted by uncertain measurements can
be reconstructed as in \cite{shuster1981three,markley1988attitude}
or for simplicity visit the Appendix in \cite{hashim2018SO3Stochastic,hashim2018SE3Stochastic}.
From \eqref{eq:SE3STCH_Center_Landmark_I} and \eqref{eq:SE3STCH_Center_Landmark_B}
$P_{y}$ is reconstructed in the following manner
\begin{align}
P_{y} & =\frac{1}{\sum_{i=1}^{N_{{\rm L}}}k_{i}^{{\rm L}}}\sum_{i=1}^{N_{{\rm L}}}k_{i}^{{\rm L}}\left({\rm v}_{i}^{\mathcal{I}\left({\rm L}\right)}-R_{y}{\rm v}_{i}^{\mathcal{B}\left({\rm L}\right)}\right)\nonumber \\
& =\mathcal{G}_{c}^{\mathcal{I}}-R_{y}\mathcal{G}_{c}^{\mathcal{B}}\label{eq:SE3PPF_Py}
\end{align}
Consider the following pose filter design
\begin{align}
\dot{\hat{R}}= & \hat{R}\left[\Omega_{m}-\hat{b}_{\Omega}-\hat{R}^{\top}W_{\Omega}\right]_{\times}\label{eq:SE3PPF_Rest_dot_Ty}\\
\dot{\hat{P}}= & \hat{R}(V_{m}-\hat{b}_{V}-W_{V})\label{eq:SE3PPF_Pest_dot_Ty}\\
\dot{\hat{b}}_{\Omega}= & \frac{\gamma}{2}\boldsymbol{\Psi}_{R}\mathcal{E}_{R}\hat{R}^{\top}\mathbf{vex}(\boldsymbol{\mathcal{P}}_{a}(\tilde{R}))\nonumber \\
& +\gamma\hat{R}^{\top}\left[\tilde{P}-\hat{P}\right]_{\times}\boldsymbol{\Psi}_{P}\mathcal{E}_{P}\label{eq:SE3PPF_b1est_dot_Ty}\\
\dot{\hat{b}}_{V}= & \gamma\hat{R}^{\top}\boldsymbol{\Psi}_{P}\mathcal{E}_{P}\label{eq:SE3PPF_b2est_dot_Ty}\\
W_{\Omega}= & 2\frac{k_{w}\boldsymbol{\Psi}_{R}\mathcal{E}_{R}-\boldsymbol{\Lambda}_{R}/4}{1-||\tilde{R}||_{I}}\mathbf{vex}(\boldsymbol{\mathcal{P}}_{a}(\tilde{R}))\label{eq:SE3PPF_W1est_dot_Ty}\\
W_{V}= & \hat{R}^{\top}\left(k_{w}\boldsymbol{\Psi}_{P}\mathcal{E}_{P}+\left[\tilde{P}-\hat{P}\right]_{\times}W_{\Omega}-\boldsymbol{\Lambda}_{P}\tilde{P}\right)\label{eq:SE3PPF_W2est_dot_Ty}
\end{align}
with $\tilde{R}=\hat{R}R_{y}^{\top}$, $\tilde{P}=\hat{P}-\tilde{R}P_{y}$,
$\mathcal{E}_{R}$, $\mathcal{E}_{P}$, $\boldsymbol{\Psi}_{R}$ and
$\boldsymbol{\Psi}_{P}$ being defined in \eqref{eq:SE3PPF_mu}, and
\eqref{eq:SE3PPF_Trans_dot1}, respectively, $k_{w}$ and $\gamma$
being positive constants, and each of $\hat{b}_{\Omega}$ and $\hat{b}_{V}$
being the estimates of $b_{\Omega}$ and $b_{V}$, respectively. The equivalent quaternion representation and complete implementation
steps of the semi-direct filter are given in \nameref{sec:SO3_PPF_STCH_AppendixB}.

Define the error between the true and the estimated bias by
\begin{align}
\tilde{b}_{\Omega} & =b_{\Omega}-\hat{b}_{\Omega}\label{eq:SE3PPF_b1_tilde}\\
\tilde{b}_{V} & =b_{V}-\hat{b}_{V}\label{eq:SE3PPF_b2_tilde}
\end{align}
where $\tilde{b}=\left[\tilde{b}_{\Omega}^{\top},\tilde{b}_{V}^{\top}\right]^{\top}\in\mathbb{R}^{6}$
is the group error bias vector.
\begin{thm}
	\textbf{\label{thm:SE3PPF_1} }Consider the pose dynamics in \eqref{eq:SE3PPF_T_Dynamics},
	the group of noise-free velocity measurements in \eqref{eq:SE3PPF_Angular}
	and \eqref{eq:SE3PPF_V_Trans} such that $\Omega_{m}=\Omega+b_{\Omega}$
	and $V_{m}=V+b_{V}$, in addition to other vector measurements given
	in \eqref{eq:SE3STCH_Set_R_Norm} and \eqref{eq:SE3STCH_Set_L} coupled
	with the filter kinematics in \eqref{eq:SE3PPF_Rest_dot_Ty}, \eqref{eq:SE3PPF_Pest_dot_Ty},
	\eqref{eq:SE3PPF_b1est_dot_Ty}, \eqref{eq:SE3PPF_b2est_dot_Ty},
	\eqref{eq:SE3PPF_W1est_dot_Ty}, and \eqref{eq:SE3PPF_W2est_dot_Ty}.
	Let Assumption \ref{Assum:SE3STCH_1} hold. Define $\mathcal{U}\subseteq\mathbb{SE}\left(3\right)\times\mathbb{R}^{6}$
	by $\mathcal{U}:=\left\{ \left.(\tilde{\boldsymbol{T}}\left(0\right),\tilde{b}\left(0\right))\right|{\rm Tr}\{\tilde{R}\left(0\right)\}=-1,\tilde{P}\left(0\right)=\underline{\mathbf{0}}_{3},\tilde{b}\left(0\right)=\underline{\mathbf{0}}_{6}\right\} $.
	From almost any initial condition such that ${\rm Tr}\{\tilde{R}\left(0\right)\}\notin\mathcal{U}$
	and $\mathcal{E}\left(0\right)\in\mathcal{L}_{\infty}$, all signals
	in the closed loop are bounded, $\lim_{t\rightarrow\infty}\mathcal{E}\left(t\right)=0$,
	and $\tilde{\boldsymbol{T}}$ asymptotically approaches $\mathbf{I}_{4}$. 
\end{thm}
Theorem \ref{thm:SE3PPF_1} guarantees that the observer pose dynamics
in \eqref{eq:SE3PPF_Rest_dot_Ty}, \eqref{eq:SE3PPF_Pest_dot_Ty},
\eqref{eq:SE3PPF_b1est_dot_Ty}, \eqref{eq:SE3PPF_b2est_dot_Ty},
\eqref{eq:SE3PPF_W1est_dot_Ty}, and \eqref{eq:SE3PPF_W2est_dot_Ty}
are stable with $\mathcal{E}\left(t\right)$ asymptotically approaching
the origin. Since $\mathcal{E}\left(t\right)$ is bounded, the error
vector $\boldsymbol{e}$ in \eqref{eq:SE3PPF_Vec_error} is constrained
by the transient and steady-state boundaries introduced in \eqref{eq:SE3PPF_Presc}.

\textbf{Proof. }Consider the error in the homogeneous transformation
matrix from body-frame to estimator-frame defined as \eqref{eq:SE3PPF_Terr_matrix}.
From \eqref{eq:SE3PPF_R_Dynamics} and \eqref{eq:SE3PPF_Rest_dot_Ty}
the error dynamics are
\begin{align}
\dot{\tilde{R}} & =\hat{R}\left[\tilde{b}_{\Omega}-\hat{R}^{\top}W_{\Omega}\right]_{\times}R^{\top}=\left[\hat{R}\tilde{b}_{\Omega}-W_{\Omega}\right]_{\times}\tilde{R}\label{eq:SE3PPF_Rtilde_dot}
\end{align}
where $\left[\hat{R}\tilde{b}_{\Omega}\right]_{\times}=\hat{R}\left[\tilde{b}_{\Omega}\right]_{\times}\hat{R}^{\top}$
as given in identity \eqref{eq:SO3PPF_Identity2}. In view of \eqref{eq:SE3PPF_R_Dynamics}
and \eqref{eq:SE3PPF_NormR_dynam}, one can express the error dynamics
in \eqref{eq:SE3PPF_Rtilde_dot} in terms of normalized Euclidean
distance as
\begin{align}
\frac{d}{dt}||\tilde{R}||_{I} & =\frac{d}{dt}\frac{1}{4}{\rm Tr}\{\mathbf{I}_{3}-\tilde{R}\}\nonumber \\
& =-\frac{1}{4}{\rm Tr}\left\{ \left[\hat{R}\tilde{b}_{\Omega}-W_{\Omega}\right]_{\times}\boldsymbol{\mathcal{P}}_{a}(\tilde{R})\right\} \nonumber \\
& =\frac{1}{2}\mathbf{vex}(\boldsymbol{\mathcal{P}}_{a}(\tilde{R}))^{\top}(\hat{R}\tilde{b}_{\Omega}-W_{\Omega})\label{eq:SE3PPF_NormRtilde_dot}
\end{align}
with ${\rm Tr}\left\{ \tilde{R}\left[\tilde{b}-W\right]_{\times}\right\} =-2\mathbf{vex}(\boldsymbol{\mathcal{P}}_{a}(\tilde{R}))^{\top}(\tilde{b}-W)$
being defined in \eqref{eq:SO3PPF_Identity7}. Since the position
error is given by $\tilde{P}=\hat{P}-\tilde{R}P$ in \eqref{eq:SE3PPF_Terr_matrix},
one can find the derivative of $\tilde{P}$ to be 
\begin{align}
\dot{\tilde{P}} & =\dot{\hat{P}}-\dot{\tilde{R}}P-\tilde{R}\dot{P}\nonumber \\
& =\dot{\hat{P}}-\left[\hat{R}\tilde{b}_{\Omega}-W_{\Omega}\right]_{\times}\tilde{R}P-\tilde{R}R(V_{m}-b_{V})\nonumber \\
& =\hat{R}(\tilde{b}_{V}-W_{V})+\left[\hat{P}-\tilde{P}\right]_{\times}(\hat{R}\tilde{b}_{\Omega}-W_{\Omega})\label{eq:SE3PPF_Ptilde_dot}
\end{align}
with $\left[\hat{R}\tilde{b}_{\Omega}\right]_{\times}\hat{P}=-\left[\hat{P}\right]_{\times}\hat{R}\tilde{b}_{\Omega}$.
From \eqref{eq:SE3PPF_NormRtilde_dot} and \eqref{eq:SE3PPF_Ptilde_dot},
and in view of \eqref{eq:SE3PPF_T_VEC_Dyn}, the dynamics of the error
vector in \eqref{eq:SE3PPF_Vec_error} become
\begin{align}
\left[\begin{array}{c}
||\dot{\tilde{R}}||_{I}\\
\dot{\tilde{P}}
\end{array}\right] & =\left[\begin{array}{cc}
\frac{1}{2}\mathbf{vex}(\boldsymbol{\mathcal{P}}_{a}(\tilde{R}))^{\top} & \underline{\mathbf{0}}_{3}^{\top}\\
\left[\hat{P}-\tilde{P}\right]_{\times} & \hat{R}
\end{array}\right]\left[\begin{array}{c}
\hat{R}\tilde{b}_{\Omega}-W_{\Omega}\\
\tilde{b}_{V}-W_{V}
\end{array}\right]\label{eq:SE3PPF_VEC_tilde_dot}
\end{align}
Accordingly, the derivative of the transformed error in \eqref{eq:SE3PPF_Trans_dot}
can be represented with direct substitution of $\boldsymbol{e}=\left[||\tilde{R}||_{I},\tilde{P}^{\top}\right]^{\top}$
in addition to the result in \eqref{eq:SE3PPF_VEC_tilde_dot}. Now,
consider the following candidate Lyapunov function

\begin{align}
V(\mathcal{E},\tilde{b}_{\Omega},\tilde{b}_{V}) & =\frac{1}{2}||\mathcal{E}||^{2}+\frac{1}{2\gamma}||\tilde{b}_{\Omega}||^{2}+\frac{1}{2\gamma}||\tilde{b}_{V}||^{2}\label{eq:SE3PPF_V_Ry}
\end{align}
Differentiating $V:=V(\mathcal{E},\tilde{b}_{\Omega},\tilde{b}_{V})$
in \eqref{eq:SE3PPF_V_Ry} results in 
\begin{align}
\dot{V}= & \mathcal{E}^{\top}\dot{\mathcal{E}}-\frac{1}{\gamma}\tilde{b}_{\Omega}^{\top}\dot{\hat{b}}_{\Omega}-\frac{1}{\gamma}\tilde{b}_{V}^{\top}\dot{\hat{b}}_{V}\nonumber \\
= & \mathcal{E}_{R}\boldsymbol{\Psi}_{R}\left(\frac{1}{2}\mathbf{vex}(\boldsymbol{\mathcal{P}}_{a}(\tilde{R}))^{\top}(\hat{R}\tilde{b}_{\Omega}-W_{\Omega})-\boldsymbol{\Lambda}_{R}||\tilde{R}||_{I}\right)\nonumber \\
& +\mathcal{E}_{P}^{\top}\boldsymbol{\Psi}_{P}\left(\hat{R}(\tilde{b}_{V}-W_{V})+\left[\hat{P}-\tilde{P}\right]_{\times}(\hat{R}\tilde{b}_{\Omega}-W_{\Omega})\right)\nonumber \\
& -\mathcal{E}_{P}^{\top}\boldsymbol{\Psi}_{P}\boldsymbol{\Lambda}_{P}\tilde{P}-\frac{1}{\gamma}\tilde{b}_{\Omega}^{\top}\dot{\hat{b}}_{\Omega}-\frac{1}{\gamma}\tilde{b}_{V}^{\top}\dot{\hat{b}}_{V}\label{eq:SE3PPF_Vdot_Ry}
\end{align}
Consider $||\tilde{R}||_{I}=\frac{1}{4}\frac{||\mathbf{vex}(\boldsymbol{\mathcal{P}}_{a}(\tilde{R}))||^{2}}{1-||\tilde{R}||_{I}}$
as defined in \eqref{eq:SE3PPF_lemm1_1}. Using the result in \eqref{eq:SE3PPF_Vdot_Ry}
and directly substituting $\dot{\hat{b}}_{\Omega}$, $\dot{\hat{b}}_{V}$,
$W_{\Omega}$ and $W_{V}$ with their definitions in \eqref{eq:SE3PPF_b1est_dot_Ty},
\eqref{eq:SE3PPF_b2est_dot_Ty}, \eqref{eq:SE3PPF_W1est_dot_Ty},
and \eqref{eq:SE3PPF_W2est_dot_Ty}, respectively, one obtains
\begin{align}
\dot{V} & =-\frac{1}{4}k_{w}\mathcal{E}_{R}^{2}\boldsymbol{\Psi}_{R}^{2}\frac{||\mathbf{vex}(\boldsymbol{\mathcal{P}}_{a}(\tilde{R}))||^{2}}{1-||\tilde{R}||_{I}}-k_{w}\mathcal{E}_{P}^{\top}\boldsymbol{\Psi}_{P}^{2}\mathcal{E}_{P}\nonumber \\
& =-k_{w}\mathcal{E}_{R}^{2}\boldsymbol{\Psi}_{R}^{2}||\tilde{R}||_{I}-k_{w}\mathcal{E}_{P}^{\top}\boldsymbol{\Psi}_{P}^{2}\mathcal{E}_{P}\label{eq:SE3PPF_Vdot_Ry_Final}
\end{align}
The result obtained in \eqref{eq:SE3PPF_Vdot_Ry_Final} indicates
that $V\left(t\right)\leq V\left(0\right),\forall t\geq0$. Given
that $V\left(t\right)\leq V\left(0\right),\forall t\geq0$, $\tilde{R}\left(0\right)\notin\mathcal{U}$
and $\mathcal{E}\left(0\right)\in\mathbb{R}^{4}$, $\tilde{b}$ remains
bounded, and $\mathcal{E}$ is bounded and well defined for all $t\geq0$.
Consequently, $\tilde{P}$, $||\tilde{R}||_{I}$ and $\mathbf{vex}(\boldsymbol{\mathcal{P}}_{a}(\tilde{R}))$
are bounded, which in turn signifies that $\dot{\tilde{P}}$, $||\dot{\tilde{R}}||_{I}$,
$\dot{\mathcal{E}}_{R}$ and $\dot{\mathcal{E}}_{P}$ are bounded
as well. From the result in \eqref{eq:SE3PPF_Vdot_Ry_Final} it follows
that
\begin{align}
\ddot{V}= & -k_{w}\left(2\mathcal{E}_{R}\boldsymbol{\Psi}_{R}(\dot{\mathcal{E}}_{R}\boldsymbol{\Psi}_{R}+\mathcal{E}_{R}\dot{\boldsymbol{\Psi}}_{R})||\tilde{R}||_{I}+\mathcal{E}_{R}^{2}\boldsymbol{\Psi}_{R}^{2}||\dot{\tilde{R}}||_{I}\right)\nonumber \\
& -2k_{w}\mathcal{E}_{P}^{\top}\boldsymbol{\Psi}_{P}^{2}\dot{\mathcal{E}}_{P}-2k_{w}\mathcal{E}_{P}^{\top}\boldsymbol{\Psi}_{P}\dot{\boldsymbol{\Psi}}_{P}\mathcal{E}_{P}\label{eq:SE3PPF_Vddot_Ry_Final}
\end{align}
Since $\boldsymbol{\Psi}_{R}=\mu_{1}$ and $\boldsymbol{\Psi}_{P}={\rm diag}(\mu_{2},\mu_{3},\mu_{4})$
defined in \eqref{eq:SE3PPF_mu}, $\dot{\mu}_{i}$ can be expressed
as follows for all $i=1,2,\ldots,4$ 
\begin{align}
\dot{\mu}_{i}= & -\frac{1}{2}\frac{\underline{\delta}_{i}\dot{\xi}_{i}+\dot{\boldsymbol{e}}_{i}}{(\underline{\delta}_{i}\xi_{i}+\boldsymbol{e}_{i})^{2}}-\frac{1}{2}\frac{\bar{\delta}_{i}\dot{\xi}_{i}-\dot{\boldsymbol{e}}_{i}}{(\bar{\delta}_{i}\xi_{i}-\boldsymbol{e}_{i})^{2}}\label{eq:SE3PPF_mu_dot}
\end{align}
with $\dot{\xi}_{i}=-\ell_{i}(\xi_{i}^{0}-\xi_{i}^{\infty})\exp(-\ell_{i}t)$.
Due to the fact that $\dot{\boldsymbol{e}}_{i}$ is bounded for all
$i=1,2,\ldots,4$, $\dot{\mu}_{i}$ is bounded and $\ddot{V}$ in
\eqref{eq:SE3PPF_Vddot_Ry_Final} is uniformly bounded for all $t\geq0$.
It should be remarked that $\mathcal{E}_{1}>0$ for all $||\tilde{R}||_{I}>0$,
and $\mathcal{E}_{1}\rightarrow0$ as $||\tilde{R}||_{I}\rightarrow0$
and vice versa as stated in property (ii) of Proposition \ref{Prop:SE3PPF_1}.
In addition, $\mathcal{E}_{i}\neq0\forall\boldsymbol{e}_{i}\neq0$
and $\mathcal{E}_{i}=0$ if and only if $\boldsymbol{e}_{i}=0$ as
indicated in property (iii) of Proposition \ref{Prop:SE3PPF_1}. Therefore,
$\dot{V}$ is uniformly continuous, and in consistence with Barbalat
Lemma, $\dot{V}\rightarrow0$ as $t\rightarrow\infty$ signifies that
$\mathcal{E}_{i}\rightarrow0$ and $\boldsymbol{e}_{i}\rightarrow0$.
As mentioned by property (iv) of Proposition \ref{Prop:SE3PPF_1},
$\mathcal{E}\rightarrow0$ implies that $\tilde{\boldsymbol{T}}$
asymptotically approaches $\mathbf{I}_{4}$ which completes the proof.

\subsection{Direct Pose Filter with Prescribed Performance\label{subsec:SE3PPF_Explicit-Filter}}

The reconstructed homogeneous transformation matrix $\boldsymbol{T}_{y}$
defined in Subsection \ref{subsec:SE3PPF_Passive-Filter} consists
of two elements: $R_{y}$ and $P_{y}$. Although, $R_{y}$ can be
statically reconstructed applying, for example, QUEST \cite{shuster1981three},
or SVD \cite{markley1988attitude}, the aforementioned methods of
static reconstruction could significantly increase processing cost
\cite{hashim2018Conf1,hashim2019SO3Det}.
Thus, the pose filter proposed in this Subsection avoids the necessity
of attitude reconstruction and instead uses measurements from the
inertial and body-frame units directly. Let us define
\begin{align}
\boldsymbol{\mathcal{M}}_{{\rm T}}=\left[\begin{array}{cc}
\mathbf{M}_{{\rm T}} & \mathbf{m}_{{\rm v}}\\
\mathbf{m}_{{\rm v}}^{\top} & \mathbf{m}_{{\rm c}}
\end{array}\right]= & \sum_{i=1}^{N_{{\rm R}}}k_{i}^{{\rm R}}\left[\begin{array}{c}
\upsilon_{i}^{\mathcal{I}\left({\rm R}\right)}\\
0
\end{array}\right]\left[\begin{array}{c}
\upsilon_{i}^{\mathcal{I}\left({\rm R}\right)}\\
0
\end{array}\right]^{\top}\nonumber \\
& +\sum_{j=1}^{N_{{\rm L}}}k_{j}^{{\rm L}}\left[\begin{array}{c}
{\rm v}_{j}^{\mathcal{I}\left({\rm L}\right)}\\
1
\end{array}\right]\left[\begin{array}{c}
{\rm v}_{j}^{\mathcal{I}\left({\rm L}\right)}\\
1
\end{array}\right]^{\top}\label{eq:SE3PPF_MI}
\end{align}
such that $\mathbf{M}_{{\rm T}}=\mathbf{M}_{{\rm R}}+\mathbf{M}_{{\rm L}}$
with
\begin{align}
\mathbf{M}_{{\rm R}}= & \sum_{i=1}^{N_{{\rm R}}}k_{i}^{{\rm R}}\upsilon_{i}^{\mathcal{I}\left({\rm R}\right)}\left(\upsilon_{i}^{\mathcal{I}\left({\rm R}\right)}\right)^{\top}\nonumber \\
\mathbf{M}_{{\rm L}}= & \sum_{j=1}^{N_{{\rm L}}}k_{j}^{{\rm L}}{\rm v}_{j}^{\mathcal{I}\left({\rm L}\right)}\left({\rm v}_{j}^{\mathcal{I}\left({\rm L}\right)}\right)^{\top}\nonumber \\
\mathbf{m}_{{\rm v}}= & \sum_{j=1}^{N_{{\rm L}}}k_{j}^{{\rm L}}{\rm v}_{j}^{\mathcal{I}\left({\rm L}\right)}\nonumber \\
\mathbf{m}_{{\rm c}}= & \sum_{j=1}^{N_{{\rm L}}}k_{j}^{{\rm L}}\label{eq:SE3PPF_MI_elements}
\end{align}
where $k_{i}^{{\rm R}}$ and $k_{j}^{{\rm L}}$ are constant gains
of the confidence level of $i$th and $j$th sensor measurements,
respectively. Define
\begin{align}
\boldsymbol{\mathcal{K}}_{{\rm T}}=\left[\begin{array}{cc}
\mathbf{K}_{{\rm T}} & \mathbf{k}_{{\rm v}}\\
\mathbf{m}_{{\rm v}}^{\top} & \mathbf{m}_{{\rm c}}
\end{array}\right]= & \sum_{i=1}^{N_{{\rm R}}}k_{i}^{{\rm R}}\left[\begin{array}{c}
\upsilon_{i}^{\mathcal{B}\left({\rm R}\right)}\\
0
\end{array}\right]\left[\begin{array}{c}
\upsilon_{i}^{\mathcal{I}\left({\rm R}\right)}\\
0
\end{array}\right]^{\top}\nonumber \\
& +\sum_{j=1}^{N_{{\rm L}}}k_{j}^{{\rm L}}\left[\begin{array}{c}
{\rm v}_{j}^{\mathcal{B}\left({\rm L}\right)}\\
1
\end{array}\right]\left[\begin{array}{c}
{\rm v}_{j}^{\mathcal{I}\left({\rm L}\right)}\\
1
\end{array}\right]^{\top}\label{eq:SE3PPF_KIB}
\end{align}
such that $\mathbf{m}_{{\rm v}}=\sum_{j=1}^{N_{{\rm L}}}k_{j}^{{\rm L}}{\rm v}_{j}^{\mathcal{I}\left({\rm L}\right)}$
and $\mathbf{m}_{{\rm c}}=\sum_{j=1}^{N_{{\rm L}}}k_{j}^{{\rm L}}$
as defined in \eqref{eq:SE3PPF_MI_elements}, and
\begin{align}
\mathbf{K}_{{\rm T}} & =\sum_{i=1}^{N_{{\rm R}}}k_{i}^{{\rm R}}\upsilon_{i}^{\mathcal{B}\left({\rm R}\right)}\left(\upsilon_{i}^{\mathcal{I}\left({\rm R}\right)}\right)^{\top}+\sum_{j=1}^{N_{{\rm L}}}k_{i}^{{\rm L}}{\rm v}_{j}^{\mathcal{B}\left({\rm L}\right)}\left({\rm v}_{j}^{\mathcal{I}\left({\rm L}\right)}\right)^{\top}\nonumber \\
\mathbf{k}_{{\rm v}} & =\sum_{j=1}^{N_{{\rm L}}}k_{j}^{{\rm L}}{\rm v}_{j}^{\mathcal{B}\left({\rm L}\right)}\label{eq:SE3PPF_KIB_elements-1}
\end{align}
In this work $k_{i}^{{\rm R}}$ is selected such that $\sum_{i=1}^{N_{{\rm R}}}k_{i}^{{\rm R}}=3$.
It can be easily deduced that $\mathbf{M}_{{\rm R}}$ is symmetric.
Assuming that Assumption \ref{Assum:SE3STCH_1} holds, $\mathbf{M}_{{\rm R}}$
is nonsingular with ${\rm rank}(\mathbf{M}_{{\rm R}})=3$. Accordingly,
the three eigenvalues of $\mathbf{M}_{{\rm R}}$ are greater than
zero. Define $\bar{\mathbf{M}}_{{\rm R}}={\rm Tr}\{\mathbf{M}_{{\rm R}}\}\mathbf{I}_{3}-\mathbf{M}_{{\rm R}}\in\mathbb{R}^{3\times3}$,
provided that ${\rm rank}(\mathbf{M}_{{\rm R}})=3$, then, the following
three statements hold (\cite{bullo2004geometric} page. 553): 
\begin{enumerate}
	\item $\mathbf{M}_{{\rm R}}$ is a positive-definite matrix.
	\item The eigenvectors of $\mathbf{M}_{{\rm R}}$ coincide with the eigenvectors
	of $\bar{\mathbf{M}}_{{\rm R}}$. 
	\item Assuming that the three eigenvalues of $\mathbf{M}_{{\rm R}}$ are
	$\lambda(\mathbf{M}_{{\rm R}})=\{\lambda_{1},\lambda_{2},\lambda_{3}\}$,
	then $\lambda(\bar{\mathbf{M}}_{{\rm R}})=\{\lambda_{3}+\lambda_{2},\lambda_{3}+\lambda_{1},\lambda_{2}+\lambda_{1}\}$
	with the minimum singular value $\underline{\lambda}(\bar{\mathbf{M}}_{{\rm R}})>0$. 
\end{enumerate}
In the remainder of this Subsection, it is considered that ${\rm rank}(\mathbf{M}_{{\rm R}})=3$
in order to ensure that the above-mentioned statements are true. Define
\begin{equation}
\hat{\upsilon}_{i}^{\mathcal{B}\left({\rm R}\right)}=\hat{R}^{\top}\upsilon_{i}^{\mathcal{I}\left({\rm R}\right)}\label{eq:SE3PPF_vB_hat}
\end{equation}
Defining the error in the homogeneous transformation matrix as in
\eqref{eq:SE3PPF_Terr_matrix}, the attitude error can be expressed
as $\tilde{R}=\hat{R}R^{\top}$ and the position error is defined
by $\tilde{P}=\hat{P}-\tilde{R}P$. Also, let the bias error be as
in \eqref{eq:SE3PPF_b1_tilde} and \eqref{eq:SE3PPF_b2_tilde}. In
order to derive the direct pose filter, it is necessary to introduce
the following series of equations written in terms of vectorial measurements.
According to identity \eqref{eq:SO3PPF_Identity1} and \eqref{eq:SO3PPF_Identity2},
one has
\begin{align*}
& \left[\hat{R}\sum_{i=1}^{N_{{\rm R}}}\frac{k_{i}^{{\rm R}}}{2}\hat{\upsilon}_{i}^{\mathcal{B}\left({\rm R}\right)}\times\upsilon_{i}^{\mathcal{B}\left({\rm R}\right)}\right]_{\times}\\
& \hspace{2em}=\hat{R}\left[\sum_{i=1}^{N_{{\rm R}}}\frac{k_{i}^{{\rm R}}}{2}\hat{\upsilon}_{i}^{\mathcal{B}\left({\rm R}\right)}\times\upsilon_{i}^{\mathcal{B}\left({\rm R}\right)}\right]_{\times}\hat{R}^{\top}\\
& \hspace{2em}=\hat{R}\sum_{i=1}^{N_{{\rm R}}}\frac{k_{i}^{{\rm R}}}{2}\left(\upsilon_{i}^{\mathcal{B}\left({\rm R}\right)}\left(\hat{\upsilon}_{i}^{\mathcal{B}\left({\rm R}\right)}\right)^{\top}-\hat{\upsilon}_{i}^{\mathcal{B}\left({\rm R}\right)}\left(\upsilon_{i}^{\mathcal{B}\left({\rm R}\right)}\right)^{\top}\right)\hat{R}^{\top}\\
& \hspace{2em}=\frac{1}{2}\hat{R}R^{\top}\mathbf{M}_{{\rm R}}-\frac{1}{2}\mathbf{M}_{{\rm R}}R\hat{R}^{\top}\\
& \hspace{2em}=\boldsymbol{\mathcal{P}}_{a}(\tilde{R}\mathbf{M}_{{\rm R}})
\end{align*}
such that
\begin{equation}
\mathbf{vex}(\boldsymbol{\mathcal{P}}_{a}(\tilde{R}\mathbf{M}_{{\rm R}}))=\hat{R}\sum_{i=1}^{N_{{\rm R}}}\left(\frac{k_{i}^{{\rm R}}}{2}\hat{\upsilon}_{i}^{\mathcal{B}\left({\rm R}\right)}\times\upsilon_{i}^{\mathcal{B}\left({\rm R}\right)}\right)\label{eq:SE3PPF_VEX_VM}
\end{equation}
Thus, $\tilde{R}\mathbf{M}_{{\rm R}}$ is defined in terms of vectorial
measurements by
\begin{equation}
\tilde{R}\mathbf{M}_{{\rm R}}=\hat{R}\sum_{i=1}^{N_{{\rm R}}}\left(k_{i}^{{\rm R}}\upsilon_{i}^{\mathcal{B}\left({\rm R}\right)}\left(\upsilon_{i}^{\mathcal{I}\left({\rm R}\right)}\right)^{\top}\right)\label{eq:SE3PPF_RM_VM}
\end{equation}
The normalized Euclidean distance of $\tilde{R}\mathbf{M}_{{\rm R}}$
is found to be
\begin{align}
||\tilde{R}\mathbf{M}_{{\rm R}}||_{I} & =\frac{1}{4}{\rm Tr}\{(\mathbf{I}_{3}-\tilde{R})\mathbf{M}_{{\rm R}}\}\nonumber \\
& =\frac{1}{4}{\rm Tr}\left\{ \mathbf{I}_{3}-\hat{R}\sum_{i=1}^{N_{{\rm R}}}\left(k_{i}^{{\rm R}}\upsilon_{i}^{\mathcal{B}\left({\rm R}\right)}\left(\upsilon_{i}^{\mathcal{I}\left({\rm R}\right)}\right)^{\top}\right)\right\} \nonumber \\
& =\frac{1}{4}\sum_{i=1}^{N_{{\rm R}}}\left(1-\left(\hat{\upsilon}_{i}^{\mathcal{B}\left({\rm R}\right)}\right)^{\top}\upsilon_{i}^{\mathcal{B}\left({\rm R}\right)}\right)\label{eq:SE3PPF_RI_VM}
\end{align}
Let us introduce the following variable
\begin{align}
\boldsymbol{\Upsilon}(\mathbf{M}_{{\rm R}},\tilde{R})= & {\rm Tr}\left\{ \tilde{R}\mathbf{M}_{{\rm R}}\mathbf{M}_{{\rm R}}^{-1}\right\} \nonumber \\
= & {\rm Tr}\left\{ \left(\sum_{i=1}^{N_{{\rm R}}}k_{i}^{{\rm R}}\upsilon_{i}^{\mathcal{B}\left({\rm R}\right)}\left(\upsilon_{i}^{\mathcal{I}\left({\rm R}\right)}\right)^{\top}\right)\right.\nonumber \\
& \left.\hspace{1em}\bullet\left(\sum_{i=1}^{N_{{\rm R}}}k_{i}^{{\rm R}}\hat{\upsilon}_{i}^{\mathcal{B}\left({\rm R}\right)}\left(\upsilon_{i}^{\mathcal{I}\left({\rm R}\right)}\right)^{\top}\right)^{-1}\right\} \label{eq:SE3PPF_Gamma_VM}
\end{align}
where $\bullet$ is a multiplication operator of the two matrices. From \eqref{eq:SE3PPF_MI}
and \eqref{eq:SE3PPF_MI_elements}, one obtains
\begin{align}
\tilde{\boldsymbol{T}}\boldsymbol{\mathcal{M}}^{\mathcal{I}} & =\left[\begin{array}{cc}
\tilde{R}\mathbf{M}_{{\rm T}}+\tilde{P}\mathbf{m}_{{\rm v}}^{\top} & \tilde{R}\mathbf{m}_{{\rm v}}+\mathbf{m}_{{\rm c}}\tilde{P}\\
\mathbf{m}_{{\rm v}}^{\top} & \mathbf{m}_{{\rm c}}
\end{array}\right]\label{eq:SE3PPF_VM_Part1}
\end{align}
The above-mentioned result can be additionally expressed as
\begin{align}
\tilde{\boldsymbol{T}}\boldsymbol{\mathcal{M}}^{\mathcal{I}} & =\left[\begin{array}{cc}
\hat{R} & \hat{P}\\
\underline{\mathbf{0}}_{3}^{\top} & 1
\end{array}\right]\left[\begin{array}{cc}
\mathbf{K}_{{\rm T}} & \mathbf{k}_{{\rm v}}\\
\mathbf{m}_{{\rm v}}^{\top} & \mathbf{m}_{{\rm c}}
\end{array}\right]\nonumber \\
& =\left[\begin{array}{cc}
\hat{R}\mathbf{K}_{{\rm T}}+\hat{P}\mathbf{m}_{{\rm v}}^{\top} & \hat{R}\mathbf{k}_{{\rm v}}+\mathbf{m}_{{\rm c}}\hat{P}\\
\mathbf{m}_{{\rm v}}^{\top} & \mathbf{m}_{{\rm c}}
\end{array}\right]\label{eq:SE3PPF_VM_Part2}
\end{align}
As such, from \eqref{eq:SE3PPF_VM_Part1} and \eqref{eq:SE3PPF_VM_Part2},
the position error can be reformulated with respect to vectorial measurements
as
\begin{equation}
\tilde{P}=\hat{P}+\frac{1}{\mathbf{m}_{{\rm c}}}\left(\hat{R}\mathbf{k}_{{\rm v}}-\tilde{R}\mathbf{M}_{{\rm R}}\mathbf{M}_{{\rm R}}^{-1}\mathbf{m}_{{\rm v}}\right)\label{eq:SE3PPF_Ptil_VM}
\end{equation}
with $\tilde{R}\mathbf{M}_{{\rm R}}$ being calculated as in \eqref{eq:SE3PPF_RM_VM} and $\mathbf{m}_{{\rm c}}\neq 0$ for at least one landmark.
Consequently, $\mathbf{vex}(\boldsymbol{\mathcal{P}}_{a}(\tilde{R}\mathbf{M}_{{\rm R}}))$,
$\tilde{R}\mathbf{M}_{{\rm R}}$, $||\tilde{R}\mathbf{M}_{{\rm R}}||_{I}$,
$\boldsymbol{\Upsilon}(\mathbf{M}_{{\rm R}},\tilde{R})$, and $\tilde{P}$
will be obtained through a set of vectorial measurements as defined
in \eqref{eq:SE3PPF_VEX_VM}, \eqref{eq:SE3PPF_RM_VM}, \eqref{eq:SE3PPF_RI_VM},
\eqref{eq:SE3PPF_Gamma_VM}, and \eqref{eq:SE3PPF_Ptil_VM}, respectively,
in all the subsequent derivations and calculations. Let us modify
the vector error in \eqref{eq:SE3PPF_Vec_error} to be 
\begin{equation}
\boldsymbol{e}=\left[\boldsymbol{e}_{1},\boldsymbol{e}_{2},\boldsymbol{e}_{3},\boldsymbol{e}_{4}\right]^{\top}=\left[||\tilde{R}\mathbf{M}_{{\rm R}}||_{I},\tilde{P}^{\top}\right]^{\top}\label{eq:SE3PPF_e_VM}
\end{equation}
with $||\tilde{R}\mathbf{M}_{{\rm R}}||_{I}$ and $\tilde{P}$ being
defined in \eqref{eq:SE3PPF_RI_VM} and \eqref{eq:SE3PPF_Ptil_VM},
respectively. Thus, all the discussion in Subsection \ref{subsec:SE3PPF_Prescribed-Performance}
is to be reformulated using the error vector in \eqref{eq:SE3PPF_e_VM}
instead of \eqref{eq:SE3PPF_Vec_error}. Define the minimum eigenvalue
of $\bar{\mathbf{M}}_{{\rm R}}$ as $\underline{\lambda}:=\underline{\lambda}(\bar{\mathbf{M}}_{{\rm R}})$,
and consider the following filter design
\begin{align}
\dot{\hat{R}}= & \hat{R}\left[\Omega_{m}-\hat{b}_{\Omega}-\hat{R}^{\top}W_{\Omega}\right]_{\times}\label{eq:SE3PPF_Rest_dot_T_VM}\\
\dot{\hat{P}}= & \hat{R}(V_{m}-\hat{b}_{V}-W_{V})\label{eq:SE3PPF_Pest_dot_T_VM}\\
\dot{\hat{b}}_{\Omega}= & \frac{\gamma}{2}\boldsymbol{\Psi}_{R}\mathcal{E}_{R}\hat{R}^{\top}\mathbf{vex}(\boldsymbol{\mathcal{P}}_{a}(\tilde{R}\mathbf{M}_{{\rm R}}))\nonumber \\
& +\gamma\hat{R}^{\top}\left[\tilde{P}-\hat{P}\right]_{\times}\boldsymbol{\Psi}_{P}\mathcal{E}_{P}\label{eq:SE3PPF_b1est_dot_T_VM}\\
\dot{\hat{b}}_{V}= & \gamma\hat{R}^{\top}\boldsymbol{\Psi}_{P}\mathcal{E}_{P}\label{eq:SE3PPF_b2est_dot_T_VM}\\
W_{\Omega}= & \frac{4}{\underline{\lambda}}\frac{k_{w}\boldsymbol{\Psi}_{R}\mathcal{E}_{R}-\boldsymbol{\Lambda}_{R}}{1+\boldsymbol{\Upsilon}(\mathbf{M}_{{\rm R}},\tilde{R})}\mathbf{vex}(\boldsymbol{\mathcal{P}}_{a}(\tilde{R}\mathbf{M}_{{\rm R}}))\label{eq:SE3PPF_W1est_dot_T_VM}\\
W_{V}= & \hat{R}^{\top}\left(k_{w}\boldsymbol{\Psi}_{P}\mathcal{E}_{P}+\left[\tilde{P}-\hat{P}\right]_{\times}W_{\Omega}-\boldsymbol{\Lambda}_{P}\tilde{P}\right)\label{eq:SE3PPF_W2est_dot_T_VM}
\end{align}
with $\boldsymbol{\Upsilon}(\mathbf{M}_{{\rm R}},\tilde{R})$ and
$\mathbf{vex}(\boldsymbol{\mathcal{P}}_{a}(\tilde{R}\mathbf{M}_{{\rm R}}))$
being specified in \eqref{eq:SE3PPF_Gamma_VM} and \eqref{eq:SE3PPF_VEX_VM},
respectively, $\mathcal{E}=[\mathcal{E}_{R},\mathcal{E}_{P}^{\top}]^{\top}=[\mathcal{E}_{1},\mathcal{E}_{2},\mathcal{E}_{3},\mathcal{E}_{4}]^{\top}$,
 $\mathcal{E}_{i}:=\mathcal{E}_{i}(\boldsymbol{e}_{i},\xi_{i})$ and
$\mu_{i}:=\mu_{i}(\boldsymbol{e}_{i},\xi_{i})$ being defined in \eqref{eq:SE3PPF_trans3}
and \eqref{eq:SE3PPF_mu}, respectively, while $\boldsymbol{e}$ is
as in \eqref{eq:SE3PPF_e_VM}, $k_{w}$ and $\gamma$ are positive
constants, and $\hat{b}_{\Omega}$ and $\hat{b}_{V}$ are the estimates
of $b_{\Omega}$ and $b_{V}$, respectively. The equivalent quaternion representation and complete implementation
steps of the direct filter are given in \nameref{sec:SO3_PPF_STCH_AppendixB}.
\begin{thm}
	\textbf{\label{thm:SE3PPF_2}} Consider coupling the pose filter in
	\eqref{eq:SE3PPF_Rest_dot_T_VM}, \eqref{eq:SE3PPF_Pest_dot_T_VM},
	\eqref{eq:SE3PPF_b1est_dot_T_VM}, \eqref{eq:SE3PPF_b2est_dot_T_VM},
	\eqref{eq:SE3PPF_W1est_dot_T_VM}, and \eqref{eq:SE3PPF_W2est_dot_T_VM}
	with the set of vector measurements in \eqref{eq:SE3STCH_Set_R_Norm}
	and \eqref{eq:SE3STCH_Set_L}, %
	{} and the velocity measurements in \eqref{eq:SE3PPF_Angular} and \eqref{eq:SE3PPF_V_Trans}
	where $\Omega_{m}=\Omega+b_{\Omega}$ and $V_{m}=V+b_{V}$. Let Assumption
	\ref{Assum:SE3STCH_1} hold. Define $\mathcal{U}\subseteq\mathbb{SE}\left(3\right)\times\mathbb{R}^{6}$
	by{\small{} $\mathcal{U}:=\left\{ \left.(\tilde{\boldsymbol{T}}\left(0\right),\tilde{b}\left(0\right))\right|{\rm Tr}\{\tilde{R}\left(0\right)\}=-1,\tilde{P}\left(0\right)=\underline{\mathbf{0}}_{3},\tilde{b}\left(0\right)=\underline{\mathbf{0}}_{6}\right\} $}.
	If $\tilde{R}\left(0\right)\notin\mathcal{U}$ and $\mathcal{E}\left(0\right)\in\mathcal{L}_{\infty}$,
	then, all error signals are bounded, $\mathcal{E}\left(t\right)$
	asymptotically approaches $0$, and $\tilde{\boldsymbol{T}}$ asymptotically
	approaches $\mathbf{I}_{4}$. 
\end{thm}
Theorem \ref{thm:SE3PPF_2} guarantees the observer dynamics in \eqref{eq:SE3PPF_Rest_dot_T_VM},
\eqref{eq:SE3PPF_Pest_dot_T_VM}, \eqref{eq:SE3PPF_b1est_dot_T_VM},
\eqref{eq:SE3PPF_b2est_dot_T_VM}, \eqref{eq:SE3PPF_W1est_dot_T_VM},
and \eqref{eq:SE3PPF_W2est_dot_T_VM} to be stable. In consistence
with Remark \ref{rem:SO3PPF_1} boundedness of $\mathcal{E}\left(t\right)$
indicates that $\boldsymbol{e}$ follows the dynamic decreasing boundaries
in \eqref{eq:SE3PPF_Presc}.

\textbf{Proof. }Consider the error in the homogeneous transformation
matrix and bias defined as in \eqref{eq:SE3PPF_Terr_matrix}, \eqref{eq:SE3PPF_b1_tilde}
and \eqref{eq:SE3PPF_b2_tilde}, respectively. From \eqref{eq:SE3PPF_R_Dynamics}
and \eqref{eq:SE3PPF_Rest_dot_T_VM}, the error dynamics of $\tilde{R}$
can be found to be analogous to \eqref{eq:SE3PPF_Rtilde_dot}. The
$i$th inertial measurements ${\rm v}_{i}^{\mathcal{I}\left({\rm R}\right)}$
and ${\rm v}_{i}^{\mathcal{I}\left({\rm L}\right)}$ are constant,
thus, $\dot{\mathbf{M}}_{{\rm R}}=\mathbf{0}_{3\times3}$. Consequently,
from \eqref{eq:SE3PPF_Rtilde_dot}, the derivative of $||\tilde{R}\mathbf{M}_{{\rm R}}||_{I}$
is equivalent to
\begin{align}
\frac{d}{dt}||\tilde{R}\mathbf{M}_{{\rm R}}||_{I}= & -\frac{1}{4}{\rm Tr}\left\{ \left[\hat{R}\tilde{b}_{\Omega}-W_{\Omega}\right]_{\times}\tilde{R}\mathbf{M}_{{\rm R}}\right\} \nonumber \\
= & -\frac{1}{4}{\rm Tr}\left\{ \left[\hat{R}\tilde{b}_{\Omega}-W_{\Omega}\right]_{\times}\boldsymbol{\mathcal{P}}_{a}(\tilde{R}\mathbf{M}_{{\rm R}})\right\} \nonumber \\
= & \frac{1}{2}\mathbf{vex}(\boldsymbol{\mathcal{P}}_{a}(\tilde{R}\mathbf{M}_{{\rm R}}))^{\top}(\hat{R}\tilde{b}_{\Omega}-W_{\Omega})\label{eq:SE3PPF_NormMIRtilde_dot}
\end{align}
where ${\rm Tr}\left\{ \left[W_{\Omega}\right]_{\times}\tilde{R}\mathbf{M}_{{\rm R}}\right\} =-2\mathbf{vex}(\boldsymbol{\mathcal{P}}_{a}(\tilde{R}\mathbf{M}_{{\rm R}}))^{\top}W_{\Omega}$
as given in \eqref{eq:SO3PPF_Identity7}. One could find that the
derivative of $\tilde{P}$ is equivalent to \eqref{eq:SE3PPF_Ptilde_dot}.
From \eqref{eq:SE3PPF_NormMIRtilde_dot} and \eqref{eq:SE3PPF_Ptilde_dot},
and in view of \eqref{eq:SE3PPF_T_VEC_Dyn}, the derivative of $\boldsymbol{e}$
given in \eqref{eq:SE3PPF_e_VM}, becomes
\begin{align}
\dot{\boldsymbol{e}} & =\left[\begin{array}{cc}
\frac{1}{2}\mathbf{vex}(\boldsymbol{\mathcal{P}}_{a}(\tilde{R}\mathbf{M}_{{\rm R}}))^{\top} & 0_{1\times3}\\
\left[\hat{P}-\tilde{P}\right]_{\times} & \hat{R}
\end{array}\right]\left[\begin{array}{c}
\hat{R}\tilde{b}_{\Omega}-W_{\Omega}\\
\tilde{b}_{V}-W_{V}
\end{array}\right]\label{eq:SE3PPF_VEC_MI_tilde_dot}
\end{align}
The derivative of the transformed error in \eqref{eq:SE3PPF_Trans_dot}
be acquired by direct substitution of $\boldsymbol{e}$ as in \eqref{eq:SE3PPF_e_VM},
in addition to the result in \eqref{eq:SE3PPF_VEC_MI_tilde_dot}.
Consider the candidate Lyapunov function 

\begin{align}
V(\mathcal{E},\tilde{b}_{\Omega},\tilde{b}_{V}) & =\frac{1}{2}||\mathcal{E}||^{2}+\frac{1}{2\gamma}||\tilde{b}_{\Omega}||^{2}+\frac{1}{2\gamma}||\tilde{b}_{V}||^{2}\label{eq:SE3PPF_V_VM}
\end{align}
The derivative of $V:=V(\mathcal{E},\tilde{b}_{\Omega},\tilde{b}_{V})$
is as follows
\begin{align}
\dot{V}= & \mathcal{E}^{\top}\dot{\mathcal{E}}-\frac{1}{\gamma}\tilde{b}_{\Omega}^{\top}\dot{\hat{b}}_{\Omega}-\frac{1}{\gamma}\tilde{b}_{V}^{\top}\dot{\hat{b}}_{V}\nonumber \\
= & \frac{1}{2}\mathcal{E}_{R}\boldsymbol{\Psi}_{R}\mathbf{vex}(\boldsymbol{\mathcal{P}}_{a}(\tilde{R}\mathbf{M}_{{\rm R}}))^{\top}(\hat{R}\tilde{b}_{\Omega}-W_{\Omega})\nonumber \\
& +\mathcal{E}_{P}^{\top}\boldsymbol{\Psi}_{P}\left(\hat{R}(\tilde{b}_{V}-W_{V})+\left[\hat{P}-\tilde{P}\right]_{\times}(\hat{R}\tilde{b}_{\Omega}-W_{\Omega})\right)\nonumber \\
& -\mathcal{E}_{R}\boldsymbol{\Psi}_{R}\boldsymbol{\Lambda}_{R}||\tilde{R}\mathbf{M}_{{\rm R}}||_{I}-\mathcal{E}_{P}^{\top}\boldsymbol{\Psi}_{P}\boldsymbol{\Lambda}_{P}\tilde{P}\nonumber \\
& -\frac{1}{\gamma}\tilde{b}_{\Omega}^{\top}\dot{\hat{b}}_{\Omega}-\frac{1}{\gamma}\tilde{b}_{V}^{\top}\dot{\hat{b}}_{V}\label{eq:SE3PPF_Vdot_VM-1}
\end{align}
Directly substituting for $\dot{\hat{b}}_{\Omega}$, $\dot{\hat{b}}_{V}$,
$W_{\Omega}$ and $W_{V}$ in \eqref{eq:SE3PPF_b1est_dot_T_VM}, \eqref{eq:SE3PPF_b2est_dot_T_VM},
\eqref{eq:SE3PPF_W1est_dot_T_VM}, and \eqref{eq:SE3PPF_W2est_dot_T_VM},
respectively, results in
\begin{align}
\dot{V}\leq & \boldsymbol{\Lambda}_{R}\left(\frac{2}{\underline{\lambda}}\frac{||\mathbf{vex}(\boldsymbol{\mathcal{P}}_{a}(\tilde{R}\mathbf{M}_{{\rm R}}))||^{2}}{1+\boldsymbol{\Upsilon}(\mathbf{M}_{{\rm R}},\tilde{R})}-||\tilde{R}\mathbf{M}_{{\rm R}}||_{I}\right)\mathcal{E}_{R}\boldsymbol{\Psi}_{R}\nonumber \\
& -\frac{2}{\underline{\lambda}}\frac{k_{w}\mathcal{E}_{R}^{2}\boldsymbol{\Psi}_{R}^{2}}{1+\boldsymbol{\Upsilon}(\mathbf{M}_{{\rm R}},\tilde{R})}\left\Vert \mathbf{vex}(\boldsymbol{\mathcal{P}}_{a}(\tilde{R}\mathbf{M}_{{\rm R}}))\right\Vert ^{2}\nonumber \\
& -k_{w}\mathcal{E}_{P}^{\top}\boldsymbol{\Psi}_{P}^{2}\mathcal{E}_{P}\label{eq:SE3PPF_Vdot_VM-2}
\end{align}
It can be easily found that{\small{}
	\begin{equation}
	\boldsymbol{\Lambda}_{R}\left(\frac{2}{\underline{\lambda}}\frac{\left\Vert \mathbf{vex}(\boldsymbol{\mathcal{P}}_{a}(\tilde{R}\mathbf{M}_{{\rm R}}))\right\Vert ^{2}}{1+\boldsymbol{\Upsilon}(\mathbf{M}_{{\rm R}},\tilde{R})}-||\tilde{R}\mathbf{M}_{{\rm R}}||_{I}\right)\mathcal{E}_{R}\boldsymbol{\Psi}_{R}\leq0\label{eq:SE3PPF_Factor}
	\end{equation}
}where $\mathcal{E}_{R}>0\forall||\tilde{R}\mathbf{M}_{{\rm R}}||_{I}\neq0$
and $\mathcal{E}_{R}=0$ at $||\tilde{R}\mathbf{M}_{{\rm R}}||_{I}=0$
as presented in (ii) Proposition \ref{Prop:SE3PPF_1}, and $\boldsymbol{\Psi}_{R}>0\forall t\geq0$
as given in \eqref{eq:SE3PPF_mu}. Also, $\dot{\xi}_{i}$ is negative
and strictly increasing that satisfies $\dot{\xi}_{i}\rightarrow0$
as $t\rightarrow\infty$, and $\xi_{i}:\mathbb{R}_{+}\to\mathbb{R}_{+}$
such that $\xi_{i}\rightarrow\xi_{i}^{\infty}$ as $t\rightarrow\infty$.
Thus, $\dot{\xi}_{i}/\xi_{i}\leq0$ which means that $\boldsymbol{\Lambda}_{R}\leq0$.
Considering \eqref{eq:SE3PPF_lemm1_2} in Lemma \ref{Lemm:SE3PPF_1},
thus, the expression in \eqref{eq:SE3PPF_Factor} is negative semi-definite.
As such, the inequality in \eqref{eq:SE3PPF_Vdot_VM-2} can be expressed
as
\begin{align}
\dot{V}\leq & -k_{w}\mathcal{E}_{R}^{2}\boldsymbol{\Psi}_{R}^{2}||\tilde{R}\mathbf{M}_{{\rm R}}||_{I}-k_{w}\mathcal{E}_{P}^{\top}\boldsymbol{\Psi}_{P}^{2}\mathcal{E}_{P}\label{eq:SO3PPF_Vdot_VM_Final}
\end{align}
This signifies that $V\left(t\right)\leq V\left(0\right),\forall t\geq0$.
From almost any initial conditions such that ${\rm Tr}\left\{ \tilde{R}\left(0\right)\right\} \neq-1$
and $\mathcal{E}\left(0\right)\in\mathbb{R}^{4}$, $\mathcal{E}$
and $\tilde{b}$ are bounded for all $t\geq0$. Thereby, $\mathcal{E}$
is bounded and well-defined for all $t\geq0$. $\tilde{P}$, $||\tilde{R}\mathbf{M}_{{\rm R}}||_{I}$,
and $\mathbf{vex}(\boldsymbol{\mathcal{P}}_{a}(\tilde{R}\mathbf{M}_{{\rm R}}))$
are also bounded which indicates that $\dot{\tilde{P}}$, $||\dot{\tilde{R}}\mathbf{M}_{{\rm R}}||_{I}$,
$\dot{\mathcal{E}}_{R}$ and $\dot{\mathcal{E}}_{P}$ are bounded
as well. In order to prove asymptotic convergence of $\mathcal{E}$
to the origin and $\tilde{\boldsymbol{T}}$ to the identity, it is
necessary to show that the second derivative of \eqref{eq:SE3PPF_V_VM}
is
\begin{align}
\ddot{V}\leq & -2k_{w}\mathcal{E}_{R}\boldsymbol{\Psi}_{R}(\dot{\mathcal{E}}_{R}\boldsymbol{\Psi}_{R}+\mathcal{E}_{R}\dot{\boldsymbol{\Psi}}_{R})||\tilde{R}\mathbf{M}_{{\rm R}}||_{I}\nonumber \\
& -k_{w}\mathcal{E}_{R}^{2}\boldsymbol{\Psi}_{R}^{2}||\dot{\tilde{R}}\mathbf{M}_{{\rm R}}||_{I}\nonumber \\
& -2k_{w}\mathcal{E}_{P}^{\top}\boldsymbol{\Psi}_{P}(\boldsymbol{\Psi}_{P}\dot{\mathcal{E}}_{P}+\dot{\boldsymbol{\Psi}}_{P}\mathcal{E}_{P})\label{eq:SO3PPF_Vddot_VM_Final}
\end{align}
Recall that $\boldsymbol{\Psi}_{R}=\mu_{1}$ and $\boldsymbol{\Psi}_{P}={\rm diag}(\mu_{2},\mu_{3},\mu_{4})$,
where $\dot{\mu}_{i}$ was defined in \eqref{eq:SE3PPF_mu_dot} for
all $i=1,2,\ldots,4$. Since $\dot{\boldsymbol{e}}_{i}$ is bounded,
$\dot{\mu}_{i}$ is bounded as well and $\ddot{V}$ in \eqref{eq:SO3PPF_Vddot_VM_Final}
is bounded for all $t\geq0$. From property (ii) of Proposition \ref{Prop:SE3PPF_1},
$||\mathcal{E}_{1}||\rightarrow0$ indicates that $||\tilde{R}\mathbf{M}_{{\rm R}}||_{I}\rightarrow0$,
while $\mathcal{E}_{1}\neq0\forall||\tilde{R}\mathbf{M}_{{\rm R}}||_{I}\neq0$
and according to property (iii) of Proposition \ref{Prop:SE3PPF_1},
$\mathcal{E}_{i}\neq0\forall\boldsymbol{e}_{i}\neq0$ and $\mathcal{E}_{i}=0$
if and only if $\boldsymbol{e}_{i}=0$ for all $i=1,\ldots,4$. Therefore,
$\dot{V}$ is uniformly continuous, and on the basis of Barbalat Lemma,
$\dot{V}\rightarrow0$ implies that $||\mathcal{E}||\rightarrow0$
and $\left\Vert \boldsymbol{e}\right\Vert \rightarrow0$ as $t\rightarrow\infty$.
This means that $\tilde{\boldsymbol{T}}$ approaches $\mathbf{I}_{4}$
asymptotically in accordance with (iv) of Proposition \ref{Prop:SE3PPF_1},
which completes the proof.

The estimates $\dot{\hat{b}}_{\Omega}$ and $\dot{\hat{b}}_{V}$ and
the correction factors $W_{\Omega}$ and $W_{V}$ are functions of
the transformed error $\mathcal{E}$ and the auxiliary component $\mu$.
$\mathcal{E}$ and $\mu$ rely on the error $\boldsymbol{e}$ such
that their values become increasingly aggressive as $||\tilde{R}||_{I}$
approaches the unstable equilibria $||\tilde{R}||_{I}\rightarrow+1$
and $\tilde{P}\rightarrow\infty$. Their dynamic behavior is essential
for forcing the proposed filters to obey the prescribed performance
constraints. On the other side $\mathcal{E}\rightarrow0$ as $\boldsymbol{e}\rightarrow0$.
This significant advantage was not offered in literature, such as
\cite{rehbinder2003pose,baldwin2007complementary,baldwin2009nonlinear,hua2011observer,vasconcelos2010nonlinear,hashim2018SE3Stochastic}. 
\begin{rem}
	\label{rem:SO3PPF_3}\textbf{(Design parameters)} The dynamic boundaries
	of $\boldsymbol{e}$ are described by $\bar{\delta}$, $\underline{\delta}$,
	$\xi_{\infty}$, and $\xi_{0}$ where $\xi_{0}$ and $\xi_{\infty}$
	define the large and small sets, respectively. The rate of convergence
	from the given large set to the small set is controlled by $\ell$.
	The initial value of $\boldsymbol{e}\left(0\right)$ in \eqref{eq:SE3PPF_Vec_error}
	or \eqref{eq:SE3PPF_e_VM} can be easily obtained. When applying semi-direct
	pose filter, $R_{y}\left(0\right)$ can be reconstructed, for example,
	using \cite{shuster1981three,markley1988attitude}%
	, $P_{y}\left(0\right)$ can be evaluated by $P_{y}\left(0\right)=\mathcal{G}_{c}^{\mathcal{I}}-R_{y}\left(0\right)\mathcal{G}_{c}^{\mathcal{B}}$
	as in \eqref{eq:SE3PPF_Py}, and finally $||\tilde{R}\left(0\right)||_{I}=\frac{1}{4}{\rm Tr}\{\mathbf{I}_{3}-\hat{R}\left(0\right)R_{y}^{\top}\left(0\right)\}$
	and $\tilde{P}\left(0\right)=\hat{P}\left(0\right)-\tilde{R}\left(0\right)P_{y}\left(0\right)$.
	In case when the direct pose filter is used, $||\tilde{R}\left(0\right)\mathbf{M}_{{\rm R}}||_{I}$
	can be defined from \eqref{eq:SE3PPF_RI_VM} and $\tilde{P}\left(0\right)$
	can be easily obtained in the form of a vectorial measurement based
	on \eqref{eq:SE3PPF_Ptil_VM}. Next, the user can select $\bar{\delta}$,
	$\underline{\delta}$, and $\xi_{0}$ to be greater than $\boldsymbol{e}\left(0\right)$.
\end{rem}

\subsection{Simplified steps of the proposed pose filters}

\textit{The implementation of the proposed nonlinear pose filters
	on $\mathbb{SE}\left(3\right)$ with prescribed performance given
	in Subsections \ref{subsec:SE3PPF_Passive-Filter} and \ref{subsec:SE3PPF_Explicit-Filter}
	can be summarized in the following 7 simplified steps:}

\textbf{\textit{Step 1}}\textit{: Select $\gamma,k_{w}>0$, $\bar{\delta}=\underline{\delta}>\boldsymbol{e}\left(0\right)$,
	the desired speed of the convergence rate $\ell$, and the upper bound
	of the small set $\xi_{\infty}$. }

\textbf{\textit{Step 2}}\textit{: For the case of the semi-direct
	pose filter, define $\boldsymbol{e}=\left[||\tilde{R}||_{I},\tilde{P}^{\top}\right]^{\top}$
	with $\tilde{R}=\hat{R}R_{y}^{\top}$ and $\tilde{P}=\hat{P}-\tilde{R}P_{y}$
	where $P_{y}$ is given in \eqref{eq:SE3PPF_Py} and $R_{y}$ is reconstructed
	(for example \cite{shuster1981three,markley1988attitude}}%
\textit{). For the case of the direct pose filter, define $\boldsymbol{e}=\left[||\tilde{R}\mathbf{M}_{{\rm R}}||_{I},\tilde{P}^{\top}\right]^{\top}$
	with $||\tilde{R}\mathbf{M}_{{\rm R}}||_{I}$ and $\tilde{P}$ being
	specified as in \eqref{eq:SE3PPF_RI_VM} and \eqref{eq:SE3PPF_Ptil_VM},
	respectively. }

\textbf{\textit{Step 3}}\textit{: For the case of the semi-direct
	pose filter, evaluate $\mathbf{vex}(\boldsymbol{\mathcal{P}}_{a}(\tilde{R}))$,
	whereas, for the case of the direct pose filter, define $\mathbf{vex}(\boldsymbol{\mathcal{P}}_{a}(M^{\mathcal{B}}\tilde{R}))$
	and $\boldsymbol{\Upsilon}(\mathbf{M}_{{\rm R}},\tilde{R})$ from
	\eqref{eq:SE3PPF_VEX_VM}, and \eqref{eq:SE3PPF_Gamma_VM}, respectively.
}

\textbf{\textit{Step 4}}\textit{: Find the PPF $\xi$ from \eqref{eq:SE3PPF_Presc}.
}

\textbf{\textit{Step 5}}\textit{: Evaluate the transformed error $\mathcal{E}$,
	$\boldsymbol{\Lambda}_{R},$ $\boldsymbol{\Psi}_{R}$, $\boldsymbol{\Lambda}_{P}$,
	and $\boldsymbol{\Psi}_{P}$ from \eqref{eq:SE3PPF_trans3} and \eqref{eq:SE3PPF_mu},
	respectively.  }

\textbf{\textit{Step 6}}\textit{: Obtain the filter kinematics $\dot{\hat{R}}$,
	$\dot{\hat{P}}$, $\dot{\hat{b}}_{\Omega}$, $\dot{\hat{b}}_{V}$,
	$W_{\Omega}$, and $W_{V}$ from \eqref{eq:SE3PPF_Rest_dot_Ty}, \eqref{eq:SE3PPF_Pest_dot_Ty},
	\eqref{eq:SE3PPF_b1est_dot_Ty}, \eqref{eq:SE3PPF_b2est_dot_Ty},
	\eqref{eq:SE3PPF_W1est_dot_Ty}, and \eqref{eq:SE3PPF_W2est_dot_Ty},
	respectively, for the semi-direct pose filter, or from \eqref{eq:SE3PPF_Rest_dot_T_VM},
	\eqref{eq:SE3PPF_Pest_dot_T_VM}, \eqref{eq:SE3PPF_b1est_dot_T_VM},
	\eqref{eq:SE3PPF_b2est_dot_T_VM}, \eqref{eq:SE3PPF_W1est_dot_T_VM},
	and \eqref{eq:SE3PPF_W2est_dot_T_VM}, respectively, for the direct
	pose filter.  }

\textbf{\textit{Step 7}}\textit{: Go to }\textbf{\textit{Step 2}}\textit{.
} 

\section{Simulations \label{sec:SO3PPF_Simulations}}

This section illustrates the robustness of the proposed pose filters
on $\mathbb{SE}\left(3\right)$ with prescribed performance against
large error in initialization of $\tilde{\boldsymbol{T}}\left(0\right)$
and high levels of bias and noise inherent to the measurement process.
Let the dynamics of the homogeneous transformation matrix $\boldsymbol{T}$
follow \eqref{eq:SE3PPF_T_Dynamics}. Define the true angular velocity
$\left({\rm rad/sec}\right)$ by
\[
\Omega=\left[{\rm sin}\left(0.5t\right),0.7{\rm sin}\left(0.4t+\pi\right),0.5{\rm sin}\left(0.35t+\frac{\pi}{3}\right)\right]^{\top}
\]
with $R\left(0\right)=\mathbf{I}_{3}$. Consider the following true
translational velocity $\left({\rm m/sec}\right)$ 
\[
V=\left[0.3{\rm sin}\left(0.6t\right),0.18{\rm sin}\left(0.4t+\frac{\pi}{2}\right),0.3{\rm sin}\left(0.1t+\frac{\pi}{4}\right)\right]^{\top}
\]
and the initial position $P\left(0\right)=\underline{\mathbf{0}}_{3}$.
Let the measurements of angular and translational velocities be $\Omega_{m}=\Omega+b_{\Omega}+\omega_{\Omega}$
and $V_{m}=V+b_{V}+\omega_{V}$, respectively, with $b_{\Omega}=0.1\left[1,-1,1\right]^{\top}$
and $b_{V}=0.1\left[2,5,1\right]^{\top}$. $\omega_{\Omega}$ and
$\omega_{V}$ represent random noise process at each time instant
with zero mean and standard deviation (STD) equal to $0.15\left({\rm rad/sec}\right)$
and $0.3\left({\rm m/sec}\right)$, respectively. Assume that one
landmark is available for measurement $\left(N_{{\rm L}}=1\right)$
\[
{\rm v}_{1}^{\mathcal{I}\left({\rm L}\right)}=\left[\frac{1}{2},\sqrt{2},1\right]^{\top}
\]
where the body-frame measurements are defined as \eqref{eq:SE3STCH_Vec_Landmark}
such that ${\rm v}_{1}^{\mathcal{B}\left({\rm L}\right)}=R^{\top}\left({\rm v}_{1}^{\mathcal{I}\left({\rm L}\right)}-P\right)+{\rm b}_{1}^{\mathcal{B}\left({\rm L}\right)}+\omega_{1}^{\mathcal{B}\left({\rm L}\right)}$.
The bias vector is %
${\rm b}_{1}^{\mathcal{B}\left({\rm L}\right)}=0.1\left[0.3,0.2,-0.2\right]^{\top}$
while $\omega_{1}^{\mathcal{B}\left({\rm L}\right)}$ is a Gaussian
noise vector with zero mean and ${\rm STD}=0.1$. Assume that two
non-collinear inertial-frame vectors $\left(N_{{\rm R}}=2\right)$
are available with
\begin{align*}
{\rm v}_{1}^{\mathcal{I}\left({\rm R}\right)} & =\frac{1}{\sqrt{3}}\left[1,-1,1\right]^{\top},\hspace{1em}{\rm v}_{2}^{\mathcal{I}\left({\rm R}\right)}=\left[0,0,1\right]^{\top}
\end{align*}
while the two body-frame vectors are defined as in \eqref{eq:SE3STCH_Vect_R}
${\rm v}_{i}^{\mathcal{B}\left({\rm R}\right)}=R^{\top}{\rm v}_{i}^{\mathcal{I}\left({\rm R}\right)}+{\rm b}_{i}^{\mathcal{B}\left({\rm R}\right)}+\omega_{i}^{\mathcal{B}\left({\rm R}\right)}$
for $i=1,2$ such that ${\rm b}_{1}^{\mathcal{B}\left({\rm R}\right)}=0.1\left[-1,1,0.5\right]^{\top}$
and ${\rm b}_{2}^{\mathcal{B}\left({\rm R}\right)}=0.1\left[0,0,1\right]^{\top}$.
In addition, $\omega_{1}^{\mathcal{B}\left({\rm R}\right)}$ and $\omega_{2}^{\mathcal{B}\left({\rm R}\right)}$
are Gaussian noise vectors with zero mean and ${\rm STD}=0.1$. The
third vector is obtained using ${\rm v}_{3}^{\mathcal{I}\left({\rm R}\right)}={\rm v}_{1}^{\mathcal{I}\left({\rm R}\right)}\times{\rm v}_{2}^{\mathcal{I}\left({\rm R}\right)}$
and ${\rm v}_{3}^{\mathcal{B}\left({\rm R}\right)}={\rm v}_{1}^{\mathcal{B}\left({\rm R}\right)}\times{\rm v}_{2}^{\mathcal{B}\left({\rm R}\right)}$.
This step is followed by the normalization of ${\rm v}_{i}^{\mathcal{B}\left({\rm R}\right)}$
and ${\rm v}_{i}^{\mathcal{I}\left({\rm R}\right)}$ to $\upsilon_{i}^{\mathcal{B}\left({\rm R}\right)}$
and $\upsilon_{i}^{\mathcal{I}\left({\rm R}\right)}$, respectively,
for $i=1,2,3$ as given in \eqref{eq:SE3STCH_Vector_norm}. Thus,
Assumption \ref{Assum:SE3STCH_1} holds. For the semi-direct pose
filter with prescribed performance, $R_{y}$ is obtained by SVD \cite{markley1988attitude},
or for simplicity visit the Appendix in \cite{hashim2018SO3Stochastic}
with $\tilde{R}=\hat{R}R_{y}^{\top}$. The total simulation time is
30 seconds.

Initial attitude error is set to be considerably large. Initial attitude
estimate is given by $\hat{R}\left(0\right)=\mathcal{R}_{\alpha}\left(\alpha,u/||u||\right)$
according to angle-axis parameterization as in \eqref{eq:SE3STCH_att_ang}
with $\alpha=175\left({\rm deg}\right)$ and $u$= $\left[3,10,8\right]^{\top}$.
It is worth noting that the value of $||\tilde{R}||_{I}\approx0.999$
is fairly close to the unstable equilibria $\left(+1\right)$ and
the initial position is $\hat{P}\left(0\right)=\left[4,-3,5\right]^{\top}$.
In brief, we have{\small{}
	\[
	\boldsymbol{T}\left(0\right)=\mathbf{I}_{4},\hspace{1em}\hat{\boldsymbol{T}}\left(0\right)=\left[\begin{array}{cccc}
	-0.8923 & 0.2932 & 0.3432 & 4\\
	0.3992 & 0.1577 & 0.9032 & -3\\
	0.2107 & 0.9430 & -0.2577 & 5\\
	0 & 0 & 0 & 1
	\end{array}\right]
	\]
}The design parameters of the proposed filters are chosen as $\gamma=1$,
$k_{w}=5$, $\bar{\delta}=\underline{\delta}=\left[1.3,5,4,6\right]^{\top}$,
$\xi^{0}=\left[1.3,5,-4,6\right]^{\top}$, $\xi^{\infty}=\left[0.07,0.3,0.3,0.3\right]^{\top}$,
and $\ell=\left[4,4,4,4\right]^{\top}$. The initial bias estimates
are $\hat{b}_{\Omega}\left(0\right)=\left[0,0,0\right]^{\top}$ and
$\hat{b}_{V}\left(0\right)=\left[0,0,0\right]^{\top}$.

Color notation used in the plots is: black center-lines and green
solid-lines refer to the true values, red illustrates the performance
of the nonlinear semi-direct pose filter (S-DIR) on $\mathbb{SE}\left(3\right)$
proposed in Subsection \ref{subsec:SE3PPF_Passive-Filter}, and blue
demonstrates the performance of the direct filter (DIR) on $\mathbb{SE}\left(3\right)$
presented in Subsection \ref{subsec:SE3PPF_Explicit-Filter}. Also,
magenta depicts a measured value while orange and black dashed lines
refer to the prescribed performance response.

\begin{figure}[h!]
	\centering{}\includegraphics[scale=0.27]{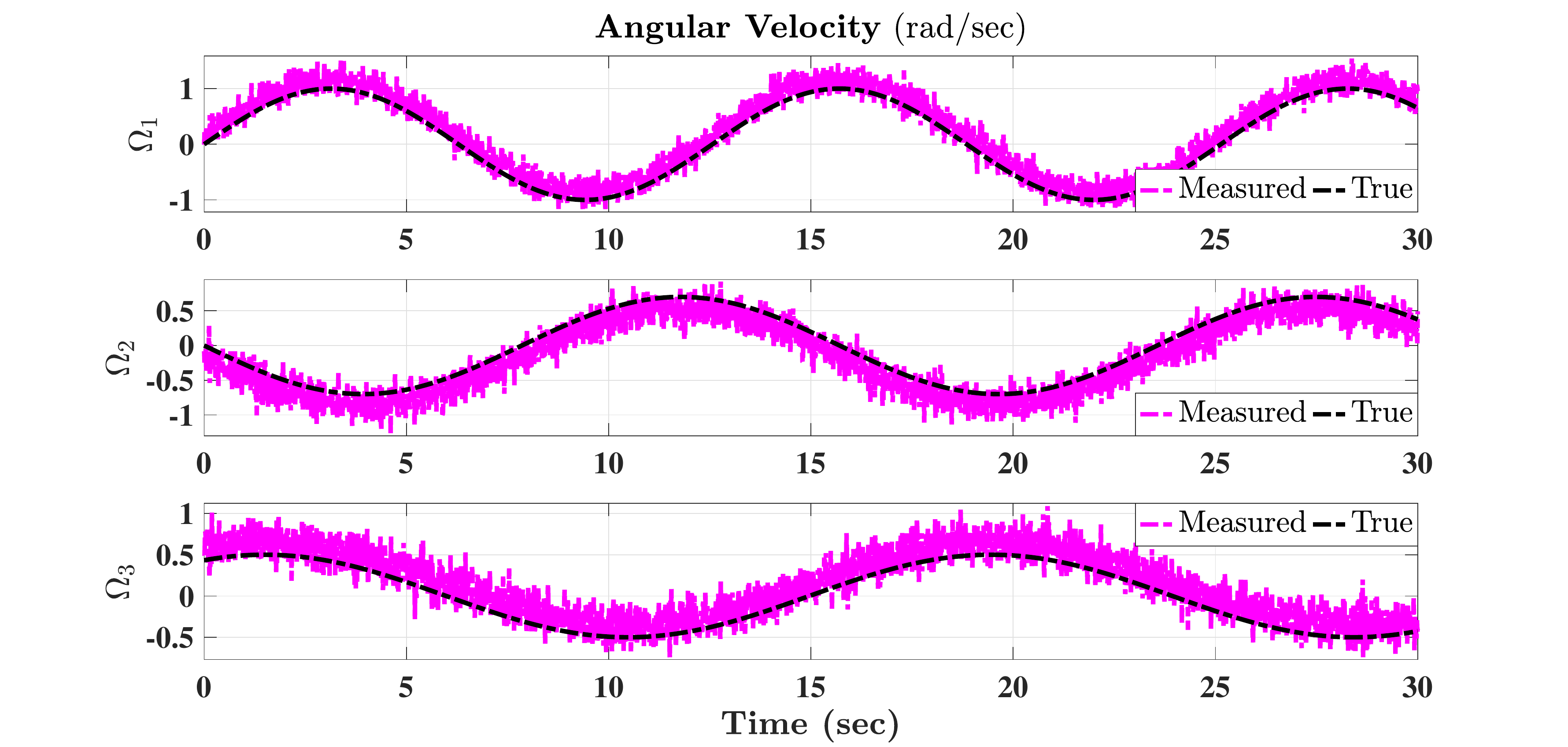}\caption{Measured and true values of angular velocities.}
	\label{fig:SE3PPF_Simulation1} 
\end{figure}

\begin{figure}[h!]
	\centering{}\includegraphics[scale=0.27]{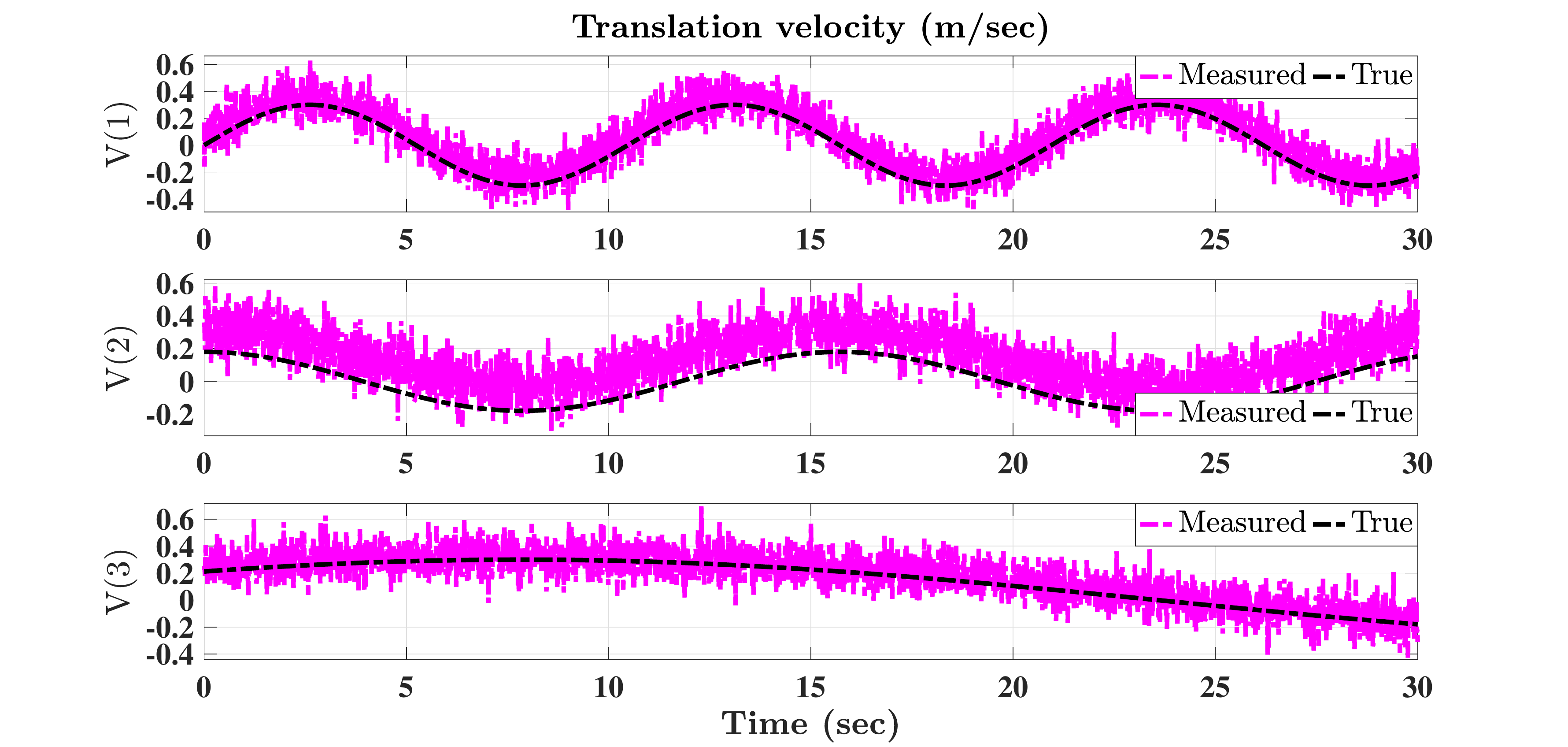}\caption{Measured and true values of translational velocities.}
	\label{fig:SE3PPF_Simulation2} 
\end{figure}

\begin{figure}[h!]
	\centering{}\includegraphics[scale=0.27]{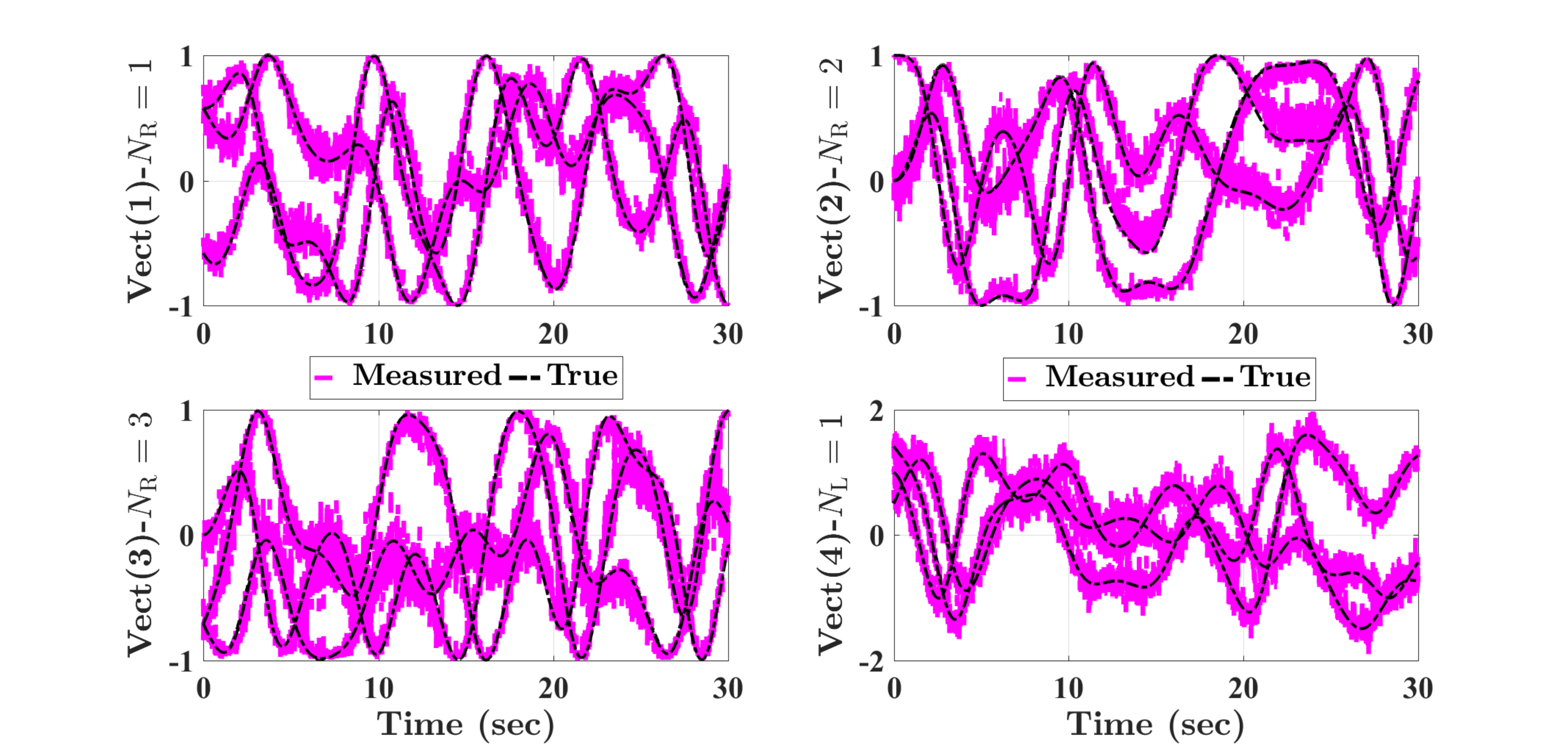}\caption{True and measured body-frame vectorial measurements.}
	\label{fig:SE3PPF_Simulation3} 
\end{figure}

\begin{figure}[h!]
	\centering{}\includegraphics[scale=0.27]{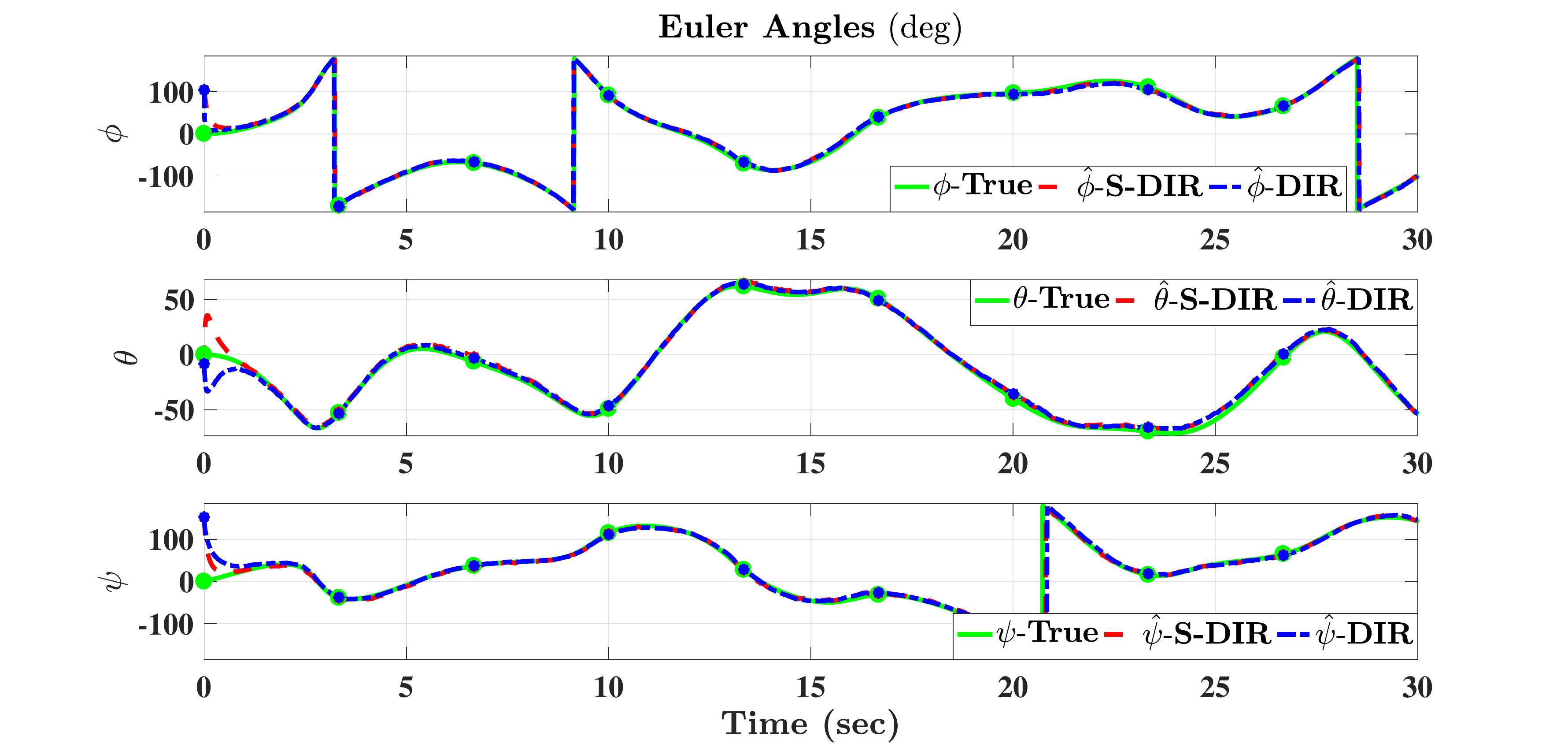}\caption{True and estimated Euler angles of the rigid-body.}
	\label{fig:SE3PPF_Simulation4} 
\end{figure}

\begin{figure}[h!]
	\centering{}\includegraphics[scale=0.27]{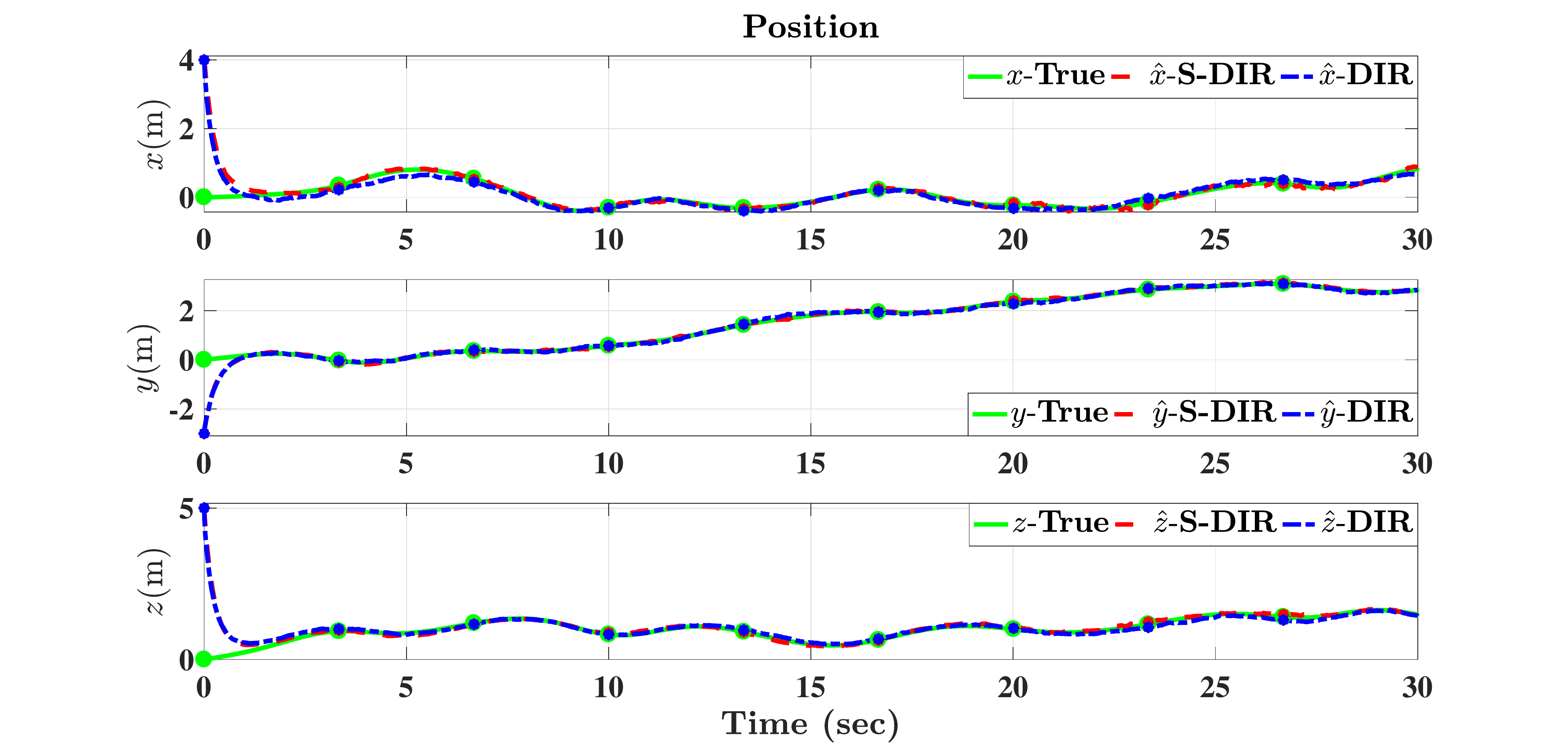}\caption{True and estimated rigid-body positions in 3D space.}
	\label{fig:SE3PPF_Simulation5} 
\end{figure}

\begin{figure*}[h!]
	\centering{}\includegraphics[scale=0.49]{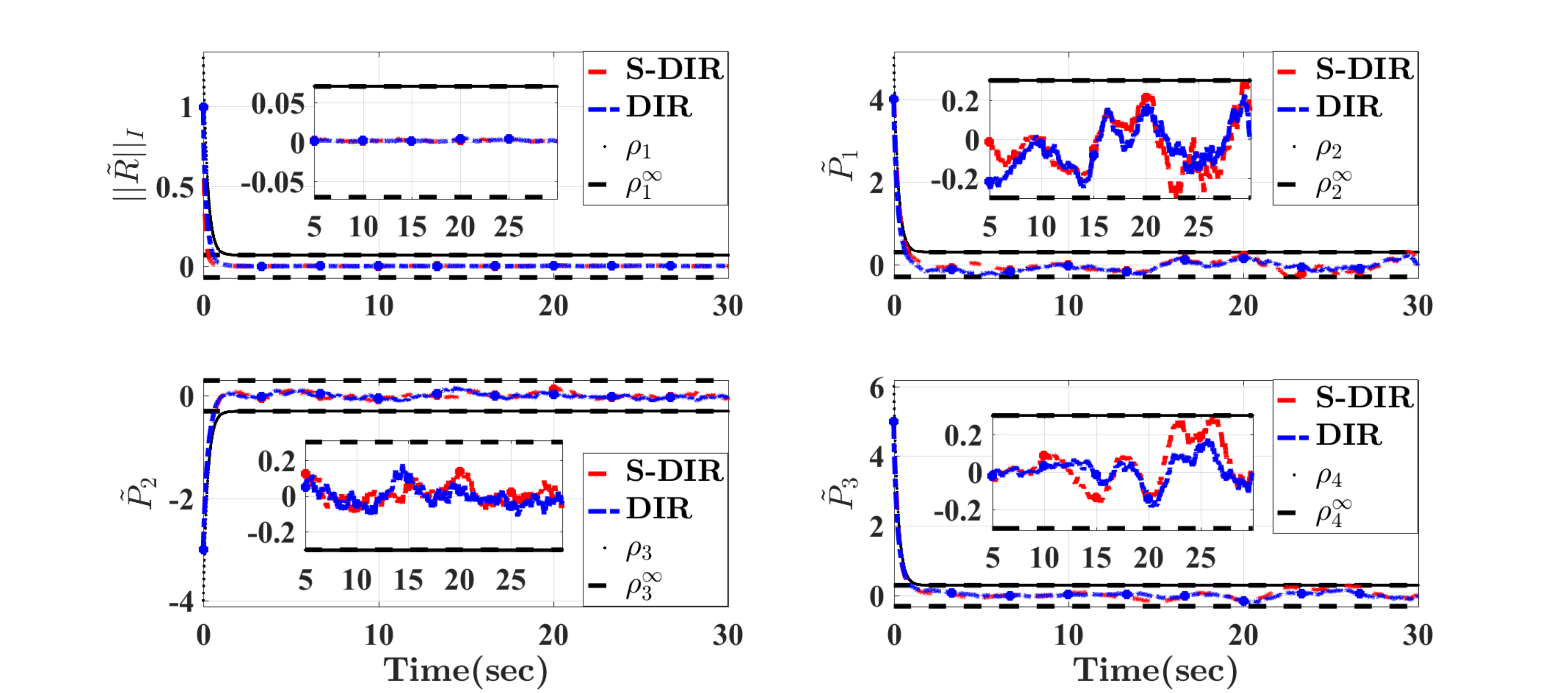}\caption{Systematic convergence of the error trajectories within the prescribed
		performance boundaries.}
	\label{fig:SE3PPF_Simulation6} 
\end{figure*}

\begin{figure}[h!]
	\centering{}\includegraphics[scale=0.27]{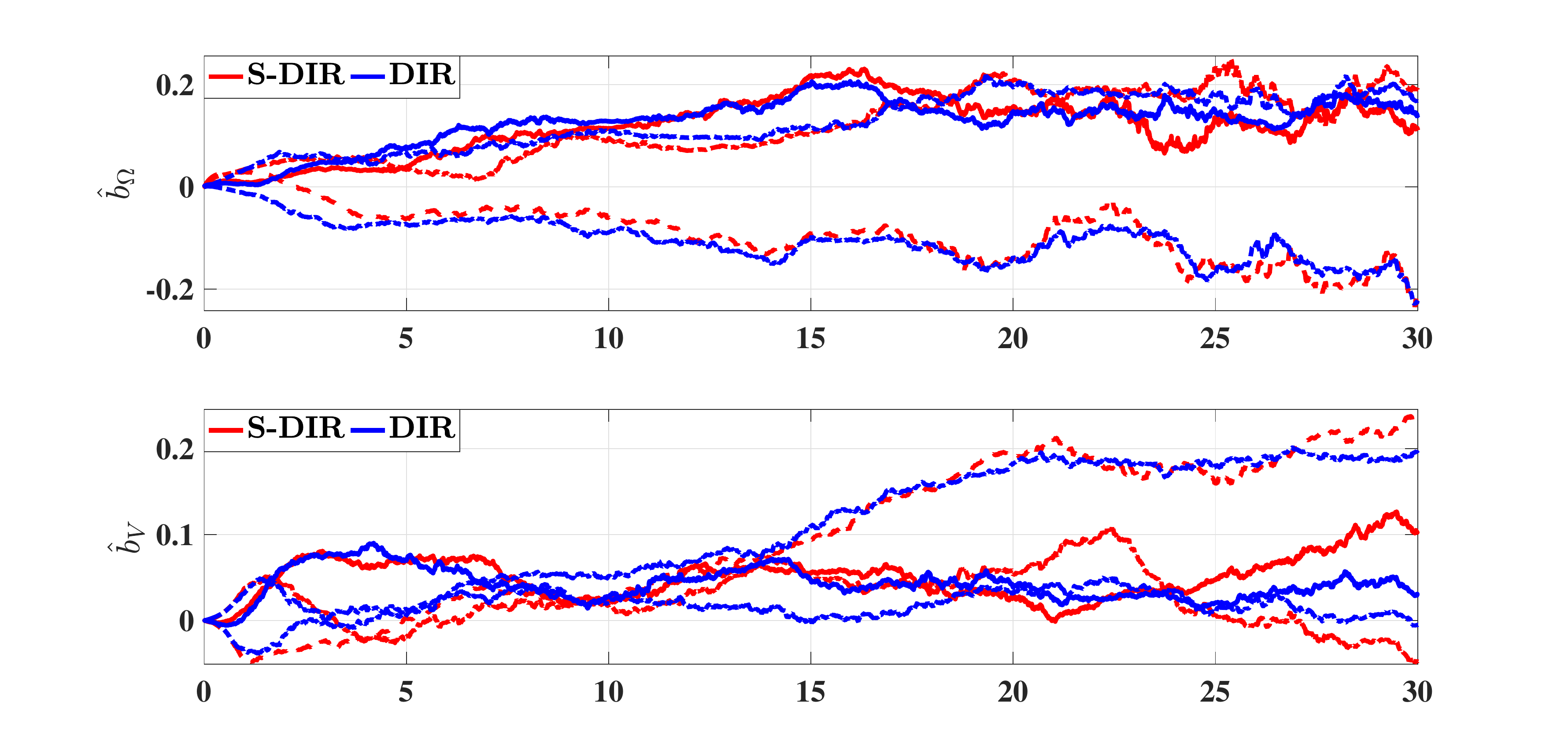}\caption{The estimated bias of the proposed filters.}
	\label{fig:SE3PPF_Simulation7} 
\end{figure}

Fig. \ref{fig:SE3PPF_Simulation1}, \ref{fig:SE3PPF_Simulation2}
and \ref{fig:SE3PPF_Simulation3} depict high values of noise and
bias components attached to velocity and body-frame vector measurements
plotted against the true values. Fig. \ref{fig:SE3PPF_Simulation4}
and \ref{fig:SE3PPF_Simulation5} show the output performance of the
proposed filters described in terms of Euler angles $\left(\phi,\theta,\psi\right)$
and the true position in 3D space, respectively. Fig. \ref{fig:SE3PPF_Simulation4}
and \ref{fig:SE3PPF_Simulation5} present remarkable tracking performance
with fast convergence to the true Euler angles and $xyz$-positions
3D space. The systematic and smooth convergence of the error vector
$\boldsymbol{e}$ is depicted in Fig. \ref{fig:SE3PPF_Simulation6}.
It can be clearly observed how $||\tilde{R}||_{I}$ in Fig. \ref{fig:SE3PPF_Simulation6}
started very near to the unstable equilibria while $\tilde{P}_{1}$,
$\tilde{P}_{2}$, and $\tilde{P}_{3}$ started remarkably far from
the origin within the predefined large set and decayed smoothly and
systematically to the predefined small set guided by the dynamic boundaries
of the PPF such that $\tilde{R}=\hat{R}R^{\top}$ and $\tilde{P}=\hat{P}-\tilde{R}P$.
Finally, the estimated bias $\hat{b}$ is bounded as depicted in Fig.
\ref{fig:SE3PPF_Simulation7}.

The simulation results establish the strong filtering capability of
the two proposed pose filters and their robustness against uncertain
measurements and large initialized errors making them perfectly fit
for the measurements obtained from low quality sensors such as IMU.
The two filters conform to the dynamic constraints imposed by the
user referring guaranteed prescribed performance measures in transient
as well as steady-state performance. The pose filters previously proposed
in the literature \cite{rehbinder2003pose,baldwin2007complementary,baldwin2009nonlinear,hua2011observer,vasconcelos2010nonlinear}
lack this remarkable quality. Semi-direct pose filter with prescribed
performance demands pose reconstruction, in this case attitude has
been extracted using SVD \cite{markley1988attitude,hashim2018SO3Stochastic}.
This adds complexity, and therefore the semi-direct pose filter requires
more computational power in comparison with the direct pose filter
with prescribed performance. Nevertheless, the two proposed pose filters
are robust and demonstrate impressive convergence capabilities. %

\section{Conclusion\label{sec:SO3PPF_Conclusion}}

Two nonlinear pose filters evolved directly on $\mathbb{SE}\left(3\right)$
with prescribed performance characteristics have been considered.
Pose error has been defined in terms of position error and normalized
Euclidean distance error, and the innovation term has been selected
to guarantee predefined measures of transient and steady-state performance.
As a result, the proposed filters exhibit superior convergence properties
with transient error being bounded by a predefined dynamically decreasing
constrained function and steady-state error being less than a predefined
lower bound. %
{} The proposed pose filters are deterministic and the stability analysis
ensure boundedness of all closed loop signals with asymptotic convergence
of the homogeneous transformation matrix to the origin. Simulation
results established the strong ability of the proposed filters to
impose the predefined constraints on the pose error considering large
initial pose error and high level of uncertainties in the measurements.

\section*{Acknowledgment}

The authors would like to thank \textbf{Maria Shaposhnikova} for proofreading the article.

\section*{Appendix A \label{sec:SO3STCH_EXPL_AppendixA} }
\begin{center}
	\textbf{Proof of Lemma }\textbf{\large{}\ref{Lemm:SE3PPF_1}}{\large\par}
	\par\end{center}

{} Let $R\in\mathbb{SO}\left(3\right)$ be the attitude of a rigid-body
in 3D space. The attitude could be extracted for a given Rodriguez
parameters vector $\rho\in\mathbb{R}^{3}$. The mapping from Rodriguez
vector to $\mathbb{SO}\left(3\right)$ is defined by $\mathcal{R}_{\rho}:\mathbb{R}^{3}\rightarrow\mathbb{SO}\left(3\right)$
\cite{shuster1993survey}
\begin{align}
\mathcal{R}_{\rho}\left(\rho\right)= & \frac{1}{1+||\rho||^{2}}\left(\left(1-||\rho||^{2}\right)\mathbf{I}_{3}+2\rho\rho^{\top}+2\left[\rho\right]_{\times}\right)\label{eq:SE3PPF_SO3_Rodr}
\end{align}
With direct substitution of \eqref{eq:SE3PPF_SO3_Rodr} in \eqref{eq:SE3STCH_Ecul_Dist}
one easily obtains \cite{hashim2018SO3Stochastic}
\begin{equation}
||R||_{I}=\frac{||\rho||^{2}}{1+||\rho||^{2}}\label{eq:SE3PPF_TR2}
\end{equation}
Additionally, for $\mathcal{R}_{\rho}=\mathcal{R}_{\rho}\left(\rho\right)$
the anti-symmetric projection operator of the attitude in \eqref{eq:SE3PPF_SO3_Rodr}
is equivalent to
\begin{align}
\boldsymbol{\mathcal{P}}_{a}\left(R\right)=\frac{1}{2}\left(\mathcal{R}_{\rho}-\mathcal{R}_{\rho}^{\top}\right)= & 2\frac{1}{1+||\rho||^{2}}\left[\rho\right]_{\times}\label{eq:SE3PPF_Pa}
\end{align}
Thus, the vex operator of \eqref{eq:SE3PPF_Pa} becomes
\begin{equation}
\mathbf{vex}\left(\boldsymbol{\mathcal{P}}_{a}\left(R\right)\right)=2\frac{\rho}{1+||\rho||^{2}}\label{eq:SE3PPF_VEX_Pa}
\end{equation}
From the result in \eqref{eq:SE3PPF_TR2} one can obtain
\begin{equation}
\left(1-||R||_{I}\right)||R||_{I}=\frac{||\rho||^{2}}{\left(1+||\rho||^{2}\right)^{2}}\label{eq:SE3PPF_append1}
\end{equation}
and from \eqref{eq:SE3PPF_VEX_Pa} it is easily shown that
\begin{equation}
||\mathbf{vex}\left(\boldsymbol{\mathcal{P}}_{a}\left(R\right)\right)||^{2}=4\frac{||\rho||^{2}}{\left(1+||\rho||^{2}\right)^{2}}\label{eq:SE3PPF_append2}
\end{equation}
Therefore, \eqref{eq:SE3PPF_append1} and \eqref{eq:SE3PPF_append2}
prove \eqref{eq:SE3PPF_lemm1_1} in Lemma \ref{Lemm:SE3PPF_1}. From
Section \ref{subsec:SE3PPF_Explicit-Filter} $\sum_{i=1}^{N_{{\rm R}}}k_{i}^{{\rm R}}=3$
which indicates that ${\rm Tr}\left\{ \mathbf{M}_{{\rm R}}\right\} =3$.
Recall that the normalized Euclidean distance of $R\mathbf{M}_{{\rm R}}$
is $\left\Vert R\mathbf{M}_{{\rm R}}\right\Vert _{I}=\frac{1}{4}{\rm Tr}\left\{ \left(\mathbf{I}_{3}-R\right)\mathbf{M}_{{\rm R}}\right\} $.
From the angle-axis parameterization in \eqref{eq:SE3STCH_att_ang},
one finds
\begin{align}
\left\Vert R\mathbf{M}_{{\rm R}}\right\Vert _{I} & =\frac{1}{4}{\rm Tr}\left\{ -\left(\sin(\theta)\left[u\right]_{\times}+\left(1-\cos(\theta)\right)\left[u\right]_{\times}^{2}\right)\mathbf{M}_{{\rm R}}\right\} \nonumber \\
& =-\frac{1}{4}{\rm Tr}\left\{ \left(1-\cos(\theta)\right)\left[u\right]_{\times}^{2}\mathbf{M}_{{\rm R}}\right\} \label{eq:SO3STCH_EXPL_append3}
\end{align}
where ${\rm Tr}\left\{ \left[u\right]_{\times}\mathbf{M}_{{\rm R}}\right\} =0$
as in identity \eqref{eq:SO3PPF_Identity6}. One has \cite{murray1994mathematical}
\begin{equation}
\left\Vert R\right\Vert _{I}=\frac{1}{4}{\rm Tr}\left\{ \mathbf{I}_{3}-R\right\} ={\rm sin}^{2}\left(\theta/2\right)\label{eq:SE3PPF_EXPL_append4}
\end{equation}
The Rodriguez vector can be expressed in terms of angle-axis parameterization
as \cite{shuster1993survey}
\begin{equation}
u={\rm cot}\left(\theta/2\right)\rho\label{eq:SE3PPF_EXPL_append5}
\end{equation}
From identity \eqref{eq:SO3PPF_Identity3} %
{} and \eqref{eq:SE3PPF_EXPL_append5}, the expression in \eqref{eq:SO3STCH_EXPL_append3}
becomes
\begin{align*}
\left\Vert R\mathbf{M}_{{\rm R}}\right\Vert _{I} & =\frac{1}{2}\left\Vert R\right\Vert _{I}u^{\top}\bar{\mathbf{M}}_{{\rm R}}u=\frac{1}{2}\left\Vert R\right\Vert _{I}{\rm cot}^{2}\left(\frac{\theta}{2}\right)\rho^{\top}\bar{\mathbf{M}}_{{\rm R}}\rho
\end{align*}
Also, from \eqref{eq:SE3PPF_EXPL_append4}, ${\rm cos}^{2}\left(\frac{\theta}{2}\right)=1-\left\Vert R\right\Vert _{I}$
which implies that
\[
{\rm tan}^{2}\left(\frac{\theta}{2}\right)=\frac{\left\Vert R\right\Vert _{I}}{1-\left\Vert R\right\Vert _{I}}
\]
Accordingly, the normalized Euclidean distance of $R\mathbf{M}_{{\rm R}}$
could be formulated in the sense of Rodriguez vector 
\begin{align}
\left\Vert R\mathbf{M}_{{\rm R}}\right\Vert _{I} & =\frac{1}{2}\left(1-\left\Vert R\right\Vert _{I}\right)\rho^{\top}\bar{\mathbf{M}}_{{\rm R}}\rho=\frac{1}{2}\frac{\rho^{\top}\bar{\mathbf{M}}_{{\rm R}}\rho}{1+\left\Vert \rho\right\Vert ^{2}}\label{eq:SE3PPF_EXPL_append_MBR_I}
\end{align}
The anti-symmetric projection operator of $R\mathbf{M}_{{\rm R}}$
can be defined in terms of Rodriquez vector using identity \eqref{eq:SO3PPF_Identity1}
and \eqref{eq:SO3PPF_Identity4} by
\begin{align*}
\boldsymbol{\mathcal{P}}_{a}\left(R\mathbf{M}_{{\rm R}}\right) & =\frac{\rho\rho^{\top}\mathbf{M}_{{\rm R}}-\mathbf{M}_{{\rm R}}\rho\rho^{\top}+\mathbf{M}_{{\rm R}}\left[\rho\right]_{\times}+\left[\rho\right]_{\times}\mathbf{M}_{{\rm R}}}{1+\left\Vert \rho\right\Vert ^{2}}\\
& =\frac{\left[\left({\rm Tr}\left\{ \mathbf{M}_{{\rm R}}\right\} \mathbf{I}_{3}-\mathbf{M}_{{\rm R}}-\left[\rho\right]_{\times}\mathbf{M}_{{\rm R}}\right)\rho\right]_{\times}}{1+\left\Vert \rho\right\Vert ^{2}}
\end{align*}
Thereby, the vex operator of the above expression is
\begin{align}
\mathcal{\mathbf{vex}}\left(\boldsymbol{\mathcal{P}}_{a}\left(R\mathbf{M}_{{\rm R}}\right)\right) & =\frac{\left(\mathbf{I}_{3}+\left[\rho\right]_{\times}\right)}{1+\left\Vert \rho\right\Vert ^{2}}\bar{\mathbf{M}}_{{\rm R}}\rho\label{eq:SE3PPF_EXPL_append_MBR_VEX}
\end{align}
Hence, the 2-norm of \eqref{eq:SE3PPF_EXPL_append_MBR_VEX} is equivalent
to
\begin{align*}
\left\Vert \mathcal{\mathbf{vex}}\left(\boldsymbol{\mathcal{P}}_{a}\left(R\mathbf{M}_{{\rm R}}\right)\right)\right\Vert ^{2} & =\frac{\rho^{\top}\bar{\mathbf{M}}_{{\rm R}}\left(\mathbf{I}_{3}-\left[\rho\right]_{\times}^{2}\right)\bar{\mathbf{M}}_{{\rm R}}\rho}{\left(1+\left\Vert \rho\right\Vert ^{2}\right)^{2}}
\end{align*}
From the identity in \eqref{eq:SO3PPF_Identity3}, $\left[\rho\right]_{\times}^{2}=-||\rho||^{2}\mathbf{I}_{3}+\rho\rho^{\top}$
such that
\begin{align}
\left\Vert \mathcal{\mathbf{vex}}\left(\boldsymbol{\mathcal{P}}_{a}\left(R\mathbf{M}_{{\rm R}}\right)\right)\right\Vert ^{2} & =\frac{\rho^{\top}\bar{\mathbf{M}}_{{\rm R}}\left(\mathbf{I}_{3}-\left[\rho\right]_{\times}^{2}\right)\bar{\mathbf{M}}_{{\rm R}}\rho}{\left(1+\left\Vert \rho\right\Vert ^{2}\right)^{2}}\nonumber \\
& =\frac{\rho^{\top}\left(\bar{\mathbf{M}}_{{\rm R}}\right)^{2}\rho}{1+\left\Vert \rho\right\Vert ^{2}}-\frac{\left(\rho^{\top}\bar{\mathbf{M}}_{{\rm R}}\rho\right)^{2}}{\left(1+\left\Vert \rho\right\Vert ^{2}\right)^{2}}\nonumber \\
& \geq\underline{\lambda}\left(1-\frac{\left\Vert \rho\right\Vert ^{2}}{1+\left\Vert \rho\right\Vert ^{2}}\right)\frac{\rho^{\top}\bar{\mathbf{M}}_{{\rm R}}\rho}{1+||\rho||^{2}}\nonumber \\
& \geq2\underline{\lambda}\left(1-\left\Vert R\right\Vert _{I}\right)\left\Vert R\mathbf{M}_{{\rm R}}\right\Vert _{I}\label{eq:SE3PPF_EXPL_append_VEX_MI2}
\end{align}
where $\underline{\lambda}=\underline{\lambda}\left(\bar{\mathbf{M}}_{{\rm R}}\right)$
is the minimum singular value of $\bar{\mathbf{M}}_{{\rm R}}$ and
$\left\Vert R\right\Vert _{I}=\left\Vert \rho\right\Vert ^{2}/\left(1+\left\Vert \rho\right\Vert ^{2}\right)$
as in \eqref{eq:SE3PPF_TR2}. One can find
\begin{align}
1-\left\Vert R\right\Vert _{I} & =\frac{1}{4}\left(1+{\rm Tr}\left\{ R\mathbf{M}_{{\rm R}}\mathbf{M}_{{\rm R}}^{-1}\right\} \right)\label{eq:SE3PPF_EXPL_append_rho2}
\end{align}
Hence, from \eqref{eq:SE3PPF_EXPL_append_VEX_MI2} and \eqref{eq:SE3PPF_EXPL_append_rho2}
the following inequality holds
\begin{align*}
\left\Vert \mathcal{\mathbf{vex}}\left(\boldsymbol{\mathcal{P}}_{a}\left(R\mathbf{M}_{{\rm R}}\right)\right)\right\Vert ^{2} & \geq\frac{\underline{\lambda}}{2}\left(1+{\rm Tr}\left\{ R\mathbf{M}_{{\rm R}}\mathbf{M}_{{\rm R}}^{-1}\right\} \right)\left\Vert R\mathbf{M}_{{\rm R}}\right\Vert _{I}
\end{align*}
This validates \eqref{eq:SE3PPF_lemm1_2} and completes the proof
of Lemma \ref{Lemm:SE3PPF_1}.

\section*{Appendix B\label{sec:SO3_PPF_STCH_AppendixB} }
\begin{center}
	\textbf{\large{}{}{}{}{}{}{}{}{}{}{}{}Quaternion Representation}{\large{}{}{}
	} 
	\par\end{center}

\noindent Define $Q=[q_{0},q^{\top}]^{\top}\in\mathbb{S}^{3}$ as
a unit-quaternion with $q_{0}\in\mathbb{R}$ and $q\in\mathbb{R}^{3}$
such that $\mathbb{S}^{3}=\{\left.Q\in\mathbb{R}^{4}\right|||Q||=\sqrt{q_{0}^{2}+q^{\top}q}=1\}$.
$Q^{-1}=[\begin{array}{cc}
q_{0} & -q^{\top}\end{array}]^{\top}\in\mathbb{S}^{3}$ denotes the inverse of $Q$. Define $\odot$ as a quaternion product
where the quaternion multiplication of $Q_{1}=[\begin{array}{cc}
q_{01} & q_{1}^{\top}\end{array}]^{\top}\in\mathbb{S}^{3}$ and $Q_{2}=[\begin{array}{cc}
q_{02} & q_{2}^{\top}\end{array}]^{\top}\in\mathbb{S}^{3}$ is $Q_{1}\odot Q_{2}=[q_{01}q_{02}-q_{1}^{\top}q_{2},q_{01}q_{2}+q_{02}q_{1}+[q_{1}]_{\times}q_{2}]$.
The mapping from unit-quaternion ($\mathbb{S}^{3}$) to $\mathbb{SO}\left(3\right)$
is described by $\mathcal{R}_{Q}:\mathbb{S}^{3}\rightarrow\mathbb{SO}\left(3\right)$
\begin{align}
\mathcal{R}_{Q} & =(q_{0}^{2}-||q||^{2})\mathbf{I}_{3}+2qq^{\top}+2q_{0}\left[q\right]_{\times}\in\mathbb{SO}\left(3\right)\label{eq:NAV_Append_SO3}
\end{align}
The quaternion identity is described by $Q_{{\rm I}}=[1,0,0,0]^{\top}$
with $\mathcal{R}_{Q_{{\rm I}}}=\mathbf{I}_{3}$. Define the estimate
of $Q=[q_{0},q^{\top}]^{\top}\in\mathbb{S}^{3}$ as $\hat{Q}=[\hat{q}_{0},\hat{q}^{\top}]^{\top}\in\mathbb{S}^{3}$
with $\mathcal{R}_{\hat{Q}}=(\hat{q}_{0}^{2}-||\hat{q}||^{2})\mathbf{I}_{3}+2\hat{q}\hat{q}^{\top}+2\hat{q}_{0}\left[\hat{q}\right]_{\times}$,
see the map in \eqref{eq:NAV_Append_SO3}. For any $x\in\mathbb{R}^{3}$
and $Q\in\mathbb{S}^{3}$, define the map
\begin{align*}
\overline{x} & =[0,x^{\top}]^{\top}\in\mathbb{R}^{4}\\
\overline{\mathbf{Y}(Q^{-1},x)} & =\left[\begin{array}{c}
0\\
\mathbf{Y}(Q^{-1},x)
\end{array}\right]=Q^{-1}\odot\left[\begin{array}{c}
0\\
x
\end{array}\right]\odot Q\\
\overline{\mathbf{Y}(Q,x)} & =\left[\begin{array}{c}
0\\
\mathbf{Y}(Q,x)
\end{array}\right]=Q\odot\left[\begin{array}{c}
0\\
x
\end{array}\right]\odot Q^{-1}
\end{align*}
The equivalent quaternion representation and complete implementation
steps of the filter in \eqref{eq:SE3PPF_Rest_dot_Ty}, \eqref{eq:SE3PPF_Pest_dot_Ty},
\eqref{eq:SE3PPF_b1est_dot_Ty}, \eqref{eq:SE3PPF_b2est_dot_Ty},
\eqref{eq:SE3PPF_W1est_dot_Ty}, and \eqref{eq:SE3PPF_W2est_dot_Ty}
is:
\[
\begin{cases}
\upsilon_{i}^{\mathcal{B}} & =\mathbf{Y}(Q^{-1},\upsilon_{i}^{\mathcal{I}})\\
Q_{y} & :\text{Reconstructed by QUEST algorithm}\\
\tilde{Q} & =[\tilde{q}_{0},\tilde{q}^{\top}]^{\top}=\hat{Q}\odot Q_{y}^{-1}\\
||\tilde{R}||_{I} & =1-\tilde{q}_{0}^{2}\\
P_{y} & =\frac{1}{\sum_{i=1}^{N_{{\rm L}}}k_{i}^{{\rm L}}}\sum_{i=1}^{N_{{\rm L}}}s_{i}^{{\rm L}}\left({\rm v}_{i}^{\mathcal{I}\left({\rm L}\right)}-\mathbf{Y}\left(Q_{y},{\rm v}_{i}^{\mathcal{B}\left({\rm L}\right)}\right)\right)\\
\tilde{P} & =\hat{P}-\mathbf{Y}\left(\tilde{Q},P_{y}\right)\\
\Gamma & =\Omega_{m}-\hat{b}-W\\
\dot{\hat{Q}} & =\frac{1}{2}\left[\begin{array}{cc}
0 & -\Gamma^{\top}\\
\Gamma & -\left[\Gamma\right]_{\times}
\end{array}\right]\hat{Q}\\
\dot{\hat{P}} & =\mathbf{Y}\left(\hat{Q},V_{m}-\hat{b}_{V}-W_{V}\right)\\
\dot{\hat{b}}_{\Omega} & =\gamma\boldsymbol{\Psi}_{R}\mathcal{E}_{R}\tilde{q}_{0}\mathbf{Y}\left(\hat{Q}^{-1},\tilde{q}\right)\\
& \hspace{1em}+\gamma\left[\mathbf{Y}\left(\hat{Q}^{-1},\tilde{P}-\hat{P}\right)\right]_{\times}\mathbf{Y}\left(\hat{Q}^{-1},\boldsymbol{\Psi}_{P}\mathcal{E}_{P}\right)\\
\dot{\hat{b}}_{V} & =\gamma\mathbf{Y}\left(\hat{Q}^{-1},\boldsymbol{\Psi}_{P}\mathcal{E}_{P}\right)\\
W_{\Omega} & =4\frac{k_{w}\boldsymbol{\Psi}_{R}\mathcal{E}_{R}-\boldsymbol{\Lambda}_{R}/4}{\tilde{q}_{0}}\tilde{q}\\
W_{V} & =\mathbf{Y}\left(\hat{Q}^{-1},k_{w}\boldsymbol{\Psi}_{P}\mathcal{E}_{P}+\left[\tilde{P}-\hat{P}\right]_{\times}W_{\Omega}-\boldsymbol{\Lambda}_{P}\tilde{P}\right)
\end{cases}
\]

\noindent The equivalent quaternion representation and complete implementation
steps of the filter in \eqref{eq:SE3PPF_Rest_dot_T_VM}, \eqref{eq:SE3PPF_Pest_dot_T_VM},
\eqref{eq:SE3PPF_b1est_dot_T_VM}, \eqref{eq:SE3PPF_b2est_dot_T_VM},
\eqref{eq:SE3PPF_W1est_dot_T_VM}, and \eqref{eq:SE3PPF_W2est_dot_T_VM}
is:
\[
\begin{cases}
\left[\begin{array}{c}
0\\
\upsilon_{i}^{\mathcal{B}}
\end{array}\right] & =\left[\begin{array}{c}
0\\
\mathbf{Y}(Q^{-1},\upsilon_{i}^{\mathcal{I}})
\end{array}\right]=Q^{-1}\odot\left[\begin{array}{c}
0\\
\upsilon_{i}^{\mathcal{I}}
\end{array}\right]\odot Q\\
\left[\begin{array}{c}
0\\
\hat{\upsilon}_{i}^{\mathcal{B}}
\end{array}\right] & =\left[\begin{array}{c}
0\\
\mathbf{Y}(\hat{Q}^{-1},\upsilon_{i}^{\mathcal{I}})
\end{array}\right]=\hat{Q}^{-1}\odot\left[\begin{array}{c}
0\\
\upsilon_{i}^{\mathcal{I}}
\end{array}\right]\odot\hat{Q}\\
\boldsymbol{\Upsilon} & =\mathbf{Y}\left(\hat{Q},\sum_{i=1}^{N_{{\rm R}}}\left(\frac{s_{i}^{{\rm R}}}{2}\hat{\upsilon}_{i}^{\mathcal{B}\left({\rm R}\right)}\times\upsilon_{i}^{\mathcal{B}\left({\rm R}\right)}\right)\right)\\
||\tilde{R}\mathbf{M}_{{\rm R}}||_{I} & =\frac{1}{4}\sum_{i=1}^{N_{{\rm R}}}\left(1-\left(\hat{\upsilon}_{i}^{\mathcal{B}\left({\rm R}\right)}\right)^{\top}\upsilon_{i}^{\mathcal{B}\left({\rm R}\right)}\right)\\
M_{1} & =\sum_{i=1}^{N_{{\rm R}}}s_{i}^{{\rm R}}\upsilon_{i}^{\mathcal{B}\left({\rm R}\right)}\left(\upsilon_{i}^{\mathcal{I}\left({\rm R}\right)}\right)^{\top}\\
M_{2} & =\left(\sum_{i=1}^{N_{{\rm R}}}s_{i}^{{\rm R}}\hat{\upsilon}_{i}^{\mathcal{B}\left({\rm R}\right)}\left(\upsilon_{i}^{\mathcal{I}\left({\rm R}\right)}\right)^{\top}\right)^{-1}\\
\tilde{P} & =\hat{P}+\frac{1}{\mathbf{m}_{{\rm c}}}\left(\mathbf{Y}\left(\hat{Q},\mathbf{k}_{{\rm v}}\right)-M_{1}M_{2}\mathbf{m}_{{\rm v}}\right)\\
\Gamma & =\Omega_{m}-\hat{b}-W\\
\dot{\hat{Q}} & =\frac{1}{2}\left[\begin{array}{cc}
0 & -\Gamma^{\top}\\
\Gamma & -\left[\Gamma\right]_{\times}
\end{array}\right]\hat{Q}\\
\dot{\hat{P}} & =\mathbf{Y}\left(\hat{Q},V_{m}-\hat{b}_{V}-W_{V}\right)\\
\dot{\hat{b}}_{\Omega} & =\frac{\gamma}{2}\boldsymbol{\Psi}_{R}\mathcal{E}_{R}\mathbf{Y}\left(\hat{Q}^{-1},\boldsymbol{\Upsilon}\right)\\
& \hspace{1em}+\gamma\left[\mathbf{Y}\left(\hat{Q}^{-1},\tilde{P}-\hat{P}\right)\right]_{\times}\mathbf{Y}\left(\hat{Q}^{-1},\boldsymbol{\Psi}_{P}\mathcal{E}_{P}\right)\\
\dot{\hat{b}}_{V} & =\gamma\mathbf{Y}\left(\hat{Q}^{-1},\boldsymbol{\Psi}_{P}\mathcal{E}_{P}\right)\\
W_{\Omega} & =\frac{4}{\underline{\lambda}}\frac{k_{w}\boldsymbol{\Psi}_{R}\mathcal{E}_{R}-\boldsymbol{\Lambda}_{R}}{1+{\rm Tr\{M_{1}M_{2}\}}}\boldsymbol{\Upsilon}\\
W_{V} & =\mathbf{Y}\left(\hat{Q}^{-1},k_{w}\boldsymbol{\Psi}_{P}\mathcal{E}_{P}+\left[\tilde{P}-\hat{P}\right]_{\times}W_{\Omega}-\boldsymbol{\Lambda}_{P}\tilde{P}\right)
\end{cases}
\]

\bibliographystyle{IEEEtran}
\bibliography{bib_PPF_SE3}

\begin{thebibliography}{10}
\providecommand{\url}[1]{#1}
\csname url@samestyle\endcsname
\providecommand{\newblock}{\relax}
\providecommand{\bibinfo}[2]{#2}
\providecommand{\BIBentrySTDinterwordspacing}{\spaceskip=0pt\relax}
\providecommand{\BIBentryALTinterwordstretchfactor}{4}
\providecommand{\BIBentryALTinterwordspacing}{\spaceskip=\fontdimen2\font plus
\BIBentryALTinterwordstretchfactor\fontdimen3\font minus
  \fontdimen4\font\relax}
\providecommand{\BIBforeignlanguage}[2]{{%
\expandafter\ifx\csname l@#1\endcsname\relax
\typeout{** WARNING: IEEEtran.bst: No hyphenation pattern has been}%
\typeout{** loaded for the language `#1'. Using the pattern for}%
\typeout{** the default language instead.}%
\else
\language=\csname l@#1\endcsname
\fi
#2}}
\providecommand{\BIBdecl}{\relax}
\BIBdecl

\bibitem{shuster1981three}
M.~D. Shuster and S.~D. Oh, ``Three-axis attitude determination from vector
  observations,'' \emph{Journal of Guidance, Control, and Dynamics}, vol.~4,
  pp. 70--77, 1981.

\bibitem{markley1988attitude}
F.~L. Markley, ``Attitude determination using vector observations and the
  singular value decomposition,'' \emph{Journal of the Astronautical Sciences},
  vol.~36, no.~3, pp. 245--258, 1988.

\bibitem{hashim2018SO3Stochastic}
H.~A. Hashim, L.~J. Brown, and K.~McIsaac, ``Nonlinear stochastic attitude
  filters on the special orthogonal group 3: Ito and stratonovich,'' \emph{IEEE
  Transactions on Systems, Man, and Cybernetics: Systems}, pp. 1--13, 2018.

\bibitem{mahony2008nonlinear}
R.~Mahony, T.~Hamel, and J.-M. Pflimlin, ``Nonlinear complementary filters on
  the special orthogonal group,'' \emph{IEEE Transactions on Automatic
  Control}, vol.~53, no.~5, pp. 1203--1218, 2008.

\bibitem{mohamed2019filter}
H.~A.~H. Mohamed, ``Nonlinear attitude and pose filters with superior
  convergence properties,'' \emph{Ph.D, University of Western Ontario}, 2019.

\bibitem{choukroun2006novel}
D.~Choukroun, I.~Y. Bar-Itzhack, and Y.~Oshman, ``Novel quaternion kalman
  filter,'' \emph{IEEE Transactions on Aerospace and Electronic Systems},
  vol.~42, no.~1, pp. 174--190, 2006.

\bibitem{lefferts1982kalman}
E.~J. Lefferts, F.~L. Markley, and M.~D. Shuster, ``Kalman filtering for
  spacecraft attitude estimation,'' \emph{Journal of Guidance, Control, and
  Dynamics}, vol.~5, no.~5, pp. 417--429, 1982.

\bibitem{markley2003attitude}
F.~L. Markley, ``Attitude error representations for kalman filtering,''
  \emph{Journal of guidance, control, and dynamics}, vol.~26, no.~2, pp.
  311--317, 2003.

\bibitem{hashim2018Conf1}
H.~A. Hashim, L.~J. Brown, and K.~McIsaac, ``Nonlinear explicit stochastic
  attitude filter on {SO}(3),'' in \emph{Proceedings of the 57th {IEEE}
  conference on Decision and Control ({CDC})}.\hskip 1em plus 0.5em minus
  0.4em\relax IEEE, 2018, pp. 1210 --1216.

\bibitem{grip2012attitude}
H.~F. Grip, T.~I. Fossen, T.~A. Johansen, and A.~Saberi, ``Attitude estimation
  using biased gyro and vector measurements with time-varying reference
  vectors,'' \emph{IEEE Transactions on Automatic Control}, vol.~57, no.~5, pp.
  1332--1338, 2012.

\bibitem{liu2018complementary}
S.~Q. Liu and R.~Zhu, ``A complementary filter based on multi-sample rotation
  vector for attitude estimation,'' \emph{IEEE Sensors Journal}, 2018.

\bibitem{rehbinder2003pose}
H.~Rehbinder and B.~K. Ghosh, ``Pose estimation using line-based dynamic vision
  and inertial sensors,'' \emph{IEEE Transactions on Automatic Control},
  vol.~48, no.~2, pp. 186--199, 2003.

\bibitem{baldwin2007complementary}
G.~Baldwin, R.~Mahony, J.~Trumpf, T.~Hamel, and T.~Cheviron, ``Complementary
  filter design on the special euclidean group se (3),'' in \emph{Control
  Conference (ECC), 2007 European}.\hskip 1em plus 0.5em minus 0.4em\relax
  IEEE, 2007, pp. 3763--3770.

\bibitem{hashim2018SE3Stochastic}
H.~A. Hashim, L.~J. Brown, and K.~McIsaac, ``Nonlinear stochastic position and
  attitude filter on the special euclidean group 3,'' \emph{Journal of the
  Franklin Institute}, vol. 356, no.~7, pp. 4144--4173, 2018.

\bibitem{hashim2019Conf1}
H.~A. Hashim, L.~J. Brown, and K.~McIsaac, ``Guaranteed performance of
  nonlinear pose filter on {SE}(3),'' in \emph{Proceedings of the American
  Control Conference ({ACC})}, 2019, pp. 1--6.

\bibitem{baldwin2009nonlinear}
G.~Baldwin, R.~Mahony, and J.~Trumpf, ``A nonlinear observer for 6 dof pose
  estimation from inertial and bearing measurements,'' in \emph{Robotics and
  Automation, 2009. ICRA'09. IEEE International Conference on}.\hskip 1em plus
  0.5em minus 0.4em\relax IEEE, 2009, pp. 2237--2242.

\bibitem{hua2011observer}
M.-D. Hua, T.~Hamel, R.~Mahony, and J.~Trumpf, ``Gradient-like observer design
  on the special euclidean group se (3) with system outputs on the real
  projective space,'' in \emph{Decision and Control (CDC), 2015 IEEE 54th
  Annual Conference on}.\hskip 1em plus 0.5em minus 0.4em\relax IEEE, 2015, pp.
  2139--2145.

\bibitem{vasconcelos2010nonlinear}
J.~F. Vasconcelos, R.~Cunha, C.~Silvestre, and P.~Oliveira, ``A nonlinear
  position and attitude observer on se (3) using landmark measurements,''
  \emph{Systems \& Control Letters}, vol.~59, no.~3, pp. 155--166, 2010.

\bibitem{dominguez2017simultaneous}
S.~Dominguez, ``Simultaneous recognition and relative pose estimation of 3d
  objects using 4d orthonormal moments,'' \emph{Sensors}, vol.~17, no.~9, p.
  2122, 2017.

\bibitem{hua2017riccati}
M.-D. Hua and G.~Allibert, ``Riccati observer design for pose, linear velocity
  and gravity direction estimation using landmark position and imu
  measurements,'' in \emph{2018 IEEE Conference on Control Technology and
  Applications}, 2018.

\bibitem{tanaka2017practical}
M.~Tanaka, K.~Tanaka, and H.~O. Wang, ``Practical model construction and stable
  control of an unmanned aerial vehicle with a parafoil-type wing,'' \emph{IEEE
  Transactions on Systems, Man, and Cybernetics: Systems}, 2017.

\bibitem{santoso2018robust}
F.~Santoso, M.~A. Garratt, S.~G. Anavatti, and I.~Petersen, ``Robust hybrid
  nonlinear control systems for the dynamics of a quadcopter drone,''
  \emph{IEEE Transactions on Systems, Man, and Cybernetics: Systems}, no.~99,
  pp. 1--13, 2018.

\bibitem{sun2018disturbance}
L.~Sun, W.~Huo, and Z.~Jiao, ``Disturbance observer-based robust relative pose
  control for spacecraft rendezvous and proximity operations under input
  saturation,'' \emph{IEEE Transactions on Aerospace and Electronic Systems},
  2018.

\bibitem{lee2018geometric}
T.~Lee, ``Geometric control of quadrotor uavs transporting a cable-suspended
  rigid body,'' \emph{IEEE Transactions on Control Systems Technology},
  vol.~26, no.~1, pp. 255--264, 2018.

\bibitem{bechlioulis2008robust}
C.~P. Bechlioulis and G.~A. Rovithakis, ``Robust adaptive control of feedback
  linearizable mimo nonlinear systems with prescribed performance,'' \emph{IEEE
  Transactions on Automatic Control}, vol.~53, no.~9, pp. 2090--2099, 2008.

\bibitem{wang2017dynamic}
M.~Wang and A.~Yang, ``Dynamic learning from adaptive neural control of robot
  manipulators with prescribed performance,'' \emph{IEEE Transactions on
  Systems, Man, and Cybernetics: Systems}, vol.~47, no.~8, pp. 2244--2255,
  2017.

\bibitem{yang2018prescribed}
Y.~Yang, J.~Tan, and D.~Yue, ``Prescribed performance tracking control of a
  class of uncertain pure-feedback nonlinear systems with input saturation,''
  \emph{IEEE Transactions on Systems, Man, and Cybernetics: Systems}, 2018.

\bibitem{na2014adaptive}
J.~Na, Q.~Chen, X.~Ren, and Y.~Guo, ``Adaptive prescribed performance motion
  control of servo mechanisms with friction compensation,'' \emph{IEEE
  Transactions on Industrial Electronics}, vol.~61, no.~1, pp. 486--494, 2014.

\bibitem{hashim2017neuro}
H.~A. Hashim, S.~El-Ferik, and F.~L. Lewis, ``Neuro-adaptive cooperative
  tracking control with prescribed performance of unknown higher-order
  nonlinear multi-agent systems,'' \emph{International Journal of Control},
  vol.~92, no.~2, pp. 445--460, 2019.

\bibitem{hashim2017adaptive}
H.~A. Hashim, S.~El-Ferik, and F.~L. Lewis, ``Adaptive synchronisation of
  unknown nonlinear networked systems with prescribed performance,''
  \emph{International Journal of Systems Science}, vol.~48, no.~4, pp.
  885--898, 2017.

\bibitem{hashim2019SO3Det}
H.~A. Hashim, L.~J. Brown, and K.~McIsaac, ``Guaranteed performance of
  nonlinear attitude filters on the special orthogonal group {SO}(3),''
  \emph{IEEE Access}, vol.~7, no.~1, pp. 3731--3745, 2019.

\bibitem{bullo2004geometric}
F.~Bullo and A.~D. Lewis, \emph{Geometric control of mechanical systems:
  modeling, analysis, and design for simple mechanical control systems}.\hskip
  1em plus 0.5em minus 0.4em\relax Springer Science \& Business Media, 2004,
  vol.~49.

\bibitem{shuster1993survey}
M.~D. Shuster, ``A survey of attitude representations,'' \emph{Navigation},
  vol.~8, no.~9, pp. 439--517, 1993.

\bibitem{murray1994mathematical}
R.~M. Murray, Z.~Li, S.~S. Sastry, and S.~S. Sastry, \emph{A mathematical
  introduction to robotic manipulation}.\hskip 1em plus 0.5em minus 0.4em\relax
  CRC press, 1994.

\end{thebibliography}

\vspace{20pt}

\begin{IEEEbiography}
	[{\includegraphics[scale=0.11,clip,keepaspectratio]{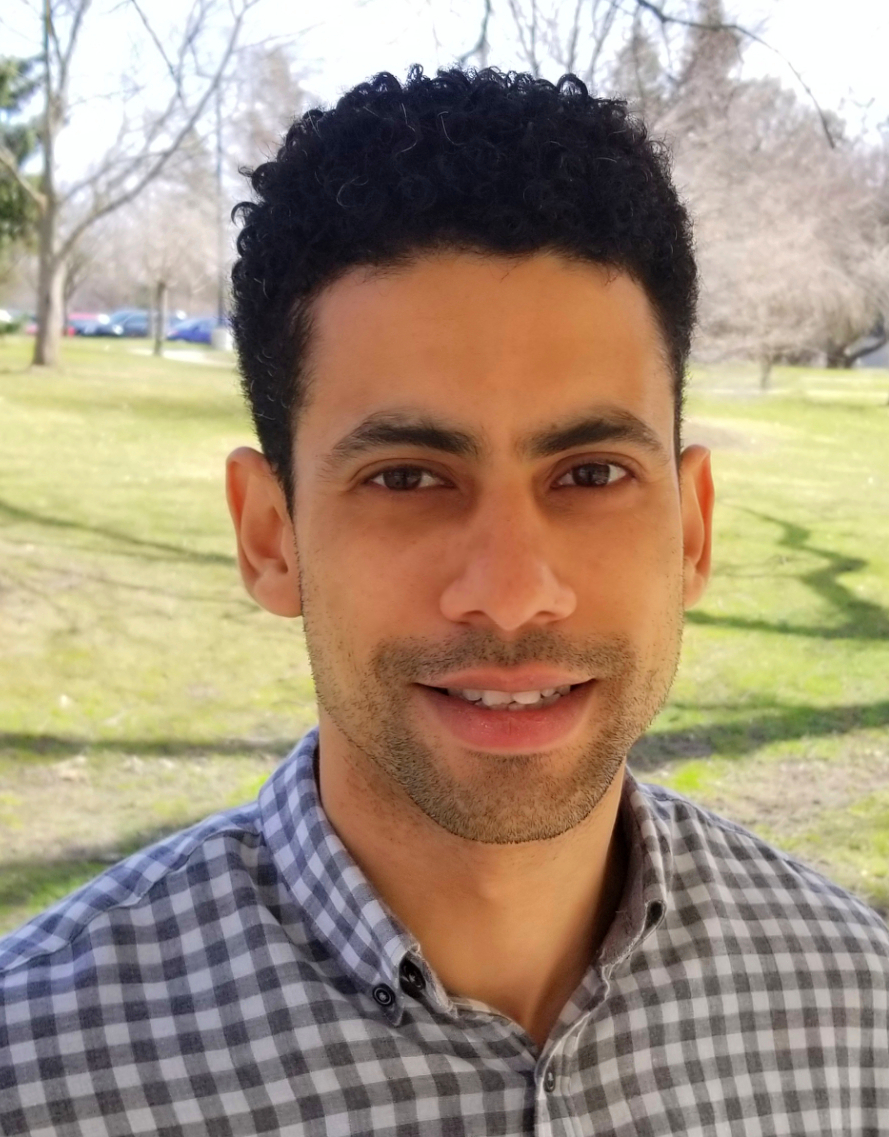}}]{Hashim A. Hashim}
	(S’18) is a Ph.D. candidate in Robotics and Control, Department of Electrical and Computer Engineering at Western University, Ontario, Canada. He received the Bachelor's degree in Mechatronics, Department of Mechanical Engineering from Helwan University, Cairo, Egypt and the M.Sc. in Systems and Control Engineering, Department of Systems Engineering from King Fahd University of Petroleum \& Minerals, Dhahran, Saudi Arabia.\\
	His current research interests include stochastic and deterministic attitude and pose filters, control of multi-agent systems, control applications and optimization techniques.
\end{IEEEbiography}

\begin{IEEEbiography}
	[{\includegraphics[scale=0.28,clip,keepaspectratio]{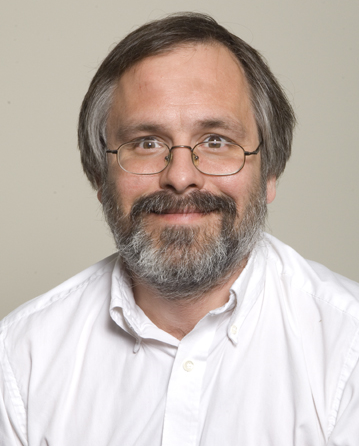}}]{Lyndon J. Brown}
	received the B.Sc. degree from the University of Waterloo, Canada in 1988 and the M.Sc. and PhD. degrees from the University of Illinois, Urbana-Champaign in 1991 and 1996, respectively. He is an associate professor in the department of electrical and computer engineering at Western University, Canada. He worked in industry for Honeywell Aerospace Canada and E.I. DuPont de Nemours.\\
	His current research includes the identification and control of predictable signals, biological control systems, welding control systems, and attitude and pose estimation.
\end{IEEEbiography}

\begin{IEEEbiography}
	[{\includegraphics[scale=0.49,clip,keepaspectratio]{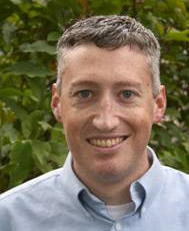}}]{Kenneth McIsaac}
	(M’99) received the B.Sc. degree from the University of Waterloo, Canada, in 1996, and the M.Sc. and Ph.D. degrees from the University of Pennsylvania, in 1998 and 2001, respectively. He is currently an Associate Professor and the Chair of Electrical and Computer Engineering with Western University, ON, Canada.\\
	His current research interests include computer vision and signal processing, mostly in the context of using machine intelligence in robotics and assistive systems, and attitude and pose estimation.
\end{IEEEbiography}

\end{document}